\documentclass{amsart}
\usepackage{amssymb, amsfonts, amsthm}
\usepackage[margin=1 in]{geometry}
\usepackage{graphicx}
\usepackage{xcolor}
\usepackage{ytableau}
\usepackage[all]{xy}
\usepackage{tikz}
\usetikzlibrary{arrows.meta}
\pgfdeclarelayer{bg}    
\pgfsetlayers{bg,main}  
\usepackage{hyperref}

\title{Richardson varieties, projected Richardson varieties and positroid varieties}

\author{David E Speyer}

\thanks{}

\theoremstyle{plain}
\newtheorem{theorem}{Theorem}[section]
\newtheorem{Theorem}[theorem]{Theorem}
\newtheorem{TD}[theorem]{Theorem/Definition}

\newtheorem{prop}[theorem]{Proposition}
\newtheorem{lemma}[theorem]{Lemma}

\newtheorem{cor}[theorem]{Corollary}

\theoremstyle{remark}
\newtheorem{definition}[theorem]{Definition}
\newtheorem{Problem}[theorem]{Problem}
\newtheorem{problem}[theorem]{Problem}
\newtheorem{prob}[theorem]{Problem}
\newtheorem{rem}[theorem]{Remark}
\newtheorem{remark}[theorem]{Remark}
\newtheorem{Remark}[theorem]{Remark}

\newtheorem{eg}[theorem]{Example}

\newcommand{\T}{\mathbf{T}}

\newenvironment{sbm}
    {\left[ \begin{smallmatrix}
    }
    { 
     \end{smallmatrix} \right]
    }

\newcommand{\newword}[1]{\textbf{\emph{#1}}}

\newcommand{\rev}[1]{{#1}^R}

\def\leftquot{\delimiter"526E30F\mathopen{}}

\newcommand{\into}{\hookrightarrow}
\newcommand{\onto}{\twoheadrightarrow}

\newcommand{\Fl}{\mathcal{F}\ell}
\newcommand{\Mat}{\mathrm{Mat}}
\newcommand{\Id}{\mathrm{Id}}
\newcommand{\GL}{\mathrm{GL}}
\newcommand{\SL}{\mathrm{SL}}
\newcommand{\Inv}{\mathrm{Inv}}
\newcommand{\BS}{\mathrm{BS}}
\newcommand{\BSo}{\mathrm{BS}^{\circ}}
\newcommand{\Pluck}{\texttt{Pluck}}
\newcommand{\SLPluck}{\texttt{SLPluck}}

\newcommand{\seq}{\text{seq}}

\DeclareMathOperator{\Ker}{Ker}

\DeclareMathOperator{\Spec}{Spec}

\DeclareMathOperator{\Proj}{Proj}
\DeclareMathOperator{\MultiProj}{MultiProj}

\DeclareMathOperator{\Span}{Span}
\DeclareMathOperator{\In}{in}
\DeclareMathOperator{\displace}{\texttt{disp}}
\newcommand{\Bound}{\text{Bound}}
\newcommand{\Edge}{\text{Edge}}
\newcommand{\Vertices}{\text{Vert}}

\newcommand{\Alt}{\bigwedge\nolimits}

\newcommand{\Xo}{\mathring{X}}
\newcommand{\Zo}{\mathring{Z}}
\newcommand{\Ro}{\mathring{R}}
\newcommand{\Pio}{\mathring{\Pi}}

\renewcommand{\AA}{\mathbb{A}}

\newcommand{\CC}{\mathbb{C}}

\newcommand{\FF}{\mathbb{F}}
\newcommand{\GG}{\mathbb{G}}

\newcommand{\PP}{\mathbb{P}}

\newcommand{\RR}{\mathbb{R}}

\newcommand{\ZZ}{\mathbb{Z}}

\newcommand{\cB}{\mathcal{B}}

\newcommand{\cD}{\mathcal{D}}

\newcommand{\cI}{\mathcal{I}}

\newcommand{\cL}{\mathcal{L}}
\newcommand{\cM}{\mathcal{M}}

\newcommand{\cO}{\mathcal{O}}

\newcommand{\cQ}{\mathcal{Q}}
\newcommand{\cR}{\mathcal{R}}
\newcommand{\cS}{\mathcal{S}}

\newcommand{\cX}{\mathcal{X}}

\newcommand{\fP}{\mathfrak{P}}
\newcommand{\fQ}{\mathfrak{Q}}

\newcommand{\tS}{\widetilde{S}}
\newcommand{\ts}{\widetilde{s}}
\newcommand{\tA}{\widetilde{A}}

\newcommand{\tcQ}{\widehat{\cQ}} 

\newcommand{\bu}{\overline{u}}
\newcommand{\bw}{\overline{w}}

\newcommand{\covers}{\gtrdot}
\newcommand{\covered}{\lessdot}

\begin{document}

\maketitle
\vspace{-0.5 in}
\tableofcontents

\section{Background}

\subsection{Symmetric groups and related combinatorics} \label{sec:CombinatorialNotation}
We write $[n]$ for $\{ 1,2,3, \ldots, n \}$ and $[a,b]$ for $\{ a,a+1, a+2, \ldots, b \}$. 
We write $\binom{[n]}{k}$ for the set of $k$-element subsets of $[n]$ and $2^{[n]}$ for the set of all subsets of $[n]$.
We put a partial order $\preceq$ on $\binom{[n]}{k}$ by defining $\{ i_1 < i_2 < \cdots < i_k \}\preceq \{ j_1< j_2< \ldots < j_k \}$ if and only if, for $1 \leq a \leq k$, we have $i_a \leq j_a$.

\begin{eg}
Here is the Hasse diagram of the $\preceq$ order on $\binom{[4]}{2}$:
\[ \xymatrix@R=0.2 in@C=0.2 in{
& 34 \ar@{-}[d] & \\
& 24 \ar@{-}[dl]\ar@{-}[dr] & \\
23 \ar@{-}[dr] && 14 \ar@{-}[dl] \\
& 13 \ar@{-}[d] & \\
& 12 & \\
} \]
\end{eg}

We write $S_n$ for the group of permutations of $[n]$ and write the action of $S_n$ on $[n]$ on the left.
In particular, for $w \in S_n$, we have an element $w[k]$ of $\binom{[n]}{k}$, and we have an increasing chain $w[1] \subset w[2] \subset \cdots \subset w[n-1]$ of subsets of $[n]$.
We write $e$ for the identity of $S_n$ and $w_0$ for the element $w_0(j) = n+1-j$.

We write $(i j)$ for the permutation in $S_n$ which switches $i$ and $j$ and fixes the elements of $[n] \setminus \{ i,j \}$.
We put $s_i = (i \ i+1)$. 

When we write a permutation in one line notation as $z_1 z_2 \cdots z_n$, we mean the permutation $j \mapsto z_j$.
For example, $s_1 s_2$ is $231$ because $s_1(s_2(1)) = 2$, $s_1(s_2(2)) = 3$ and $s_1(s_2(3))=1$.

For $w \in S_n$, the \newword{length} of $w$ is the smallest $a$ such that we can write $w=s_{i_1} s_{i_2} \cdots s_{i_a}$.
We denote the length of $w$ as $\ell(w)$.
A word for $w$ of length $\ell(w)$ is called \newword{reduced}.
We say that $(i,j)$ is an \newword{inversion} of $w$ if $1 \leq i <j \leq n$ and $w(i) > w(j)$.
The set of inversions of $w$ is denoted $\Inv(w)$ and we have $\ell(w) = \#\Inv(w)$.

We equip $S_n$ with the partial order known as \newword{Bruhat order} or \newword{strong order}:
\begin{TD}
Let $u$ and $v \in S_n$. The following are equivalent:
\begin{enumerate}
\item For all $0 < i,j < n$, we have $\#([i] \cap u[j]) \geq \#([i] \cap v[j])$. \label{RankCondition}
\item For all $0 < i < n$, we have $u[i] \preceq v[i]$.
\item There is a reduced word $s_{j_1} s_{j_2} \cdots s_{j_{\ell}}$ for $v$ and a subword $s_{j_{a_1}} s_{j_{a_2}} \cdots s_{j_{a_m}}$ with product $u$.
\item For every reduced word $s_{j_1} s_{j_2} \cdots s_{j_{\ell}}$ for $v$, there is a subword $s_{j_{a_1}} s_{j_{a_2}} \cdots s_{j_{a_m}}$ with product $u$.
\end{enumerate}
When these equivalent conditions hold, we say that $u \preceq v$.
\end{TD}
Condition~(\ref{RankCondition}) is easily seen to be equivalent to $\#([i+1,n] \cap u[j]) \leq \#([i+1,n] \cap v[j])$, to $\#([i] \cap u[j+1,n]) \leq \#([i] \cap v[j+1,n])$, and to $\#([i+1,n] \cap u[j+1, n]) \geq \#([i+1,n] \cap v[j+1, n])$.

\begin{eg} \label{eg:S3Bruhat}
Here is the Hasse diagram of  Bruhat order on $S_3$.
\[ \xymatrix@R=0.2 in@C=0.2 in{
& 321 \ar@{-}[dr] \ar@{-}[dl] & \\
231 \ar@{-}[dd] \ar@{-}[ddrr] && 312 \ar@{-}[dd] \ar@{-}[ddll] \\
&& \\
213 \ar@{-}[dr] && 132 \ar@{-}[dl] \\
& 123  .& \\
}\]
\end{eg}

The \newword{Demazure product} is the unique associative multiplication $\ast : S_n \times S_n \to S_n$ such that
\[ s_i \ast w = \begin{cases} s_i w & \ell(s_i w) = \ell(w)+1 \\ w & \ell(s_i w) = \ell(w) - 1 \end{cases} \qquad 
w \ast s_i = \begin{cases} w s_i & \ell(w s_i) = \ell(w)+1 \\ w & \ell(w s_i) = \ell(w) - 1 .\end{cases} \]

A \newword{partition} is a finite sequence of integers $(\lambda_1, \lambda_2, \ldots, \lambda_k)$ with $\lambda_1 \geq \lambda_2 \geq \cdots \geq \lambda_k > 0$. 
We write $|\lambda|$ for $\sum \lambda_i$.
We call $k$ the \newword{length} or \newword{number of rows} of the partition, and write $k = \ell(\lambda)$.
We will feel free to pad our partitions with additional zeroes at the end, so a \newword{partition with at most $n$ rows} can be written $(\lambda_1, \lambda_2, \ldots, \lambda_n)$ with $\lambda_1 \geq \lambda_2 \geq \cdots \geq \lambda_n \geq 0$.
We write $\omega_k$ for the partition $(1,1,\ldots,1,0,0,\ldots,0)$ where there are $k$ ones and $n-k$ zeroes. 
So the abelian semigroup of  \newword{partitions with at most $n$ rows} is freely generated by $\omega_1$, $\omega_2$, \dots, $\omega_n$.

\subsection{Algebraic geometry and group notation} We work over an arbitrary field $\kappa$.
We write $\AA^1$ for the affine line over $\kappa$ and $\AA^n$ for affine $n$-space.
We write $\GG_m$ for the multiplicative group scheme over $\kappa$ and $\GG_m^n$ for the $n$-fold power of $\GG_m$.
We write $\PP^{n-1}$ for the space of lines in $\AA^n$, considered as an algebraic variety; more generally, for a vector space $V$, we write $\PP(V)$ for the space of lines in $V$.

We write $\GL_n$ for the group of $n \times n$ invertible matrices.
We write $T$ for the subgroup of $\GL_n$ of invertible diagonal matrices, so $T \cong \GG_m^n$.
We write $B_+$ and $B_-$ for the subgroups of invertible upper and lower triangular matrices respectively, and write $N_+$ and $N_-$ for the subgroups of $B_+$ and $B_-$ where the diagonal elements are $1$.
So $N_{\pm}$ is normal in $B_{\pm}$ and we have $B_{\pm} = T N_{\pm} = N_{\pm} T$.

\begin{Remark}
Many papers shorten $B_+$ to $B$; some papers shorten $B_-$ to $B$.
To avoid confusion, we give both of these subgroups their subscript.
\end{Remark}

We embed $S_n$ into $\GL_n$ by sending $w \in S_n$ to the permutation matrix which has ones in positions $(w(j), j)$ and zeroes everywhere else, so $S_n \into \GL_n$ is a map of groups.
The group $S_n$ normalizes $T$, so we have $wT = Tw$ for $w \in S_n$.
We warn the reader that the literature on Schubert polynomials often uses a map $S_n \to \GL_n$ which is an anti-homomorphism of groups, so $w$ goes to the matrix with ones in positions $(i, w(i))$; see Remark~\ref{rem:InverseWarning}.

\subsection{Grassmannians, Flag varieties, Pl\"ucker coordinates} \label{sec:BackgroundPlucker}
For a finite dimensional vector space $V$, the \newword{Grassmannian}  $G(k,V)$ is the space of $k$-dimensional subspaces of $V$, equipped with the structure of a variety in the usual way.
In particular, the Pl\"ucker embedding is a closed embedding $G(k, V) \into \PP(\Alt^k(V))$ sending the subspace with basis $v_1$, $v_2$, \dots, $v_k$ to the tensor $v_1 \wedge v_2 \wedge \cdots \wedge v_k$.

We write $G(k, n)$ for the Grassmannian of $k$-planes in $\AA^n$.
Write $e_1$, $e_2$, \dots, $e_n$ for the standard basis vectors in $\AA^n$. 
Then a basis for $\Alt^k \AA^n$ is the vectors $e_{i_1} \wedge e_{i_2} \wedge \cdots \wedge e_{i_k}$ for $1 \leq i_1 < i_2 < \cdots < i_k \leq n$. 
Let $V$ be a $k$-plane in $\AA^n$ with basis $v_1$, $v_2$, \dots, $v_k$.
Abusing notation slightly, for $1 \leq i_1 < i_2 < \cdots < i_k \leq n$, let $\Delta_{i_1 i_2 \cdots i_k}(V)$ be the coefficient of $e_{i_1} \wedge e_{i_2} \wedge \cdots \wedge e_{i_k}$in $v_1 \wedge v_2 \wedge \cdots \wedge v_k$; the $\Delta_{i_1 i_2 \cdots i_k}(V)$ are called \newword{Pl\"ucker coordinates}.

The individual Pl\"ucker coordinates are not well-defined functions on $G(k,n)$, because changing bases in $V$ multiplies all of the Pl\"ucker coordinates by a common scalar.
However, the $\binom{n}{k}$ Pl\"ucker coordinates collectively form homogeneous projective coordinates on $G(k,n)$.
We also think of $\Delta_{i_1 i_2 \cdots i_k}$ as an actual function on spaces of matrices: For an $n \times k$ matrix $M$, we take $\Delta_{i_1 i_2 \cdots i_k}(M)$ to be the minor of $M$ in the rows indexed by $\{i_1, i_2, \ldots, i_k\}$; for an $n \times n$ matrix $g$; we take $\Delta_{i_1 i_2 \cdots i_k}(g)$ to be the minor of $g$ in the rows indexed by $\{i_1, i_2, \ldots, i_k\}$ and the $k$ leftmost columns. 
The compatibility between these notations is that the Pl\"ucker coordinates $(\Delta_{i_1 \cdots i_k}(M))$ (respectively, $(\Delta_{i_1 \cdots i_k}(g))$) are the homogeneous Pl\"ucker coordinates of the $k$-plane $M \AA^k$ (respectively, $g \Span(e_1, e_2, \ldots, e_k)$).
We note the boundary cases $\Delta_{\emptyset}(g) = 1$ and $\Delta_{[n]}(g) = \det(g)$.

We pause to address an issue of signs: If we write $\Delta_{i_1 i_2 \cdots i_k}$ where the $i$-indices are not in increasing order, what do we mean by $\Delta_{i_1 i_2 \cdots i_k}(V)$?
Our conventions are the following: If two of $i_1$, $i_2$, \dots, $i_k$ are equal, then $\Delta_{i_1 i_2 \cdots i_k}(V)$ is $0$; if the $i_1$, $i_2$, \dots, $i_k$ are distinct and $\sigma$ is the permutation in $S_k$ with $i_{\sigma(1)} < i_{\sigma(2)} < \cdots < i_{\sigma(k)}$, then $\Delta_{i_1 i_2 \cdots i_k}(V) = (-1)^{\ell(\sigma)} \Delta_{i_{\sigma(1)} i_{\sigma(2)} \cdots i_{\sigma(k)}}(V)$.
If $I$ is a $k$-element subset of $[n]$, then $\Delta_I(V)$ is defined to be $\Delta_{i_1 i_2 \cdots i_k}(V)$ where $i_1 < i_2 < \cdots < i_k$ are the elements of $I$.


Let $V$ be a vector space of dimension $n$. 
The \newword{complete flag variety} $\Fl(V)$ is the reduced subvariety of $\prod_{k=1}^{n-1} G(k,V)$ corresponding to complete flags $V_1 \subset \cdots \subset V_{n-1}$ in $V$ with $\dim V_i = i$.
We write $\Fl_n$ for $\Fl(\AA^n)$. 
So, for each subset $I \subseteq [n]$ with $0 < \#(I) < n$, there is a Pl\"ucker coordinate $\Delta_I$. 
For each cardinality $k$, the Pl\"ucker coordinates $\Delta_I$ with $\#(I)=k$ are defined up to rescaling by a common scalar.

We will want to also refer to \newword{partial flag manifolds}. 
For $k_1 < k_2 < \cdots < k_p$, the partial flag manifold $\Fl(k_1, k_2, \ldots, k_p; V)$ is the space of chains $V_1 \subset V_2 \subset \cdots \subset V_p$ of subspaces of $V$ with $\dim V_i = k_i$; we put $\Fl_n(k_1, k_2, \ldots, k_p) : = \Fl(k_1, k_2, \ldots, k_p; \AA^n)$.

 The group $\GL_n$ acts transitively on $\Fl_n$ and on $G(k,n)$. 
 In the case of $\Fl_n$, the stabilizer of the flag $(\Span(e_1),\  \Span(e_1, e_2),\ \cdots,\ \Span(e_1, e_2, \ldots, e_{n-1}))$ is $B_+$, so we can identify $\Fl_n$ with $\GL_n/B_+$; we will also label points of $\Fl_n$ by cosets $g B_+$. 
 Explicitly, $g B_+$ corresponds to the flag whose $k$-dimensional subspace is the span of the leftmost $k$ columns of $g$.
 In particular, for $I \in \binom{[n]}{k}$, the Pl\"ucker coordinate $\Delta_I$ corresponds to the minor of $g$ in rows $I$ and columns $[k]$.
 
 \begin{Remark}
Various papers in the literature identify $\Fl_n$ with $B_- \leftquot \GL_n$ or $\GL_n / B_-$ instead of $\GL_n/B_+$ and, while the author is not aware of an example, there is surely a paper somewhere which identifies $\Fl_n$ with $B_+ \leftquot \GL_n$.
These are all equivalent.
Concretely, for $g \in \GL_n$, one should take the left-most columns, top-most rows, right-most columns and bottom-most rows in order to consider  $\GL_n/B_+$, $B_- \leftquot \GL_n$, $\GL_n / B_-$ or $B_+ \leftquot \GL_n$ respectively.
Pl\"ucker coordinates are then minors which are located in these rows/columns.
 \end{Remark}


\subsection{Bruhat decomposition, Schubert cells, Schubert varieties} \label{BruhatSchubert}
All the results in this section are standard and can be found in many sources, for example \cite[Chapters 9 and 10]{YoungTableaux} and \cite[Chapter 15]{MS}.

\begin{TD}
The permutation matrices $w \in S_n$ form a complete set of  representatives for the double cosets $B_{\pm} \leftquot \GL_n / B_{\pm}$, for any of the four choices of the $\pm$ signs.
The decomposition $\GL_n = \bigsqcup_{w\in S_n} B_{\pm} w B_{\pm}$ is called the \newword{Bruhat decomposition} of $\GL_n$.
\end{TD}

The Bruhat decomposition can be described explicitly using ranks of submatrices:
\begin{prop} \label{prop:GLnRankConditions}
Let $g \in \GL_n$ be an $n \times n$ matrix. The matrix $g$ lies in $B_- w B_+$ if and only if, for each $0<  i,j < n$, the upper-left $i \times j$ submatrix of $g$ has rank $\#([i] \cap w[j])$, and the matrix $g$ lies in the closure $\overline{B_- w B_+}$ if and only if the rank of this submatrix is $\leq \#([i] \cap w[j])$.
\end{prop}

\begin{remark} \label{rem:RankMatrices}
The matrix $r(w)_{ij} := \#([i] \cap w[j])$ is called the \newword{rank matrix} of $w$. 
It is also convenient to set $r_{0j} = r_{i0} = 0$.
Rank matrices were introduced by Fulton~\cite{Fulton92}. 
Explicitly, if $r_{ij}$ is a matrix of integers with rows and columns indexed by $\{ 0,1,2,\ldots,n \}$, then $r$ is a rank matrix if and only if (a) $r_{ij} \leq r_{(i+1) j} \leq r_{ij}+1$ and $r_{ij} \leq r_{i(j+1)} \leq r_{ij}+1$, (b) $r_{k0} = r_{0k} = 0$ and $r_{kn} = r_{nk} = k$ and (c) if $r_{(i+1)j} = r_{(i+1)(j+1)} = r_{i(j+1)} = r$ then $r_{ij} =r$.
We will introduce a similar notion of ``cyclic rank matrix" in Section~\ref{sec:AffinePerms}.
\end{remark}

Proposition~\ref{prop:GLnRankConditions} says that, set-theoretically, $\overline{B_- w B_+}$ is cut out of $\GL_n$ by the vanishing of certain minors of the upper-left submatrices of $g$.
This is also true scheme-theoretically, and in $\Mat_{n \times n}$ as well as $\GL_n$. 
\begin{theorem} \label{thm:FultonGenerators}
The reduced ideal of the closure $\overline{B_- w B_+}$, inside the affine space $\Mat_{n \times n}$, is generated by the $(\#([i] \cap w[j])+1) \times (\#([i] \cap w[j])+1)$ minors which are contained in the upper-left $i \times j$ submatrix. 
\end{theorem}

\begin{proof}
See \cite[Lemma 3.11]{Fulton92} or \cite[Chapter 15]{MS}.
\end{proof}

We will return to this result in Section~\ref{sec:MatrixSchubert}. 

\begin{remark}
If we want to study $B_+ w B_+$, $B_- w B_-$ or $B_+ w B_-$, then we should look at the lower-left, upper-right and lower-right submatrices of size $(n+1-i) \times j$, $i \times (n+1-j)$ and $(n+1-i) \times (n+1-j)$ respectively. We then compare their ranks to $\#([i+1,n] \cap w[j])$, to $\#([i] \cap w[j+1,n])$, and to $\#([i+1,n] \cap w[j+1, n])$ respectively.
\end{remark}

Comparing the above relations to the definition of Bruhat order, we have an explicit description of the closure $\overline{B_{\pm} w B_{\pm}}$:
\begin{prop}
If the two $\pm$ signs are the same, then $\overline{B_{\pm} v B_{\pm}} = \bigsqcup_{u \preceq v} B_{\pm} u B_{\pm}$ and $\dim B_{\pm} v B_{\pm} = \binom{n+1}{2} + \ell(v)$.
If the two $\pm$ signs are different, then $\overline{B_{\pm} v B_{\pm}} = \bigsqcup_{w \succeq v} B_{\pm} w B_{\pm}$ and $\dim B_{\pm} v B_{\pm} = n^2 - \ell(v)$.
\end{prop}

Identifying $\Fl_n$ with $\GL_n/B_+$, we can quotient the double cosets $B_+ w B_+$ and $B_- w B_+$ by $B_+$ to form the \newword{Schubert cells} $\Xo^w := (B_+ w B_+)/B_+$ and $\Xo_w := (B_- w B_+)/B_+$ and their closures   $X^w := \overline{(B_+ w B_+)}/B_+$ and $X_w := \overline{(B_- w B_+)}/B_+$, which are called \newword{Schubert varieties}.
These are called cells because they are, in fact, affine spaces:
\begin{prop}
The Schubert cell $\Xo^w$ is an affine space of dimension $\ell(w)$. Explicitly, every coset in $\Xo^w$ is uniquely of the form $(w + X) B_+$ where $w$ is the permutation matrix and $X$ is a matrix all of whose nonzero entries are in positions $(w(j), i)$ where $1 \leq i < j \leq n$ and $w(j) < w(i)$.

The Schubert cell $\Xo_w$ is an affine space of dimension $\binom{n}{2} - \ell(w)$. Explicitly, every coset in $\Xo_w$ is uniquely of the form $(w + X) B_+$ where $w$ is the permutation matrix and $X$ is a matrix all of whose nonzero entries are in positions $(w(j), i)$ where $1 \leq i < j \leq n$ and $w(i) < w(j)$.
\end{prop}

\begin{eg}
In $\Fl_3$, there are $3! = 6$ Schubert cells $\Xo^w$. In the table below, we depict the representation of an element in each cell as $(w+X)B_+$:
\[ 
\begin{bmatrix} 1&0&0 \\ 0&1&0 \\ 0&0&1 \\ \end{bmatrix} \quad
\begin{bmatrix} \ast&1&0 \\ 1&0&0 \\ 0&0&1 \\ \end{bmatrix} \quad
\begin{bmatrix} 1&0&0 \\ 0&\ast&1 \\ 0&1&0 \\ \end{bmatrix} \quad
\begin{bmatrix} \ast&1&0 \\ \ast&0&1 \\ 1&0&0 \\ \end{bmatrix} \quad
\begin{bmatrix} \ast&\ast&1 \\ 1&0&0 \\ 0&1&0 \\ \end{bmatrix} \quad
\begin{bmatrix} \ast&\ast&1 \\ \ast&1&0 \\ 1&0&0 \\ \end{bmatrix} .
\]
Note that the ranks of each lower-left submatrix are constant throughout the cell.
\end{eg}

Containments between Schubert varieties are given by the following propositions:
\begin{prop}
$X^w$ is the closure of $\Xo^w$, and $X_u$ is the closure of $\Xo_u$. Concretely, $X^w = \bigsqcup_{v \leq w} \Xo^v$ and $X_u = \bigsqcup_{v \geq u} \Xo_v$. 
\end{prop}

The logic behind the upper and lower indices is clearest if we think in terms of the permutation flags $v B_+$: We have $v B_+ \in X^w$ if and only if $v$ is \emph{below} $w$ and $v B_+ \in X_u$ if and only if $v$ is \emph{above} $u$.

There are also Schubert cells and Schubert varieties in Grassmannians, which we will usually index by $k$-element subsets of $[n]$.
For such a $k$-element subset $I$, we write $X_I$ for the subvariety of $G(k,n)$ consisting of $k$-planes $V$ where $\Delta_J(V)=0$ for $J \not\succeq I$, and we write $\Xo_I$ for the open subvariety of $X_I$ where $\Delta_I$ is nonzero.
Here $X_I$ is the \newword{Grassmannian Schubert varierty} and $\Xo_I$ is the Grassmannian Schubert cell. 
If we represent $V$ as the image of an $n \times k$ matrix $g$, then $\Xo_I$ is the subvariety where the first $j$ rows of $g$ have rank $\# ([j] \cap I)$ and $X_I$ is the subvariety where the rank is $\leq \# ([j] \cap I)$.
The relation between the Flag Schubert variety $X_u$ and the Grassmannian Schubert varieties $X_I$ is that $X_u$ is the space of flags $F_1 \subset F_2 \subset \cdots \subset F_{n-1}$ where $F_k \in X_{u[k]}$. 
Of course, there are also upper indexed Grassmannian Schubert varieties $X^K$, cut out by the equations $\Delta_J=0$ for $J \not \preceq K$, and given by imposing rank conditions on the bottom submatrices of $g$. 

\begin{remark}
It is common in the literature to index Grassmannian Schubert varieties by partitions, rather than subsets. The correspondence is that the subset $\{ i_1 < i_2 < \cdots < i_k \}$ corresponds to the partition $(i_k-k, i_{k-1} - k+1, \ldots, i_2-2, i_1-1)$. We will not use this convention much.
\end{remark}


Finally, we introduce one more way of thinking about the Bruhat decomposition $\GL_n = \bigsqcup_{w \in S_n} B_+ w B_+$: We can also think of this as saying that there are $n!$ orbits of $\GL_n$ on $\Fl_n \times \Fl_n$, with representatives given by the ordered pairs $(e B_+, w B_+)$. 
Given two flags, $E$ and $F$ in $V$, we will say that \newword{$F$ is $w$-related to $E$}, and write $E \xrightarrow{\ w\ } F$, 
if $(E, F)$ is in the orbit of $(eB_+, w B_+)$. 
Concretely, we have
\begin{lemma}
We have $E \xrightarrow{\ w\ } F$ if and only if 
\[ \dim (E_i \cap F_j) = \#([i] \cap w([j]))\]
for all $1 \leq i,j \leq n$.
\end{lemma}
We spell out two particular cases: We have $E \xrightarrow{s_i} F$ if $E_j = F_j$ for $j \neq i$ and $E_i \neq F_i$.
For $w_0$ equal to the longest element of $S_n$, we have $E \xrightarrow{w_0} F$ if $E_i \cap F_{n-i} = \{ 0 \}$ for $1 \leq i \leq n-1$.

\subsection{Richardson varieties} We are now ready to define our main objects. Let $u \preceq w$ in $S_n$.
The \newword{open Richardson variety} $\Ro_u^w$ is the intersection of Schubert cells $\Xo_u \cap \Xo^w$ within $\Fl_n$.
The \newword{closed Richardson variety} $R_u^w$ is the intersection of Schubert varieties $X_u \cap X^w$ within $\Fl_n$.

Richardson varieties were introduced in~\cite{KL2}, where the number of points on Richardson varieties over finite fields is shown to compute the $R$-polynomials from~\cite{KL1}.
The name ``Richardson variety" is in honor of Richardson's 1992 paper ``Intersections of double cosets in algebraic groups"~\cite{Richardson92}, which studies intersections of double cosets of the form $HxL \cap KxL$; Richardson notes that the case $(H,K,L) = (B_-, B_+, B_+)$ had already been studied by Deodhar~\cite{Deodhar85}. 
The term ``Richardson variety" seems to have first been used by Lakshmibai and Littelmann~\cite{LL03}.

\begin{eg}
The flag manifold $\Fl_2$ is the projective line $\PP^1$. The Richardson $R_{12}^{21}$ is the entire projective line; $R_{12}^{12}$ and $R_{21}^{21}$ are the points $[1:0]$ and $[0:1]$ on this projective line. 
\end{eg}

\begin{eg} \label{Flag3Example}
The flag manifold $\Fl_3$ is $3$-dimensional and contains four $2$-dimensional Richardson hypersurfaces: $R_{123}^{231}$,  $R_{132}^{321}$, $R_{123}^{312}$ and $R_{213}^{321}$.  The intersections between these can be visualized using the figure below: $R_{123}^{231}$ and $R_{132}^{321}$ correspond to the shaded trapezoids on the left
and  $R_{123}^{312}$ and $R_{213}^{321}$ correspond to the shaded trapezoids on the right, with the $0$- and $1$-dimensional Richardsons corresponding to the line segments and points of the figures. 
The boundaries of the two hexagons are the same. We haven't drawn the $3$-dimensional Richardson variety $R_{123}^{321}$, but we can think of it as the interior of the $2$-sphere formed by gluing these two hexagons along their boundaries.
When we discuss total positivity in Section~\ref{sec:Positivity}, this picture will become literally correct.

\centerline{
\begin{tikzpicture}
  \node [label={[xshift=0.0cm, yshift=-0.6cm]$123$}] (123) at (0,0) {};
  \node  [label={[xshift=-0.5cm, yshift=-0.5cm]$213$}]  (213) at (-1,1)  {};
  \node [label={[xshift=0.5cm, yshift=-0.5cm]$132$}] (132) at (1,1)  {};
  \node [label={[xshift=-0.5cm, yshift=0.0cm]$231$}] (231) at (-1,3) {};
  \node  [label={[xshift=0.5cm, yshift=0.0cm]$312$}] (312) at (1,3)  {};
  \node  [label={[xshift=0.0cm, yshift=0.0cm]$321$}]  (321) at (0,4) {};
  \node (R1) at (-0.25,1.5) {$R_{123}^{231}$};
    \node (R2) at (0.25,2.5) {$R_{132}^{321}$};
     \begin{pgfonlayer}{bg} 
  \filldraw[fill=gray, thick] (123.center) -- (213.center) -- (231.center) -- (132.center) -- (123.center) ;
    \filldraw[fill=gray, thick] (321.center) -- (312.center) -- (132.center) -- (231.center) -- (321.center) ;
  \end{pgfonlayer}
\end{tikzpicture} \qquad
\begin{tikzpicture}
  \node [label={[xshift=0.0cm, yshift=-0.6cm]$123$}] (123) at (0,0) {};
  \node  [label={[xshift=-0.5cm, yshift=-0.5cm]$213$}]  (213) at (-1,1)  {};
  \node [label={[xshift=0.5cm, yshift=-0.5cm]$132$}] (132) at (1,1)  {};
  \node [label={[xshift=-0.5cm, yshift=0.0cm]$231$}] (231) at (-1,3) {};
  \node  [label={[xshift=0.5cm, yshift=0.0cm]$312$}] (312) at (1,3)  {};
  \node  [label={[xshift=0.0cm, yshift=0.0cm]$321$}]  (321) at (0,4) {};
  \node (R1) at (0.25,1.5) {$R_{123}^{312}$};
    \node (R2) at (-0.25,2.5) {$R_{213}^{321}$};
     \begin{pgfonlayer}{bg} 
  \filldraw[fill=gray, thick] (123.center) -- (132.center) -- (312.center) -- (213.center) -- (123.center) ;
    \filldraw[fill=gray, thick] (321.center) -- (231.center) -- (213.center) -- (312.center) -- (321.center) ;
  \end{pgfonlayer}
\end{tikzpicture} }
\end{eg}

\begin{remark}
The permutation flag $v B_+$ is in $R_u^w$ if and only if $u \leq v \leq w$.
\end{remark}

\begin{remark}
We have $X^w = R_e^w$ and $X_u = R_u^{w_0}$, so Schubert varieties are a special case of Richardson varieties.
However, it is not true that Schubert cells are a special case of open Richardson varieties: $\Ro_e^w$ is an open subvariety of $\Xo^w$ and $\Ro_u^{w_0}$ is an open subvariety of $\Xo_u$.
\end{remark}

Richardsons form a stratification as one would expect:
\begin{prop} \label{prop:RichStrat}
We have $R_u^w = \bigsqcup_{u \leq u' \leq w' \leq w} \Ro_{u'}^{w'}$. We have $\dim R_u^w = \dim \Ro_u^w = \ell(w) - \ell(u)$.
\end{prop}

\begin{proof}
The first statement holds by definition. For the dimension of $\Ro_u^w$ and thus of $R_u^w$, see~\cite{Deodhar85}.
\end{proof}

We state here the basic facts about Richardsons as algebraic varieties:
\begin{prop}
Let $u \preceq w$ in $S_n$. The open Richardson variety $\Ro_u^w$ is a smooth irreducible affine variety of dimension $\ell(w) - \ell(u)$. 
The Richardson variety $R_u^w$ is an irreducible projective variety of dimension $\ell(w) - \ell(u)$ containing $\Ro_u^w$ as a dense open subvariety.
$R_u^w$ is normal and Cohen-Macaulay with rational singularities.
\end{prop}

\begin{proof}
We defined $\Ro_u^w$ as $\Xo_u \cap \Xo^w$; the varieties $\Xo_u \cap \Xo^w$ are both affine and $\Fl_n$ is separated, so $\Ro_u^w$ is affine.
Smoothness and irreducibility of $\Ro_u^w$ are due to~\cite{Richardson92}; see also  Corollary~\ref{OpenRichardsonSmoothIrreducible}.

See~\cite[Lemma 1]{GeomMonTheory}, for the facts that $R_u^w$ is normal,  Cohen-Macaulay and  irreducible of dimension $\ell(w) - \ell(u)$.
Since $R_u^w$ is irreducible and $\Ro_u^w$ is an open set of the same dimension, it is immediate that $\Ro_u^w$ is dense in $R_u^w$.
Finally, see~\cite[Appendix A]{KLS1} or~\cite[Theorem 1.1]{BilleyCoskun} for the fact that $R_u^w$ has rational singularities.
\end{proof}

The local geometry of $R_u^w$ near a point $v$ with $u \leq v \leq w$ can be described using a result of Allen Knutson, Alex Woo and Alex Yong.
Consider the neighborhood $v \Xo_e \cong \AA^{\binom{n}{2}}$ of $v$. Then we can identify $v \Xo_e$ with $\Xo_v \times \Xo^v \cong \AA^{\binom{n}{2} - \ell(v)} \times \AA^{\ell(v)}$ such that each stratum $\Ro_u^w$ factors as a product of strata; see~\cite{KWY} for the details. 

\subsection{Projected Richardson varieties} \label{sec:ProjRichBackground}
Richardson varieties are subvarieties of $\Fl_n$, obtained as the intersection of a Schubert and an opposite Schubert.
We will also want to study subvarities of partial flag manifolds $\Fl_n(k_1, k_2, \ldots, k_p)$.
In this context, a \newword{Richardson variety in a partial flag manifold} is defined to be the intersection of two opposite Schubert varieties.
However, Knutson, Lam and Speyer found that it was valuable to study the larger class of \newword{projected Richardson varieties}.

Fix a partial flag manifold $\Fl_n(k_1, k_2, \ldots, k_p)$ and let $\pi: \Fl_n \to \Fl_n(k_1, k_2, \ldots, k_p)$ be the projection map. 
For $u \preceq w$ in $S_n$, we define $\Pi_u^w$ to be $\pi(R_u^w)$, and we define a subvariety of $\Fl_n(k_1, k_2, \ldots, k_p)$ of the form $\Pi_u^w$ to be a \newword{projected Richardson variety}.

\begin{eg}
Let $n=3$ and consider the partial flag variety $\Fl_3(1) = \PP^2$. There are $19$ Richardson varieties in $\Fl_3$, which project to $7$ distinct projected Richardson varieties in $\PP^2$. 

The projected Richardsons of $\PP^2$ are the coordinate subspaces. Of these, $6$ are Richardsons of $\PP^2$ and one, the line $\{ \Delta_2=0 \}$, is not. We can see $\{ \Delta_2=0 \}$ as $\pi(R_{132}^{312})$. To see why this works, note that $R_{132}^{312}$ is flags $(F_1, F_2)$ where the Pl\"ucker coordinates $\Delta_{12}(F_2)$ and $\Delta_{23}(F_2)$ are $0$. The unique possible $F_2$ is therefore $\Span(e_1, e_3)$, and $F_1$ must be a line of the form $\Span(x_1 e_1 + x_3 e_3)$, so $\Delta_2(F_1)$ is zero.

Most projected Richardsons in $\PP^2$ are the image of several Richardsons. For example, the line $\{ \Delta_1 = 0 \}$ is the image of $R_{213}^{312}$, $R_{231}^{321}$ and $R_{213}^{321}$. The map $\pi$ is an isomorphism in the first two cases, and has relative fiber dimension one in the third. 
Similarly, the entire space $\PP^2$ is the image of $R_{123}^{312}$, $R_{132}^{321}$ and $R_{123}^{321}$. 
The map $\pi$ is birational in the first two cases, and has relative fiber dimension one in the third. 
\end{eg}

We now describe the combinatorial data indexing projected Richardson varieties.
This material is from~\cite[Section 3]{KLS1}.

Once again, we fix the partial flag manifold $\Fl_n(k_1, k_2, \ldots, k_p)$.
We write $W_P$ for the Young subgroup $S_{k_1} \times S_{k_2-k_1} \times \cdots \times S_{n-k_p}$.
We write $v \covers u$ and say $v$ \newword{covers} $u$ if $v \succ u$ and $\ell(v) = \ell(u)+1$.
We write $v \covers_P u$ and say $v$ $P$-\newword{covers} $u$ if $v \covers u$ and $v W_P \neq u W_P$. 
We define the \newword{$P$-Bruhat order} to be the transitive closure of the $P$-covering relation, and we denote $P$-Bruhat order by $\preceq_P$. 
In the case of the Grassmannian $\Fl_n(k)$, this partial order was studied by Bergeron and Sottille~\cite{BS1} under the name ``$k$-Bruhat order"; 

\begin{prop}
The map $\pi : R_u^w \to \Pi_u^w$ is birational if and only if $u \preceq_P w$. 
\end{prop}

\begin{eg} \label{eg:kBruhat}
We continue with the example of $\Fl_3(1) \cong \PP^2$.
Here is the Hasse diagram of the $1$-Bruhat order on $S_3$ (compare to Example:\ref{eg:S3Bruhat}).
\[ \xymatrix@R=0.2 in@C=0.2 in{
& 321  \ar@{-}[dl] & \\
231  \ar@{-}[ddrr] && 312 \ar@{-}[dd] \ar@{-}[ddll] \\
&& \\
213 \ar@{-}[dr] && 132 \\
& 123  .& \\
}\]
Note that $123 \preceq_1 312$ and $132 \preceq_1 321$, but $123 \not\preceq_1 321$, matching that $R_{123}^{312} \longrightarrow \PP^2$ and $R_{132}^{321} \longrightarrow \PP^2$ are birational but $R_{123}^{321} \longrightarrow \PP^2$ is not birational.
\end{eg}

\begin{prop}
For every projected Richardson $\Pi_u^w$, we can find $(u', w')$ such that $\Pi_u^w = \Pi_{u'}^{w'}$ and $u \preceq u' \preceq_P w' \preceq w$.
\end{prop}

So, every projected Richardson variety can be described as a birational image of a Richardson variety, in many ways. 
We now give a unique representative for each projected Richardson, and describe all other such representations.
We write $W^P$ for the set of $w$ such that $w$ is minimal in the coset $w W_P$.
When $W_P = S_k \times S_{n-k}$, these are the so-called \newword{$k$-Grassmannian permutations}, and they can be described concretely as the permutations whose only descent is in the $k$-th position.

\begin{prop}
If $w \in W^P$, then $u \preceq w$ if and only if $u \preceq_P w$. Each projected Richardson can be represented in exactly one way as $\pi(R_u^w)$ with $w \in W^P$ and $u \preceq w$.
If $\Pi = \pi(R_u^w)$ with $w \in W^P$ and $u \preceq w$ then the other birational representatives of $\Pi$ are precisely of the form $\pi(R_{ux}^{wx})$ where $x \in W_P$ and $\ell(ux) = \ell(u) + \ell(x)$.
\end{prop}

Thus, the projected Richardsons can be indexed by pairs $(u,w)$ where $w \in W^P$ and $u \preceq w$.

\begin{eg}
We continue with the example of $\Fl_3(1) \cong \PP^2$. The elements of $W^P$ are $123$, $213$ and $312$. Thus, our standard representative for the projected Richardson $\PP^2$ is $u = 123$, $w = 312$. The other representative, $(132, 321)$, is $(123 \cdot s_2,\ 312 \cdot s_2)$. 
\end{eg}

The \newword{open projected Richardson}, $\Pio_u^w$ is the open subvariety of $\Pi_u^w$ where we remove all proper sub-projected Richardsons of $\Pi_u^w$.
We have the following surprising result:
\begin{prop} \label{OpenRichardsonModel}
If $u \preceq_P w$, so that $\pi : R_u^w \longrightarrow \Pi_u^w$ is birational, then $\pi : \Ro_u^w \longrightarrow \Pio_u^w$ is an isomorphism.
\end{prop}
Thus, each open projected Richardson is isomorphic to an open Richardson, and so open projected Richardsons are smooth and affine.

\subsection{Positroid varieties}
A particular case of a partial flag variety is a Grassmannian: $G(k,n) = \Fl_n(k)$. 
In this case, we will refer to a projected Richardson variety as a \newword{positroid variety}.
Positroid varieties have many combinatorial descriptions which are not available for other projected Richardson varieties, and we will discuss them further in Sections~\ref{sec:Positroid} and~\ref{sec:Plabic}.

\section{The Pl\"ucker algebra and homogeneous coordinate rings of Richardson varieties} \label{sec:CoordinateRings}

As we described above, the flag manifold $\Fl_n$ embeds into a product of Grassmannians $\prod_{k=1}^{n-1} G(k,n)$ and hence, by the Pl\"ucker embedding, into a product of projective spaces $\prod_{k=1}^{n-1} \PP^{\binom{n}{k} -1}$.
In this section, we will talk about the homogeneous coordinate ring of $\Fl_n$ and of Richardson varieties in $\Fl_n$.

The Pl\"ucker algebra, denoted $\Pluck$, is the $\kappa$-algebra of functions on $\GL_n$ generated by the Pl\"ucker coordinates $\Delta_I(g)$ considered as functions on $\GL_n$. 
(We remind the reader that $\Delta_{\emptyset}(g) = 1$ and $\Delta_{[n]}(g) = \det g$.)
The Pl\"ucker algebra is $\ZZ^n$-graded where, for $I$ a $k$-element subset of $[n]$, we place the Pl\"ucker coordinate $\Delta_I$ in degree $\omega_k:=(1,1,\ldots,1,0,0,\ldots,0)$ where there are $k$ ones and $n-k$ zeroes.
Then $\Fl_n$ is the multiproj of the Pl\"ucker algebra for this grading.

\begin{remark}
To the reader who prefers algebraic geometry to classical invariant theory, the Pl\"ucker algebra may seem slightly ad hoc.
To help orient this reader, we discuss line bundles on $\Fl_n$.
Let $\cL(\omega_k)$ be the pullback of $\cO(1)$ along $\Fl_n \onto G(k,n) \into \PP(\Alt^k \AA^n)$.
The Picard group of $\Fl_n$ is free of rank $n-1$, with generators $\cL(\omega_1)$, $\cL(\omega_2)$, \dots, $\cL(\omega_{n-1})$.
For $(\lambda_1, \lambda_2, \ldots, \lambda_{n-1}, \lambda_n) \in \ZZ^n$, we define a line bundle $\cL(\lambda)$ on $\Fl_n$ by $\cL(\lambda) = \bigotimes_{k=1}^{n-1} \cL(\omega_k)^{\lambda_k - \lambda_{k+1}}$; we note that adding the same constant to all the $\lambda_j$ leaves the line bundle unchanged. 
For $\lambda_1 \geq \lambda_2 \geq \cdots \geq \lambda_n \geq 0$, the global sections $H^0(\Fl_n, \cL(\lambda))$ are the degree $\lambda$ part of the Pl\"ucker algebra.

To describe this in another way,  let $0 = \cS_0 \subset \cS_1 \subset \cS_2 \subset \cdots \subset \cS_n$ be the tautological subbundles on $\Fl_n$, where $\cS_n$ is trivial of rank $n$. Then $\cL(\omega_k) = (\Alt^k \cS_k)^{-1}$ and $\cL(e_j) = (\cS_j/\cS_{j-1})^{-1}$. Because $\cS_n$ is trivial, so is $\Alt^n \cS_n$, which is why $\cL(\omega_n)$ is trivial.

If we want to list the isomorphism classes of line bundles on $\Fl_n$ without duplication, we can index them by integer vectors $(\lambda_1, \lambda_2, \ldots, \lambda_{n-1}, 0)$.
Thus, the Cox ring of $\Fl_n$ is the subring of $\Pluck$ in degrees with $\lambda_n=0$. In Section~\ref{sec:GTDegenRich}, we will call this ring $\SLPluck$. 

To motivate the inclusion of nonzero values of $\lambda_n$ from an algebraic geometry standpoint, one can consider $\GL_n$-equivariant line bundles on $\Fl_n$. 
For $1 \leq k \leq n-1$, let $\GL_n$ act on $\cL(\omega_k)$ by the obvious action on $\Alt^k \AA^n$ and call this equivariant line bundle $\cL^{\GL}(\omega_k)$.
Define $\cL^{\GL}(\omega_n)$ to be the trivial line bundle on $\Fl_n$, where $g \in \GL_n$ acts by $\det g$.
The $\GL_n$-equivariant Picard group of $\Fl_n$ is free of rank $n$ with generators $\cL^{\GL}(\omega_1)$, $\cL^{\GL}(\omega_2)$, \dots, $\cL^{\GL}(\omega_{n-1})$,  $\cL^{\GL}(\omega_{n})$.
For  $(\lambda_1, \lambda_2, \ldots, \lambda_{n-1}, \lambda_n) \in \ZZ^n$, we define the equivariant line bundle $\cL^{\GL}(\lambda)$ to be $\bigotimes_{k=1}^{n-1} \cL^{\GL}(\omega_k)^{\lambda_k - \lambda_{k+1}} \otimes \cL^{\GL}(\omega_n)^{\lambda_n}$.
Then the $\GL_n$-equivariant Cox ring of $\Fl_n$ is $\bigotimes_{\lambda_1 \geq \cdots \geq \lambda_n} H^0(\Fl_n, \cL^{\GL}(\lambda)) \cong \Pluck[\Delta_{[n]}^{-1}]$.

Thus, the Pl\"ucker algebra sits between the ordinary Cox ring $\SLPluck$ and the equivariant Cox ring $\Pluck[\Delta_{[n]}^{-1}]$. 
From a perspective of algebraic geometry, it is less natural than either, but it is extremely well suited to combinatorial commutative algebra. 
\end{remark}

\begin{remark} \label{remarkSchur}
We remind the reader of the classical connection between $\GL_n$ representation theory and the theory of symmetric functions.
The character of $H^0(\Fl_n, \cL(\lambda))$ as a $T$-representation  is the \newword{Schur polynomial}, $s_{\lambda}(z_1, z_2, \ldots, z_n)$.
The multiplicity of the $\mu$-weight space of $H^0(\Fl_n, \cL(\lambda))$  is the \newword{Kostka number}, $K_{\lambda \mu}$.
\end{remark}

Our goal is to discuss the coordinate rings of Richardsons, not of the flag manifolds.
We consider the coordinate ring of $R_u^w$ to be the ring generated by the Pl\"ucker variables $\Delta_I$, modulo the relations which hold on $R_u^w$; we denote this ring by $\Pluck_u^w$.
We will make the abbreviations $\Pluck_u := \Pluck_u^{w_0}$ and $\Pluck^w = \Pluck_e^w$; these correspond to the Schubert varieties $X_u$ and $X^w$.

So $\Pluck_u^w$ is a quotient ring of $\Pluck$, and we refer to the kernel of $\kappa[\Delta_I] \longrightarrow \Pluck_u^w$ as the \newword{ideal of $R_u^w$}.
Again, this ring can be understood in terms of algebraic geometry: We have
\[ \Pluck_u^w = \bigoplus_{\lambda_1 \geq \lambda_2 \geq \cdots \geq \lambda_n} H^0(R_u^w, \cL(\lambda)). \]
Implicit in the above displayed equation is the following result:
\begin{theorem}[{\cite[Proposition 1]{GeomMonTheory}}]
For any $u \preceq w$, and any partition $\lambda$, the map $H^0(\Fl_n, \cL(\lambda)) \longrightarrow H^0(R_u^w, \cL(\lambda))$ is surjective. 
\end{theorem}

In the remainder of Section~\ref{sec:CoordinateRings}, we will discuss what is known about bases for $\Pluck_u^w$ and about flat degenerations of $\Pluck_u^w$.
We first provide a whirlwind review of the classical theory of $\Pluck$.

\subsection{Classical theory of the Pl\"ucker algebra}
As we have stated, the Pl\"ucker algebra is generated by the Pl\"ucker coordinates $\Delta_I$. 
We now describe the generating relations for the Pl\"ucker algebra, as a quotient of the polynomial ring $\kappa[\Delta_I]_{\emptyset \subsetneq I \subseteq [n]}$. 
These are quadratic relations known as \newword{Pl\"ucker relations}, which we now describe.

Let $1 \leq r \leq a \leq b \leq n$. 
Choose elements $i_1$, $i_2$, \dots, $i_{a-r}$, $j_1$, $j_2$, \dots, $j_{r-1}$, $m_1$, $m_2$, \dots, $m_{b+1}$ of $[n]$. Then we have the Pl\"ucker relation
\[
 \sum_{\sigma}   (-1)^{\ell(\sigma)} \Delta_{i_1 i_2 \cdots i_{a-r} m_{\sigma(1)} m_{\sigma(2)} \cdots m_{\sigma(r)}} \Delta_{j_1 j_2 \cdots j_{r-1} m_{\sigma(r+1)} m_{\sigma(r+2)} \cdots m_{\sigma(b+1)}}=0.
\]
 Here, the sum ranges over a set of coset representatives for $S_{b+1}/(S_{r} \times S_{b+1-r})$ and we use the sign conventions from Section~\ref{sec:BackgroundPlucker}.

\begin{theorem} \label{thm:2ndFundTheorem}
The Pl\"ucker relations generate the homogeneous ideal of $\Pluck$ as a quotient of $\kappa[\Delta_I]$ and, in fact, form a Gr\"obner basis for certain term orders.
\end{theorem}

\begin{remark} \label{GrobnerOverview}
We will refer often to Gr\"obner bases and Gr\"obner degenerations in this text. 
We refer to~\cite{CoxLittleOshea} and~\cite{MS} for detailed introductions to Gr\"obner theory, but we provide a lightning overview here.
Let $k[x_1, x_2, \ldots, x_M]$ be a polynomial ring. A \newword{term order} is a total order $\preceq$ on the monomials of $k[x_1, x_2, \ldots, x_M]$  such that (1) $x^a \preceq x^b$ implies $x^{a+c} \preceq x^{b+c}$, for $a$, $b$ and $c$ in $\ZZ_{\geq 0}^M$ and (2) $\preceq$ has no infinite descending chains. 
Given a nonzero polynomal $g(x) \in k[x_1, x_2, \ldots, x_M]$, the \newword{initial term} $\text{in}_{\preceq}(g)$ is the monomial in $g$ which is largest according to $\preceq$. 
Given an ideal $I \subseteq k[x_1, \ldots, x_M]$, the \newword{initial ideal} $\text{in}_{\prec}(I)$ is the ideal generated (and, in fact, spanned as a $k$-vector space) by $\{ \text{in}_{\preceq}(g) : g \in I,\ g \neq 0 \}$. 
A \newword{Gr\"obner basis} of $I$ is a set of elements of $I$ whose initial terms generate $\text{in}_{\preceq}(I)$; this automatically implies that the polynomials in this set generate $I$.

There is always a flat degeneration from $\Spec k[x_1, \ldots, x_M]/I$ to $\Spec k[x_1, x_2, \ldots, x_M]/\text{in}_{\preceq}(I)$. Moreover, this degeneration preserves any grading with respect to which $I$ is homogeneous. In particular, if $I$ is homogeneous with respect to the usual grading, then we get a flat degeneration from $\Proj k[x_1, \ldots, x_M]/I$ to $\Proj k[x_1, x_2, \ldots, x_M]/\text{in}_{\preceq}(I)$.

The ideal $\text{in}_{\prec}(I)$ is always a \newword{monomial ideal}, which implies in particular that the underlying point set of  $\Spec k[x_1, x_2, \ldots, x_M]/\text{in}_{\preceq}(I)$ is a union of coordinate hyperplanes. In nice cases, $\text{in}_{\prec}(I)$ is always a \newword{reduced monomial ideal}, also called a \newword{Stanley-Reisner ideal}, in which case $\Spec k[x_1, x_2, \ldots, x_M]/\text{in}_{\preceq}(I)$ is a reduced union of coordinate hyperplanes. In this case, we encode the combinatorics of this reduced union using a simplicial complex $\Delta$ on the vertex set $[M]$, where $F \subseteq [M]$ is a face of $\Delta$ if and only if the linear space $\Span_{f \in F} (e_f)$ is contained in $\Spec k[x_1, x_2, \ldots, x_M]/\text{in}_{\preceq}(I)$.

In the above setting, the minimal generators of $I$ correspond to the minimal non-faces of $\Delta$, and the irreducible components of $\text{in}_{\preceq}(I)$ correspond to the maximal faces of $\Delta$.
\end{remark}

\begin{proof}[Proof of Theorem~\ref{thm:2ndFundTheorem}]
The first statement is closely related to the ``second fundamental theorem of invariant theory" and goes back to Alfred Young~\cite{Young}. 
It is often attributed to Hodge and Pedoe~\cite[Chapter XIV]{HodgePedoe}.
For the Gr\"obner basis statement, see~\cite{SW} or \cite[Chapter 14]{MS}.
\end{proof}

Analogous results also hold for the quotients $\Pluck_u$ and $\Pluck^w$, which are the homogeneous coordinate rings of the Schubert varieties $X_u$ and $X^w$:
\begin{theorem}
The linear relations $\Delta_I=0$ for $I \not\succeq u([\# I])$, combined with the Pl\"ucker relations, generate the saturated homogeneous ideal of $\Pluck_u$ as a quotient of $\kappa[\Delta_I]$ and form a Gr\"obner basis for certain term orders. 
Similarly, the linear relations $\Delta_I=0$ for $I \not\preceq w([\# I])$, combined with the Pl\"ucker relations,   generate the saturated  homogeneous ideal of $\Pluck^w$ as a quotient of $\kappa[\Delta_I]$ and, again, form a Gr\"obner basis.
\end{theorem}

\begin{proof}
For the fact that these linear relations generate the ideal, see~\cite{EqnsDefiningSchubertVarieties}.
The fact that they are a Gr\"obner basis follows from~\cite[Theorem 7.1]{KLS1} but was surely known earlier.
\end{proof}

We will now describe a basis for the Pl\"ucker algebra in terms of semistandard Young tableaux.
Given a partition $\lambda$, the \newword{Young diagram} of $\lambda$ is a grid of boxes with $\lambda_i$ boxes in row $i$. 
We will draw our Young diagrams in the English convention, meaning that the rows are numbered from top to bottom.

\begin{eg}
This is the Young diagram of the partition $(4,2,1)$:
\[ \ydiagram{4,2,1} \]
\end{eg}

A \newword{tableau} of shape $\lambda$ is a filling of the boxes of $\lambda$ with positive integers.
A tableau is semistandard if the entries increase weakly from left to right across the rows and increase strictly from top to bottom down the columns. 
We will frequently abbreviate ``semistandard Young tableau" to \newword{SSYT}.
A \newword{reverse semistandard Young tableau} is a filling of the boxes of $\lambda$ with positive integers which \emph{decrease} weakly from left to right across the rows and \emph{decrease} strictly from top to bottom down the columns. 
Given a tableau $\T$, let $I_1$, $I_2$, \dots, $I_k$ be the sets of labels in the columns of $\T$, each read from top to bottom.
Define the \newword{Pl\"ucker monomial} $\Delta(\T)$ to be $\prod_{i=1}^k \Delta_{I_i}$. 

\begin{eg} \label{egSSYT}
This is a semistandard Young tableau of shape $(4,2,1)$:
\[ \begin{ytableau}
1 &1 & 2 & 3\\
2 & 3  \\
4 \\
\end{ytableau}. \]
The corresponding Pl\"ucker monomial is:
\[ \Delta_{124} \Delta_{13} \Delta_2 \Delta_3 = 
\left| \begin{smallmatrix} z_{11} & z_{12} & z_{13} \\ z_{21} & z_{22} & z_{23} \\ z_{41} & z_{42} & z_{43} \\ \end{smallmatrix} \right| \cdot
\left| \begin{smallmatrix} z_{11} & z_{12} \\ z_{31} & z_{32} \\ \end{smallmatrix} \right| \cdot z_{21} \cdot z_{31} . \]
\end{eg}

\begin{remark}
In general, there is no important difference between SSYTs and reverse SSYTs, and we have used SSYTs because they are more standard in the literature.
However, as we will discuss in Remark~\ref{rem:DiagonalVersusAntidiagonal}, using SSYTs would be incompatible with standard conventions in the field of matrix Schubert varieties, so we will switch to reverse SSYTs in Sections~\ref{sec:GTDegenFlag} and~\ref{sec:GTDegenRich}.
In anticipation of this, we will make remarks on the reverse case as appropriate.
\end{remark}

\begin{remark}
The word ``tableau'' is French in origin; the plural is ``tableaux".
\end{remark}

The \newword{content} of an SSYT $\T$ is the vector $(\alpha_1, \alpha_2, \ldots, \alpha_n)$ where $\T$ contains $\alpha_j$ copies of $j$.
We note the algebraic meaning of the shape and content of $\T$:
\begin{prop}
Let $\T$ be an SSYT of shape $\lambda$ and content $\alpha$. 
Then the Pl\"ucker monomial $\Delta(\T)$ has weight $\lambda$ for the right action of $T$ on $\GL_n$ and has weight $\alpha$ for the left action of $T$ on $\GL_n$.
\end{prop}

\begin{eg}
The SSYT in Example~\ref{egSSYT} has content $(2,2,2,1)$.
\end{eg}

The following follows from the Gr\"obner results of Theorem~\ref{thm:2ndFundTheorem} and a proof can be found in any of the sources cited therein.
\begin{theorem} \label{semistandardBasis}
The set of semistandard Pl\"ucker monomials $\Delta(\T)$, as $\T$ ranges over all semistandard Young tableaux with entries in $[n]$, is a basis for the Pl\"ucker algebra. 
\end{theorem}

\begin{remark}
Of course, the same holds if $\T$ ranges over reverse semistandard Young tableaux instead.
\end{remark}

Combining Theorem~\ref{semistandardBasis} with Remark~\ref{remarkSchur}, we deduce
\begin{cor}
The Kotska number $K_{\lambda \alpha}$ is the number of SSYT of shape $\lambda$ and content $\alpha$; the Schur polynomial $s_{\lambda}(z_1, z_2, \ldots, z_n)$ is $\sum_{\alpha} K_{\lambda \alpha} z_1^{\alpha_1} z_2^{\alpha_n} \cdots z_n^{\alpha_n}$. 
\end{cor}

\begin{eg}
\ytableausetup{smalltableaux}
The semistandard Young tableaux of shape $(2,1,0)$ and content $(1,1,1)$ are $\begin{ytableau}1&2 \\ 3 \end{ytableau}$ and $\begin{ytableau} 1&3 \\ 2 \end{ytableau}$.
They correspond to $\Delta_{13} \Delta_2$ and $\Delta_{12} \Delta_3$. 
The monomial $\Delta_{23} \Delta_1$ is not of the form $\Delta(\T)$, and we can write it in terms of the standard monomials using the Pl\"ucker relation $\Delta_{12} \Delta_3 - \Delta_{13} \Delta_2 + \Delta_{23} \Delta_1=0$.
If we used reverse tableaux, our basis would be $\{ \Delta_{32} \Delta_1, \Delta_{31} \Delta_2 \}$ instead.
\ytableausetup{nosmalltableaux}
\end{eg}

We now want to discuss the analogous results for $\Pluck_u^w$. 
First, however, we need a cautionary example to show that things may not be as nice as we could hope.

\subsection{The ideal of the Richardson versus the sum of the ideals of the Schuberts}
We defined $R_u^w$ as the intersection $X_u \cap X^w$. 
We know that $X_u$ is cut out by the equations $\Delta_I=0$ for $I \not\succeq u[\#(I)]$, and that $X^w$ is cut out by the equations $\Delta_I=0$ for $I \not\preceq w[\#(I)]$.
So, set theoretically, $R_u^w$, is cut out the equations $\Delta_I=0$ for $I$ \textbf{not} obeying $u[\#(I)] \preceq I \preceq w[\#(I)]$. 
We will now give an example where this is not true on the level of saturated ideals. 

\begin{eg}\label{eg:s2SquaredInFlags4}
Consider the Richardson $R_{1324}^{4231}$. 
The only Pl\"uckers which vanish on this Richardson are $\Delta_{12}$ and $\Delta_{34}$. 
We will show now that we have the relation
$\Delta_{123} \Delta_4 - \Delta_{124} \Delta_3=0$ on this Richardson, even though $\Delta_{123} \Delta_4 - \Delta_{124} \Delta_3$ is not in the ideal of the Pl\"ucker algebra generated by $\Delta_{12}$ and $\Delta_{34}$.
We remark that the relation $\Delta_{123} \Delta_4 - \Delta_{124} \Delta_3 + \Delta_{134} \Delta_2 - \Delta_{234} \Delta_1=0$ holds on the entire flag manifold.

Let us see why $\Delta_{123} \Delta_4 - \Delta_{124} \Delta_3=0$ on this Richardson.
The Richardson $R_{1324}^{4231}$ is the space of flags $(F_1, F_2, F_3)$ where $\dim (F_2 \cap \Span(e_1, e_2)) = \dim (F_2 \cap \Span(e_3, e_4)) = 1$. In other words, $F_2$ is of the form $\Span(x_1 e_1 + x_2 e_2, x_3 e_3 + x_4 e_4)$. 
The $3$-space $F_3$ is given by the equation $\Delta_{234} z_1 - \Delta_{134} z_2 + \Delta_{124} z_3 - \Delta_{123} z_4=0$.
The $3$-space $F_3$ contains $F_2$ and hence contains $x_3 e_3 + x_4 e_4$; we deduce that $\Delta_{124} x_3 - \Delta_{123} x_4=0$. 
Meanwhile, the line $F_1$ is contained in $F_2$, so it is of the form $\Span(y_1 (x_1 e_1 + x_2 e_2) + y_2 (x_3 e_3 + x_4 e_4))$; we deduce that $(\Delta_1 : \Delta_2 : \Delta_3 : \Delta_4) = (y_1 x_1 : y_1 x_2 : y_2 x_3 : y_3 x_4)$. Combining these observations, we have $\Delta_{123} \Delta_4 - \Delta_{124} \Delta_3=0$. 
Geometrically, what we have just shown is that the orthogonal projection of $F_1$ onto $\Span(e_3, e_4)$ coincides with the intersection $F_3 \cap \Span(e_3, e_4)$. 
 
This example shows that the defining ideal of $R_u^w$ in the Pl\"ucker algebra can be larger than the sum of the ideals of $X_u$ and $X^w$: For both $X_{1324}$ and $X^{4231}$, the only relation in the $(1,1,1,1)$ weight space of $H^0(\cL(2,1,1,0))$ is $\Delta_{123} \Delta_4 - \Delta_{124} \Delta_3 + \Delta_{134} \Delta_2 - \Delta_{234} \Delta_1=0$.

We will close this example by noting that no similar issue occurs in $H^0(\cL(3,2,1,0))$:  For every $1 \leq i<j \leq 4$, the polynomial $\Delta_{ij} \left( \Delta_{123} \Delta_4 - \Delta_{124} \Delta_3 \right)$ is in the sum of the ideals of $X_{1324}$ and $X^{4231}$.
For example,
\begin{multline*}
 \Delta_{13} \left( \Delta_{123} \Delta_4 - \Delta_{124} \Delta_3 \right) = \\
\Delta_{123} (\Delta_{13} \Delta_4 - \Delta_{14} \Delta_3 + \Delta_{34} \Delta_1)  + (\Delta_{123} \Delta_{14} - \Delta_{124} \Delta_{13} + \Delta_{134} \Delta_{12}) \Delta_3  
- \Delta_{123} \Delta_{34} \Delta_1 - \Delta_{134} \Delta_{12} \Delta_3 . 
\end{multline*}

 The trinomials are Pl\"ucker relations and hence vanish on the whole flag manifold; $\Delta_{123} \Delta_{34} \Delta_1$ vanishes on $X^{4231}$ since it is divisible by $\Delta_{34}$ and  $\Delta_{134} \Delta_{12} \Delta_3$ vanishes on $X_{1324}$ since it is divisible by $\Delta_{12}$.
 Since similar formulas exist for each $\Delta_{ij} \left( \Delta_{123} \Delta_4 - \Delta_{124} \Delta_3 \right)$, the binomial $\Delta_{123} \Delta_4 - \Delta_{124} \Delta_3$ is in the saturation of the sum of the ideals of $X_u$ and $X^w$.
\end{eg}

In the previous example, we showed that the defining ideal of $R_u^w$ can be larger than the sum of the ideals of $X_u$ and $X^w$. 
However, this extent to which this can occur is limited by a theorem of Lakshmibai and Littelmann~\cite[Theorem 16]{LL03}:
\begin{theorem} \label{thm:LLSumTheorem}
Let $u \preceq w$ in $S_n$. Let $\lambda_1 > \lambda_2 > \cdots > \lambda_n \geq 0$.
Then the ideal of $R_u^w$ in degree $\lambda$ is the sum of the ideals of $X_u$ and $X^w$. 
\end{theorem}

Theorem~\ref{thm:LLSumTheorem} is best understood in terms of sheaf cohomology: 
It says that, for $\cL(\lambda)$ an ample bundle, the kernel of $H^0(\Fl_n, \cL(\lambda)) \longrightarrow H^0(R_u^w, \cL(\lambda))$ is the sum of the kernels of $H^0(\Fl_n, \cL(\lambda)) \longrightarrow H^0(X_u, \cL(\lambda))$ and $H^0(\Fl_n, \cL(\lambda)) \longrightarrow H^0(X^w, \cL(\lambda))$.
Letting $\cI$ be the reduced ideal sheaf of $X_u \cup X^w$, this comes down to verifying that $H^1(\Fl_n, \cI \otimes \cL(\lambda))=0$.
See \cite{GeomMonTheory} for a proof from this perspective.

We saw this in Example~\ref{eg:s2SquaredInFlags4}, where the sum of the Schubert ideals wasn't large enough in degree $(2,1,1,0)$, but it was large enough in degree $(3,2,1,0)$.

\subsection{Standard monomial theory and semistandard tableaux} \label{sec:SMT}
Brion and Lakshmibai~\cite{GeomMonTheory} give the following basis of $\Pluck_u^w$:
\begin{TD} \label{thm:StandardMonomialTheory}
Let $u \preceq w$ in $S_n$. Let $\mu$ be a partition and write $\mu = \omega_{k_1} + \omega_{k_2} + \cdots + \omega_{k_m}$ with $k_1 \geq k_2 \geq \cdots \geq k_m$. 
In other words, $(k_1, k_2, \ldots, k_m)$ is the transpose partition to $\mu$.

Let $I_1$, $I_2$, \dots, $I_m$ be a sequence of subsets of $[n]$ with $\#(I_j) = k_j$.
 We define $(I_1, I_2, \ldots, I_m)$ to be \newword{standard for $(u,w)$} if there is a chain $u \preceq v_1 \preceq v_2 \preceq \cdots \preceq v_m \preceq w$ in $S_n$ with $v_j([k_j]) = I_j$.
Then the Pl\"ucker monomials $\prod_{j=1}^m \Delta_{I_j}$ where $(I_1, I_2, \ldots, I_m)$ is standard for $(u,w)$ form a basis of $H^0(R_u^w, \cL(\mu))$.
\end{TD}

\begin{remark}
We have switched our name for the shape of a Young diagram from $\lambda$ to $\mu$ in order to avoid conflict with the notation $\lambda_i$ in Definition~18 of~\cite{GeomMonTheory}.
\end{remark}

\begin{remark}
We will see in Lemma~\ref{lem:StandardIsStandard} that the condition that such a chain exists implies that $I_1$, $I_2$, \dots, $I_m$ are the columns of an SSYT of shape $\mu$ and that, for $(u,w) = (e,w_0)$, the chain condition is equivalent to imposing that $I_1$, $I_2$, \dots, $I_m$ are the columns of an SSYT of shape $\mu$.
\end{remark}

\begin{remark} \label{rem:StandardIsSimplicial}
It is easy to see that, if $(I_1, I_2, \dots, I_m)$ is standard, then so is any sequence  $(I_{a_1}, I_{a_2}, \dots, I_{a_r})$ for $1 \leq a_1 \leq a_2 \leq \cdots \leq a_r \leq m$.
This means that we can think of the standard monomials as faces of a simplicial complex; see Section~\ref{sec:SimplicialResults} for more on this perspective.
\end{remark}

\begin{remark}
If we were using reverse SSYTs, then we would ask that $w \succeq v_1 \succeq v_2 \succeq \cdots \succeq v_m \succeq u$ with $v_j([k_j]) = I_j$.
The sharp eyed reader will note that this is actually the ordering convention in~\cite{GeomMonTheory}.
\end{remark}

\begin{proof}[Proof sketch]
This is Theorem~3 in~\cite{GeomMonTheory}, but that paper is written in a very high level of generality.
We explain how to extract the concrete description of standard monomials here from Definition~18 of~\cite{GeomMonTheory}.

Definition~18 speaks of a sequence of dominant characters $\lambda_1$, $\lambda_2$, \dots, $\lambda_m$; these will be  $\omega_{k_1}$, $\omega_{k_2}$, \dots, $\omega_{k_m}$.
The standard monomial basis is defined as certain products $p_{w(\lambda_1), v(\lambda_1)} p_{w(\lambda_2), v(\lambda_2)} \cdots p_{w(\lambda_m), v(\lambda_m)}$ where $p_{w(\lambda_i), v(\lambda_i)}$ is in $H^0(\Fl_n, \cL(\lambda_i))$. 
In our case, $H^0(\Fl_n, \cL(\omega_{k_i}))$ is $\Alt^{k_i} \AA^n$ and the $p_{v,w}$'s will be the Pl\"ucker coordinates $\Delta_I$ with $\#(I) = k_i$.
For a general character $\lambda$ ``of classical type", the $p_{v,w}$ are indexed by a pair $(v,w)$ of elements of $W^{\lambda}$ but, when $\lambda$ is ``minuscule'', as in our case, we always have $v=w$.

The notation $W^{\lambda}$ means the quotient of the Weyl group $W$ by the stabilizer of the character $\lambda$.
In our case, the stabilizer of $\omega_k$ is $S_k \times S_{n-k}$ and the quotient $S_n/(S_k \times S_{n-k})$ is identified with $\binom{[n]}{k}$ by sending the coset $w (S_k \times S_{n-k})$ to $w[k]$.
Thus, $p_{v,v}$ means $\Delta_{v[k_i]}$ for $v$ in $W^{\lambda_i}$.
The condition that there exists a chain $u \preceq \tilde{v}_1 \preceq \tilde{w}_1 \preceq \cdots \preceq \tilde{v}_m \preceq \tilde{w}_m \preceq w$ with $\tilde{v}_i$ and $\tilde{w}_i$ in the cosets $v_i W^{\lambda_i}$ and $w_i W^{\lambda_i}$ then reduces to asking for a single chain $u \preceq v_1 \preceq v_2 \preceq \cdots \preceq v_m \preceq w$ with $v_i([k_i]) = I_i$.
\end{proof}
%

\begin{eg}  \label{eg:Weight111In210}
Let us consider the $(1,1,1)$ weight space of $H^0(\cL(2,1,0))$, on the flag manifold and on various Richardson submanifolds of it. 
The three Pl\"ucker monomials in this weight space are $\Delta_{12} \Delta_3$, $\Delta_{13} \Delta_2$ and $\Delta_{23} \Delta_1$. 
Of these, the last does not come from an SSYT. 
We can also check directly that the last monomial cannot be lifted to a chain $v_1 \preceq v_2$ as we would have to have $v_1[1] \subset v_1[2] = 23$ so $v_1[1]$ is one of $\{ 2, 3 \}$, and we would also have to have $v_2[1] = 1$. But then the inequality $v_1[1] \preceq v_2[1]$ cannot hold.

The other two monomials do come from SSYT, and they do lift to chains. The first monomial can be lifted to any of the four chains $213 \preceq 312$, $123 \preceq 312$, $213 \preceq 321$ and $123 \preceq 321$; the second monomial can be lifted only to $132 \preceq 231$. This shows that the first monomial is standard for $(u,w)$ if and only if $u \preceq 213  \preceq 312 \preceq w$ and the second monomial is standard for $(u,w)$ if and only if  $u \preceq 132  \preceq 231 \preceq w$ .

We note that the fact that $\Delta_{13} \Delta_2$ is nonstandard on (for example) $X_{213}$ doesn't mean that $\Delta_{13} \Delta_2$ is $0$ on $X_{213}$.
Instead, we have $\Delta_{12} \Delta_3 = \Delta_{13} \Delta_2$ on $X_{213}$, as we can see by taking the Pl\"ucker relation $\Delta_{12} \Delta_3 - \Delta_{13} \Delta_2 + \Delta_{23} \Delta_1$ and plugging in the fact that $\Delta_1 =0$ on  $X_{213}$
\end{eg}


\begin{eg} \label{eg:s2SquaredInFlags4Part2}
Continuing from Example~\ref{eg:s2SquaredInFlags4}, we consider the Richardson $R_{1324}^{4231}$. 
Look first at the weight space $(1,1,1,1)$ for $\mu = (2,1,1,0)$.
There are three semistandard Young tableaux of this shape and content, corresponding to the Pl\"ucker monomials $\Delta_{123} \Delta_4$, $\Delta_{124} \Delta_3$ and $\Delta_{134} \Delta_2$.
In this Richardson, we have $\Delta_{124} \Delta_3 = \Delta_{123} \Delta_4$, so one of these two must be nonstandard.
Indeed, we claim that $\Delta_{124} \Delta_3$ is not standard. Suppose, to the contrary, that we had a chain $1324 \preceq v_1 \preceq v_2 \preceq 4231$ with $v_1[3] = 124$ and $v_2[1] = 3$. 
The condition that $1324 \preceq v_1 \preceq v_2 \preceq 4231$ implies that $v_1[2]$ and $v_2[2]$ must be one of $\{ 13, 14, 23, 24 \}$. 
Since $v_1[2] \subset v_1[2] = 124$, we deduce that $v_1[2]$ is one of $\{ 14, 24 \}$.
Since $v_2[2] \supset v_2[2] = 3$, we deduce that $v_2[2]$ is one of $\{ 13, 23 \}$. 
But now the inequality $v_1[2] \preceq v_2[2]$ cannot hold.
So there are only two standard monomials on this Richardson, and  $(1,1,1,1)$ for $\mu = (2,1,1,0)$ is two dimensional.

We also consider the $(2,1,2,1)$ weight space for $\mu=(3,2,1,0)$ on this Richardson. 
There are $4$ SSYT for this example, corresponding to the Plucker monomials $\Delta_{123} \Delta_{13} \Delta_4$,  $\Delta_{123} \Delta_{14} \Delta_3$, $\Delta_{124} \Delta_{13} \Delta_3$, $\Delta_{134} \Delta_{13} \Delta_2$.
The first and last of these are standard on $R_{1324}^{4231}$ and the others are not. 
Indeed, $\Delta_{14} \Delta_3$ is already not standard on $X^{4231}$ and $\Delta_{124} \Delta_{13}$ is already not standard on $X_{1324}$.
So, again, the $(2,1,2,1)$ weight space for $\mu=(3,2,1,0)$ is two dimensional.
 \end{eg}

It is not obvious that this definition of ``standard" collapses to the definition of a semistandard Young tableau in the case $(u,w) = (e, w_0)$. 
We now check this.
\begin{lemma} \label{lem:StandardIsStandard}
Let $n \geq k_1 \geq k_2 \geq \cdots \geq k_m \geq 1$ and let $I_1$, $I_2$, \dots, $I_m$ be a sequence of subsets of $[n]$ with $\#(I_j) = k_j$.
Then $I_1$, $I_2$, \dots, $I_m$ are the columns of an SSYT (in that order) if and only if there is a chain $v_1 \preceq v_2 \preceq \cdots \preceq v_m$ in Bruhat order with $v_j[k_j] = I_j$.
\end{lemma}

\begin{proof}[Proof sketch]
First, suppose that there is a chain $v_1 \preceq v_2 \preceq \cdots \preceq v_k$  with $v_j[k_j] = I_j$.
Enter the elements of $I_j$ into the $j$-th column of $\mu$, in increasing order. We must check that the rows are increasing.
Define $\{ x_1 < x_2 < \cdots < x_{k_j} \}:=I_j = v_j[k_j] $, $\{ y_1 < y_2 < \cdots < y_{k_j} \} := v_{j+1}[k_j]$ and $\{ z_1 < z_2 < \cdots < z_{k_{j+1}} \} := I_{j+1}= v_{j+1}[k_{j+1}]$. 
So $x_i$ is the tableau entry in row $i$, column $j$, and $z_i$ is the tableau entry in row $i$, column $j+1$.
Since $v_j \preceq v_{j+1}$, we have $v_j[k_j] \preceq  v_{j+1}[k_j]$ and thus $x_i \leq y_i$ for $1 \leq i \leq k_j$.
Since $k_j \geq k_{j+1}$, we have $v_{j+1}[k_j] \supseteq v_{j+1}[k_{j+1}]$ and thus $y_i \leq z_i$ for $1 \leq i \leq k_{j+1}$.
Concatenating these relations, we deduce that $x_i \leq z_i$ for $1 \leq i \leq k_{j+1}$.
We have now checked that the rows of our tableau are weakly increasing.

For the reverse direction, we recall that a \newword{mountain permutation} is a permutation which, in one line notation, has the form $a_1 a_2 \cdots a_{j-1} n b_1 b_2 \cdots b_{n-j}$ for $a_1 < a_2 < \cdots < a_{j-1} < n > b_1 > b_2 > \cdots > b_{n-j}$. For a subset $I$ of $[n]$, we define $\text{mount}(I)$ as follows: If $n \not\in I$, then  $\text{mount}(I)$ is the unique mountain permutation with $\{ a_1, a_2, \ldots, a_j \} = I$; if $n \in I$ then $\text{mount}(I)$ is the unique mountain permutation with $\{ a_1, a_2, \ldots, a_j \} = I \setminus \{ n \}$.
We leave it to the reader to check that, if $I_1$, $I_2$, \dots, $I_k$ are the columns of an SSYT, then defining $v_j = \text{mount}(I_j)$ gives a chain in Bruhat order with the desired property.
\end{proof}

\begin{eg}
The figure below demonstrates the strategy of the second half of the proof of Lemma~\ref{lem:StandardIsStandard}. 
The left hand side is a Hasse diagram where $I$ is below $J$ if $I$ can precede $J$ in an SSYT.
The right hand side shows the corresponding mountain permutations.
SSYT correspond to weak chains on the left hand side, and the reader can check that these become weak chains on the right hand side.

\[
\xymatrix@R=0.1in@C=0.1 in{
&&&4& \\
&&3  \ar@{-}[ur]&& \\
&2  \ar@{-}[ur] &&34  \ar@{-}[ul] & \\
1 \ar@{-}[ur]&&24\ar@{-}[ur] \ar@{-}[ul] && \\
&14 \ar@{-}[ur] \ar@{-}[ul] &&23\ar@{-}[ul]& \\
&&13 \ar@{-}[ur] \ar@{-}[ul] &&234 \ar@{-}[ul]  \\
&12 \ar@{-}[ur]&& 134 \ar@{-}[ur] \ar@{-}[ul] & \\
&& 124 \ar@{-}[ur] \ar@{-}[ul] & &\\
& 123 \ar@{-}[ur] & &&\\
 && 1234 \ar@{-}[ul] &&\\
} \qquad
\xymatrix@R=0.1in@C=0.1 in{
&&&4321& \\
&&3421 \ar@{-}[ur]&& \\
&2431  \ar@{-}[ur] &&3421  \ar@{=}[ul] & \\
1432 \ar@{-}[ur]&&2431\ar@{-}[ur] \ar@{=}[ul] && \\
&1432 \ar@{-}[ur] \ar@{=}[ul] &&2341\ar@{-}[ul]& \\
&&1342 \ar@{-}[ur] \ar@{-}[ul] &&2341 \ar@{=}[ul]  \\
&1243 \ar@{-}[ur]&& 1342 \ar@{-}[ur] \ar@{=}[ul] & \\
&& 1243 \ar@{-}[ur] \ar@{=}[ul] & &\\
& 1234 \ar@{-}[ur] & &&\\
 && 1234 \ar@{=}[ul] &&\\
} \qquad
\]
\end{eg}

As we have described it, computing whether $\T$ is standard for $(u,w)$ requires considering all lifts of $\T$ to the Bruhat order.
In fact, one can carry this out greedily, as we will now discuss. 
\begin{lemma}
Let $u \in S_n$ and let $J$ be a $k$-element subset of $[n]$ with $u[k] \preceq J$. Then there is a unique $\preceq$-minimal permutation in the set of $x$ such that $u \preceq x$ and $x[k] = J$.
\end{lemma}

\begin{proof}[Sketch of proof]
We are supposed to construct $x$ with $x[k] = J$ and $x[k+1, n] = [n] \setminus J$. We will describe in which order $x$ maps the elements of $[k]$ to those of $J$; one can similarly determine in what order $x$ maps the elements of $[k+1,n]$ to those of $[n] \setminus J$.

Let $r \leq k$. Let $X = \{ K \in \binom{J}{r} : K \succeq u[r] \}$. Now, $\binom{J}{r}$ is a distributive lattice with respect to $\preceq$, and $X$ is clearly closed under the meet operation, so $X$ has a minimal element, call it $J_r$. We claim that $J_1 \subset J_2 \subset \cdots \subset J_k = J$. Indeed, suppose that $J_r \not\subset J_{r+1}$. Since $J_{r+1} \succeq u[r+1]$, there is some $r$-element subset $J'$ of $J_{r+1}$ with $J' \succeq u[r]$; choose the minimal such $J'$ and let $j$ be the single element of $J_{r+1} \setminus J'$. Since $J' \succeq u[r]$, we have $J' \succeq J_r$, and they are not equal since $J_r \not\subset J_{r+1}$, so $J' \succ J_r$. Then $J_{r+1} \succ J_r \cup \{ j \} $ and $J_r \cup \{ j \} \succeq u[r+1]$, contradicting the minimality of $J_{r+1}$. We now know that $J_1 \subset J_2 \subset \cdots \subset J_k = J$.

Take $x_r$ to be the lone element of $J_r \setminus J_{r-1}$ (where $J_0 = \emptyset$). This determines $x_1$,  $x_2$, \dots, $x_k$; we can find $x_{k+1}$, $x_{k+2}$, \dots, $x_n$ similarly.
\end{proof}

Thus, we have the following algorithm to determine whether $\T$ is standard for $(u,w)$: Let $I_1$, $I_2$, \dots, $I_m$ be the columns of $\T$. Set $x_0 = u$ and let $x_j$, inductively, be the minimal permutation with $x_j \succeq x_{j-1}$ and $x_j[k_j] = I_j$. If at any point no such permutation exists, then $\T$ is not standard, and if $x_m \not\preceq w$, then $\T$ is not standard; otherwise, $\T$ is standard.

\begin{proof}[Sketch of proof of correctness]
Let $x_0$, $x_1$, \dots, $x_m$ be as above. If the $x_i$ are defined and $x_m \preceq w$, then $x_1$, $x_2$, \dots, $x_m$ is a chain with $u \preceq x_1 \preceq x_2 \preceq \cdots \preceq x_m \preceq w$ lifting $\T$.

Conversely, suppose that $u \preceq v_1 \preceq v_2 \preceq \cdots \preceq v_m \preceq w$ is a chain lifting $\T$. Then, inductively, we have $v_j \succeq x_j$ for all $j$ (using the minimality of $x_j$). So, in particular, the $x_j$ exist and $x_m \preceq v_m \preceq w$.
\end{proof}

\begin{remark}
If $u$ is $e$, then the minimal $x_m$ described above is called the \newword{right key of $\T$}.
Keys were introduced by Lascoux and Schutzenberger~\cite{KeysStandardBases}, and a simple way of computing them was found by Willis~\cite{SimpleKey}.
\end{remark}

\begin{eg}
We verify again that the tableau with columns $(124, 3)$ is not standard for $(1324, 4231)$.
The minimal $v_1$ with $v_1 \succeq 1324$ and $v_1[3] = 124$ is $1423$; the minimal $v_2$ with $v_2 \succeq 1423$ and $v_2[1] = 3$ is $3412$. 
We have $3412 \not\preceq 4231$. 
\end{eg}

We have discussed the part of Brion and Lakshmibai's paper which is most combinatorial and which most closely generalizes Hodge's classical standard monomial theory.
Brion and Lakshmibai, and the earlier work of Laksmibai with various co-authors under the name ``standard monomial theory", constructs many different bases of $H^0(R_u^w, \cL(\lambda))$. 
We will discuss some of these other results in Section~\ref{sec:FrobeniusSplitting}.

\subsection{The simplicial complex of standard monomials} \label{sec:SimplicialResults}

Let $u \preceq w$ and let $I_1$, $I_2$, \dots, $I_m$ be nonempty subsets of $[n]$.
As pointed out in Remark~\ref{rem:StandardIsSimplicial}, if $I_1$, $I_2$, \dots, $I_m$ are the columns of an SSYT which is standard for $[u,w]$, then $I_{a_1}$, $I_{a_2}$, \dots, $I_{a_r}$ are also the columns of an SSYT which is standard for $[u,w]$ for any increasing sequence $1 \leq a_1 \leq a_2 \leq \cdots \leq a_r \leq m$.
In other words, there is a simplicial complex $\Delta(u,w)$ whose vertices are the nonempty subsets of $[n]$ and where $I_1$, $I_2$, \dots, $I_m$ are the vertices of a face if and only if $I_1$, $I_2$, \dots, $I_m$ are the columns of an SSYT which is standard for $[u,w]$.

If we restrict to the case that all the $I_j$ have the same cardinality $k$, then this simplicial complex describes the Gr\"obner degeneration of a positroid variety \cite[Section 7]{KLS1}.
Call this complex $\Delta_k(u,w)$.
The complex $\Delta_k(u,w)$ is studied in~\cite{KLS1}. 
We describe the primary results of that paper:

\begin{prop}
Let $u \preceq_k w$ and $u' \preceq_k w'$ be two intervals in $k$-Bruhat order which index the same positroid variety, meaning that $\Pi_u^w = \Pi_{u'}^{w'}$. Then $\Delta_k(u,w) =  \Delta_k(u', w')$. 
\end{prop}

Thus, it makes sense to associate the simplicial complex $\Delta_k(u,w)$ to the positroid variety $\Pi_u^w$.

\begin{theorem}
The simplicial complex $\Delta_k(u,w)$ is pure of dimension $\ell(w) - \ell(u)$. The map from saturated $\preceq_k$-Bruhat chains to facets of $\Delta_k(u,w)$ is bijective. The simplicial complex $\Delta_k(u,w)$ is shellable and is homeomorphic to a ball. 
\end{theorem}

Here a $d$-dimensional simplical complex is called \newword{shellable} if we can order the $d$-dimensional simplices as $\sigma_1$, $\sigma_2$, \dots, $\sigma_N$ such that $\sigma_j \cap \bigcup_{i=1}^{j-1} \sigma_i$ is pure of dimension $d-1$ (or empty) for all $j$. This condition implies that the simplicial complex is homotopy equivalent to a wedge of $d$-dimensional spheres, and has other important consequences for combinatorial topology and commutative algebra.
In particular, shellability of $\Delta_k(u,w)$ implies that the homogeneous coordinate ring of $\Pi_u^w$ is Cohen-Macaulay.
See Remark~\ref{GrobnerOverview} for a general discussion of Gr\"obner bases.

\begin{eg} \label{eg:G24Complex1ForBS}
We depict the interval $[2143, 3412]$ in $\preceq_2$:
\[  \xymatrix@R=0.1 in@C=0.1 in{
& 3412 \ar@{-}[dr] \ar@{-}[dl] & \\
3142 \ar@{-}[dr] && 2413 \ar@{-}[dl] \\
& 2143 & \\
} \]
The facets of $\Delta_2(2143, 3412)$ are $\{ 12, 13, 34 \}$ and  $\{ 12, 24, 34 \}$.
 We note that, although $13 \prec 24$ in the partial order on $\binom{[4]}{2}$, there is no face of $\Delta_2(2143, 3412)$ which contains both $13$ and $24$.

The corresponding positroid is $\Delta_{14} = \Delta_{23} = 0$;  the homogeneous coordinate ring of this positroid variety is $\kappa[\Delta_{12}, \Delta_{13}, \Delta_{24}, \Delta_{34}]/(\Delta_{13} \Delta_{24} - \Delta_{12} \Delta_{34})$. As promised, the monomials of the forms $\Delta_{12}^a \Delta_{13}^b \Delta_{34}^c$ and $\Delta_{12}^a \Delta_{24}^b \Delta_{34}^c$ form a basis for this ring.
\end{eg}

\begin{eg}
 \label{eg:G24Complex2ForBS}
We depict the interval $[2134, 3412]$ in $\preceq_2$. The solid lines are covers in $\preceq_2$; the dashed line is a $\preceq$-cover which is not a $\preceq_2$ relation.
\[  \xymatrix@R=0.1 in@C=0.1 in{
& 3412 \ar@{-}[dr] \ar@{-}[dl] & \\
2413 \ar@{-}[dd] && 3214 \ar@{-}[dd] \ar@{--}[ddll] \\
&&\\
2314 \ar@{-}[dr] && 3124 \ar@{-}[dl] \\
& 2134 & \\ 
} \]
The facets of $\Delta_2(2134, 3412)$ are $\{ 12, 23, 24,  34 \}$ and $\{ 12, 13, 23,  34 \}$, corresponding to the two maximal $\preceq_2$ chains. Note that the $\preceq$ chain $(2134, 2314, 3214, 3412)$, which is not a $\preceq_2$-chain, collapses to the lower dimensional face $\{ 12, 23, 34 \}$ of $\Delta_2(2134, 3412)$.
Note also that $\Delta_2(2134, 3412)$ is shellable, whereas the order complex of $[2134, 3412]$ with respect to $\preceq_2$ is not shellable.

The corresponding positroid is $\Delta_{14} = 0$;  the homogeneous coordinate ring of this positroid variety is $\kappa[\Delta_{12}, \Delta_{13}, \Delta_{14}, \Delta_{24}, \Delta_{34}]/(\Delta_{13} \Delta_{24} - \Delta_{12} \Delta_{34})$.
As promised, the monomials of the forms $\Delta_{12}^a \Delta_{23}^b \Delta_{24}^c \Delta_{34}^d$ and $\Delta_{12}^a \Delta_{13}^b \Delta_{23}^c \Delta_{34}^d$ form a basis for this ring.
\end{eg}

Almousa, Gao and Huang~\cite{AGH} have found an explicit description of the minimal nonfaces of $\Delta_k(u,w)$. In particular, they can be of arbitrarily large cardinality. 

We do not know much about the full complex $\Delta(u,w)$.
\begin{prob}
Is $\Delta(u,w)$ shellable? Is there a simple description of its maximal faces? These will be related to maximal components in flat degenerations of $X_u^w$; see Remark~\ref{GrobnerOverview}.
\end{prob}

We do not have a good understanding of the minimal nonfaces of $\Delta(u,w)$, but we do know that they can be of arbitrarily high cardinality. 
Such minimal nonfaces will be related to minimal generators of the ideal of $X_u^w$; see Remark~\ref{GrobnerOverview}.
The following example is simplified from~\cite[Example IV.31]{Kim}:
\begin{eg}
Let  $I_{k} = [k+1] \setminus \{ k \} = \{1,2,\ldots, k-1, k+1 \}$. Let $\T$ be the SSYT with columns $I_{n-1}$, $I_{n-2}$, \dots, $I_2$, $I_1$.
Below, we depict $\T$ for $n=5$:
\[ \begin{ytableau}
1 & 1 & 1 & 2 \\
2 & 2 & 3 \\
3 & 4 \\
5 \\
\end{ytableau}. \]
Let $w = s_1 w_0 = n(n-1)(n-2) \cdots 54312 $. 
We claim that $\T$ is not standard for $X^w$, but that deleting any column from $\T$ gives a tableau which is standard for $X^w$. 
In other words, $\T$ gives an $(n-2)$-dimensional minimal nonface of $\Delta(e,w)$. 

We first verify that $\T$ is not standard for $(e,w)$. 
The minimal lift of $(I_{n-1}, I_{n-2}, \ldots, I_2, I_1)$ to $S_n$ is $v_{n-1} \prec v_{n-2} \prec \cdots \prec v_1$ where $v_k$ is $123\cdots \widehat{k} \cdots (n-1) n k$.
In particular, $v_1 = 234 \cdots 1$. Since $v_1[n-1] = [n] \setminus \{ 1 \}$ and $w[n-1] = [n] \setminus \{ 2 \}$, we have $v_1 \not\preceq w$, and $\T$ is not standard for $(e,w)$. 

We now sketch the proof that deleting any column from $\T$ gives a tableau which \textbf{is} standard for $X^w$. 
Let $(u_{n-1}, u_{n-2}, \ldots, u_{k+1}, u_{k-1}, \ldots, u_1)$ be the minimal lift of $I_{n-1} I_{n-2} \cdots \widehat{I_k} \cdots I_2 I_1$.
For $j \geq k+1$, we have $u_j = 123\cdots \widehat{j} \cdots (n-1) n j$ as in the previous paragraph; in particular, $u_{k+1}(n) = k+1$. 
Since $k+1$ does not belong to any of $I_{k-1}$, $I_2$, \dots, $I_1$, we will continue to have $u_j(n) = k+1$ for $j \leq k-1$.
In particular, $u_1(n) = k+1$ so $u_1[n-1] = [n] \setminus \{ k+1 \}$ which is $\preceq w[n-1] = [n] \setminus \{ 1 \}$.
For $m < n-1$, we have $w[m] = \{ n-m+1, n-m+2, \ldots, n \}$ so the condition $u_1[m] \preceq w[m]$ is automatic.
\end{eg}

\subsection{Gr\"obner degeneration of matrix Schubert and matrix Richardson varieties} \label{sec:MatrixSchubert}
Richardson varieties are subvarieties of the flag manifold.
Standard Gr\"obner basis techniques are designed to study subvarieties of affine (or projective) spaces.
In this section, we will discuss the extremely successful theory of Gr\"obner basis for matrix Schubert varieties, and briefly discuss the possibility of building a similar theory for ``matrix Richardson varieties".
Many of the ideas in this section will then reappear when we discuss degeneration of the Richardson varieties in the flag manifold.

Recall that we have the Bruhat decomposition $\GL_n = \bigsqcup_{u \in S_n} B_- u B_+$. 
The general linear group is open and dense inside the affine space $\Mat_{n \times n}$ of $n \times n$ matrices.
We define the \newword{matrix Schubert variety} $M_u$ to be the closure $\overline{B_- u B_+}$. 
The idea of studying the matrix Schubert varieties in place of Schubert varieties was introduced by Fulton~\cite{Fulton92}.

\begin{rem} \label{rem:InverseWarning}
We warn the reader that, in~\cite{KnutsonMiller} and~\cite{MS}, and much literature derived from them, $M_u$ is defined to be $\overline{B_- u^{-1} B_+}$. 
For example, such papers would say that $M_{s_1 s_2}$ is $\{ z_{11} = z_{21} = 0 \}$, whereas we say that $M_{s_1 s_2}$ is $\{ z_{11} = z_{12} = 0 \}$. 
Observe that the matrix $s_1 s_2 = \begin{sbm} 0&0&1 \\ 1&0&0 \\ 0&1&0 \\ \end{sbm}$ is in $M_{s_1 s_2}$ with our conventions, and not with the other conventions.

This issue arises because the authors  identify $\Fl_u$ with the quotient $B_- \leftquot \GL_n$ rather than $\GL_n/B_+$.
The double coset $B_- u B_+$ is the same subset of $\GL_n$ in either way, but it gives rise to different cells in $\Fl_n$ depending on which quotient is taken.
Any reader whose primary interest is the Knutson-Miller theory will have to make the same choice, as it is omnipresent in the literature.
However, as this is not our main topic, we have chosen to follow our convention that we identify $\Fl_n$ with $\GL_n/B_+$ and thus $M_u$ is $\overline{B_- u B_+}$.
\end{rem}

We write $z_{ij}$ for the coordinates on $\Mat_{n \times n}$.
 For $I$, $J \subseteq [n]$ with $\#(I) = \#(J)$, we write $\Delta_{I, J}$ for the minor $\det [z_{ij}]_{i \in I,\ j \in J }$.
 
 We choose an antidiagonal term order on $\kappa[z_{11}, \ldots, z_{nn}]$, meaning that, for any $I = \{ i_1 < i_2 < \ldots < i_k \}$ and $J = \{ j_1< j_2 < \ldots < j_k \}$, the leading monomial in $\Delta_{I,J}$ is $z_{i_1 j_k} z_{i_2 j_{k-1}} \cdots z_{i_k j_1}$.
 For $f \in \kappa[z_{11}, \ldots, z_{nn}]$, we write $\In(f)$ for the leading monomial of $f$ in our term order.
 For an ideal $\cI \subseteq \kappa[z_{11}, \ldots, z_{nn}]$, we write $\In(\cI)$ for the ideal generated by $\{ \In(f) : f \in \cI \}$.
 If $\cX$ is the subscheme of $\Mat_{n \times n}$ corresponding to $\cI$, then we will write $\In(\cX)$ for $\Spec \kappa[z_{ij}]/\In(I)$. 
 General Gr\"obner theory tells us there will be a flat family over $\AA^1$ whose fibers over every nonzero point are isomorphic to $\cX$ and whose fiber over the $0$ point is isomorphic to  $\In(\cX)$.
 Moreover, $\In(\cI)$ will be graded for any grading for which $\cI$ is graded, so we can work with $\Proj$ or $\MultiProj$ instead of $\Spec$.

 Let $u \in S_n$ and let $\cI_u$ be the radical ideal of the matrix Schubert variety $M_u$.
 In Theorem~\ref{thm:FultonGenerators}, we learned that $\cI_u$ is generated by the minors $\Delta_{I,J}$ for $I \subseteq [a]$, $J \subseteq [b]$, $\#(I) = \#(J) > \#([a] \cap u[b])$.
 These are known as the \newword{Fulton generators} of $\cI_u$.
Knutson and Miller proved that this list of generators is a Gr\"obner basis.
Explicitly, this means:
\begin{theorem} \label{thm:KnutsonMiller}
 Let $u \in S_n$ and let $\cI_u$ be the reduced ideal of the matrix Schubert variety $M_u$.
 The ideal $\In(\cI_u)$ is generated by the monomials $\In(\Delta_{I,J})$, for $\Delta_{(I,J)}$ a Fulton generator of $\cI_u$.
\end{theorem}

\begin{eg} \label{eg:MatrixSchubertS3}
The matrix Schubert variety $M_{s_2}$ is $\{ \Delta_{12, 12} = 0 \}$ or, in other words, $\{ z_{11} z_{22} - z_{12} z_{21} = 0 \}$.
The initial scheme is  $\{ z_{12} z_{21} = 0 \} = \{ z_{12} = 0 \} \cup \{ z_{21} = 0 \}$. 
We depict this as
\[ \In(M_{s_2}) = \begin{sbm} \ast&0&\ast \\ \ast&\ast&\ast \\ \ast&\ast&\ast \\ \end{sbm} \cup  \begin{sbm} \ast&\ast&\ast \\ 0&\ast&\ast \\ \ast&\ast&\ast \end{sbm}.\]
For all other permutations $u \in S_3$, the matrix Schubert variety is a linear space, and is thus equal to its Gr\"obner degeneration. 
(In general, $M_u$ is a linear space if and only if $u$ is a $132$-avoiding permutation, also known as a ``dominant" permutation.)
We give the explicit linear spaces below:
\[ 
M_{e} =  \begin{sbm} \ast&\ast&\ast \\ \ast&\ast&\ast \\ \ast&\ast&\ast \\ \end{sbm} \quad 
M_{s_1} =  \begin{sbm} 0&\ast&\ast \\ \ast&\ast&\ast \\ \ast&\ast&\ast \\ \end{sbm} \quad 
M_{s_1 s_2} =  \begin{sbm} 0&0&\ast \\ \ast&\ast&\ast \\ \ast&\ast&\ast \\ \end{sbm} \quad 
M_{s_2 s_1} =  \begin{sbm} 0&\ast&\ast \\ 0&\ast&\ast \\ \ast&\ast&\ast \\ \end{sbm} \quad 
M_{w_0} =  \begin{sbm} 0&0&\ast \\ 0&\ast&\ast \\ \ast&\ast&\ast \\ \end{sbm} . \]
\end{eg}

Theorem~\ref{thm:KnutsonMiller} states that $\In(\cI_u)$ is an ideal generated by squarefree monomials, a class of ideals known as \newword{Stanley-Reisner ideals}.
The scheme corresponding to a Stanley-Reisner ideal is always a reduced union of coordinate subspaces. Knutson and Miller also describe which coordinate subspaces occur for a given $u$.
Our next task will be to describe this result.
 First, we explain a phenomenon which the reader can already notice in Example~\ref{eg:MatrixSchubertS3}: Only the variables $z_{ij}$ with $i+j \leq n$ occur in the generators of our initial ideals.

\begin{lemma}
 Let $u \in S_n$ and let $\Delta_{I,J}$ be a Fulton generator of $\cI_u$. Let $z_{ij}$ lie on the antidiagonal of $\Delta_{I,J}$. Then $i+j \leq n$.
\end{lemma}

\begin{proof}
Let $I = \{ i_1 < i_2 < \cdots < i_k \}$ and $J = \{ j_1 < j_2 < \cdots < j_k \}$.
Let $i = i_p$, so that $j = j_{k+1-p}$. 
Put $a = i_k$ and $b = j_k$, so we know that $i \leq a-k+p$ and $j \leq b-k + (k+1-p) = b+1-p$.
Since $\Delta_{I,J}$ is a Fulton generator, we have $k > \#(u[a] \cap [b]) \geq a+b-n$.
Combining our inequalities, $i+j \leq (a+p-k) + (b+1-p) = a+b+1-k < (a+b+1) - (a+b-n) = n+1$ so $i+j \leq n$ as desired.
\end{proof}

Let $\fP$ be a subset of $\{ (i,j) \in \ZZ_{> 0}^2 : i+j \leq n \}$, and let $L(\fP)$ be the linear space $\{ z_{ij} = 0 : (i,j) \in \fP \}$.
Define $\bu(\fP)$ to be the Demazure product (see Section~\ref{sec:CombinatorialNotation}) of $s_{i+j-1}$, where the product ranges over $(i,j) \in \fP$, ordered as a substring of 
\[ (n-1,1)\  (n-2,1) (n-2,2) \cdots (2,1) (2,2) \cdots (2,n-2) \ (1,1) (1,2) \cdots (1,n-1) . \]
In other words, in the matrix below, read the entries corresponding to $\fP$, starting at the bottom row and moving to the top and reading each row from left to right:
\[ \begin{bmatrix}
s_1 & s_2 & s_3 & \cdots & s_{n-2} & s_{n-1} & & \\
 s_2 & s_3 & \cdots & s_{n-2} & s_{n-1}& & \\
 s_3 & \cdots & s_{n-2} & s_{n-1}& & \\
 \vdots&\reflectbox{$\ddots$}&\reflectbox{$\ddots$} &&& \\
 s_{n-2} & s_{n-1} &&& & \\
 s_{n-1} &&&& & \\ 
 &&&& & \\ 
 \end{bmatrix}.  \]
This product can be visualized graphically using so-called ``pipe dreams"; see~\cite{KnutsonMiller}.
Knutson and Miller show:
\begin{theorem} \label{thm:KnutsonMillerMore}
Let $u \in S_n$ and let $\fP \subseteq \{ (i,j) \in \ZZ_{\geq 0}^2 : i+j \leq n \}$. Then the linear space $L(\fP)$ is contained in the Gr\"obner degeneration $\In(M_u)$ if and only if $\bu(\fP) \succeq u$.
\end{theorem}

Our focus is not Schubert varieties, but Richardson varieties. 
One can introduce analogous notions of opposite matrix Schubert variety and of matrix Richardson variety, defined as $M^w = \overline{B_+ w B_+}$ and $M_u^w = M_u \cap M^w$ respectively, and one can then consider the analogous Gr\"obner degenerations.
Studying $\In(M^w)$ with respect to an antidiagonal term order is the same as studying $\In(M_u)$ with respect to a diagonal term order. 
Until recently, this was thought to be an intractable problem, because of examples like the following:

\begin{eg}
We consider initial ideals with respect to a diagonal term order.
It is almost never true that the Fulton generators are a Gr\"obner basis. For example, consider $M_{2143} = \{ \Delta_{1,1} = \Delta_{123,123} = 0 \}$. The initial monomial $\In(\Delta_{123,123}) = z_{11} z_{22} z_{33}$ is already divisible by $\In(\Delta_{1,1}) = z_{11}$, and the ideal generated by $\In(\Delta_{123,123})$ and $\In(\Delta_{1,1})$ does not contain $\In(\Delta_{123,123} - \Delta_{23,23} \Delta_1)$, which is $z_{12} z_{21} z_{33}$.
More precisely, the Fulton generators are a Gr\"obner basis if and only if $u$ is 2143-avoiding, also known as ``vexillary"; see~\cite{KMY} for this fact and \cite[Section 9]{Fulton92} for generalities on vexillary permutations.
\end{eg}

\begin{eg}
We consider initial ideals with respect to a diagonal term order.
It can occur that $\In(\cI_u)$ depends on the choice of term order, and that $\In(\cI_u)$ is not reduced.
Both phenomena occur for $u = 214365$; see \cite[Section 7]{KleinWeigandt} for details.
\end{eg}

However, despite these issues, Klein and Weigandt~\cite{KleinWeigandt} have recently succeeded in describing the irreducible components of $\In(M_u)$, and their multiplicities, for a diagonal term order,
using the technology of ``bumpless pipe dreams". 
Equivalently, Klein and Weigandt can describe the irreducible components, and their multiplicities, for $\In(M^w)$ with respect to an anti-diagonal term order.
This raises the natural problem:
\begin{prob}
Describe the irreducible components of $\In(M_u^w)$, and their multiplicities, for an anti-diagonal term order. Since $\In(M_u)$ under anti-diagonal term orders involves ordinary pipe dreams, and $\In(M^w)$ under anti-diagonal term orders involves bumpless pipe dreams, this will presumably require some sort of combination of ordinary and bumpless pipe dreams.
\end{prob}

As this problem is not solved yet, we cannot describe degenerations of matrix Richardson varieties yet.
Therefore, we will be forced to do what is perhaps more natural anyway, to degenerate the Richardson variety $X_u^w$ itself.
Fortunately, our discussion of matrix Schubert varieties will not be wasted; the product $\bu(\fP)$ will reappear and many of our examples will reappear in a new guise.
First, however, we must discuss degenerating the Pl\"ucker algebra $\Pluck$.

\subsection{Gelfand-Tsetlin degeneration of the Pl\"ucker algebra} \label{sec:GTDegenFlag}
The Pl\"ucker coordinates are polynomial functions in the entries of an $n \times n$ matrix.
Thus, the Pl\"ucker algebra is a subalgebra of the polynomial ring $\kappa[z_{ij}]$ in variables $1 \leq i,j \leq n$.
We choose an antidiagonal term order on this polynomial ring,  meaning the leading term of each minor is the product of its antidiagonal entries. 
For any polynomial $f$, we write $\In(f)$ for the leading term of $f$.

We define $\In(\Pluck)$ to be the subalgebra of $\kappa[z_{ij}]$ in variables $1 \leq i,j \leq n$ which is generated (and, in fact, spanned as a $\kappa$-vector space) by $\In(f)$ for $f \in \Pluck$.
For generalities about initial algebras of subalgebras, see~\cite{RS} and~\cite[Chapter 11]{GBCP}.

\begin{theorem} \label{PluckToGTAlg}
There is a flat family over $\AA^1$ whose fiber over $t \neq 0$ is $\Pluck$ and whose fiber over $t=0$ is $\In(\Pluck)$.
The monomials $\In(\Delta(\T))$, where $\T$ ranges over reverse SSYT, form a basis for  $\In(\Pluck)$.
\end{theorem}

Theorem~\ref{PluckToGTAlg} is due to \cite[Theorem 10.6]{GL}; see also~\cite[Theorem 3.2.9]{AlgInvariantTheory} for the analogous result in the Grassmannian case.

\begin{remark} \label{rem:DiagonalVersusAntidiagonal}
If we used a diagonal term order, we would use SSYT instead of reverse SSYT.
Intrinsically, the diagonal case is no harder or easier than the antidiagonal case.
However, we will want to reuse the ideas of pipedreams and the operator $\bu$ from the literature on matrix Schubert varieties, and this literature uniformly uses an antidiagonal term order, for the reasons discussed in Section~\ref{sec:MatrixSchubert}, so we will use one here.
\end{remark}

\begin{eg}
We depict the initial terms of Pl\"ucker monomials by depicting their exponents in an $n \times n$ matrix.
Here are the initial terms of the six nontrivial Pl\"ucker coordinates in $\Fl_3$:
\[ 
\begin{array}{r@{}c@{}l@{\quad}r@{}c@{}l@{\quad}r@{}c@{}l}
\In(\Delta_1) &=& \begin{sbm} 1&0&0 \\ 0&0&0 \\ 0&0&0 \\ \end{sbm} & \In(\Delta_2) &=& \begin{sbm} 0&0&0 \\ 1&0&0 \\ 0&0&0 \\ \end{sbm} &  \In(\Delta_3) &=& \begin{sbm} 0&0&0 \\ 0&0&0 \\ 1&0&0 \\ \end{sbm} \\[0.5 em]
\In(\Delta_{12}) &=& \begin{sbm} 0&1&0 \\ 1&0&0 \\ 0&0&0 \\ \end{sbm} & \In(\Delta_{13}) &=& \begin{sbm} 0&1&0 \\ 0&0&0 \\ 1&0&0 \\ \end{sbm} &  \In(\Delta_{23}) &=& \begin{sbm} 0&0&0 \\ 0&1&0 \\ 1&0&0 \\ \end{sbm} \\
\end{array} 
\]
So we have
\[ \In(\Delta_{32} \Delta_1) = \begin{sbm} 1&0&0 \\ 0&1&0 \\ 1&0&0 \\ \end{sbm} \ \text{and}\  \In(\Delta_{31} \Delta_2) = \In(\Delta_{21} \Delta_3) =  \begin{sbm} 0&1&0 \\ 1&0&0 \\ 1&0&0 \\ \end{sbm}.\]
We note that $\Delta_{32} \Delta_1$ and $\Delta_{31} \Delta_2$ are reverse semistandard, whereas $\Delta_{21} \Delta_3$ is not. 
\end{eg}


\begin{remark}
An excellent reference for this material is \cite[Chapter 14]{MS}. 
However, we warn the reader that their matrices are transpose to ours, because of the issue in Remark~\ref{rem:InverseWarning}.
\end{remark}

It is straightforward to verify:
\begin{prop}
Let $\T$ be a reverse SSYT of shape $\lambda$ and content $\alpha$. Then the column sums of $\In(\Delta(\T))$ are $\lambda$ and the row sums are $\alpha$.
\end{prop}

Thus, $\In(\Pluck)$ is a semigroup ring, where the corresponding semigroup is the semigroup of integer matrices generated by the exponent patterns of the $\In(\Delta_I)$ as above.
We now describe this semigroup explicitly:
\begin{definition}
A \newword{Gelfand-Tsetlin pattern} is an array of integers $g_{ij}$ for $i+j \leq n+1$, obeying the inequalities $g_{ij} \geq g_{(i+1)j} \geq g_{i(j+1)}$.
\end{definition}

\begin{theorem}
Define a linear map $\gamma$ from $n \times n$ integer matrices to $n \times n$ integer matrices by $\gamma(A)_{ij} = A_{ij}+A_{(i+1)j} + \cdots + A_{nj}$. 
Then $A$ is the exponent of some monomial in $\In(\Pluck)$ if and only if $\gamma(A)_{ij}=0$ for $i+j>n+1$ and the $\gamma(A)_{ij}$ for $i+j \leq n+1$ form a Gelfand-Tsetlin pattern.
\end{theorem}

Given $I \subseteq [n]$, we will write $\Gamma_I$ for the Gelfand-Tsetlin pattern corresponding to $I$.
Given a tableaux $\T$, we will write $\Gamma(\T)$ for the Gelfand-Tsetlin pattern corresponding to $\T$.
 
\begin{eg}
The exponent matrices $\In(\Delta_{32} \Delta_1) = \begin{sbm} 1&0&0 \\ 0&1&0 \\ 1&0&0 \\ \end{sbm}$ and $\In(\Delta_{31} \Delta_2)  =  \begin{sbm} 0&1&0 \\ 1&0&0 \\ 1&0&0 \\ \end{sbm}$ correspond to the Gelfand-Tsetlin patterns $\begin{sbm}  2&1&0 \\ 1&1& \\ 1&& \\ \end{sbm}$ and $\begin{sbm} 2&1&0 \\ 2&0 \\ 1&& \\ \end{sbm}$.
\end{eg}


We have defined the map from reverse semistandard Young tableaux to Gelfand-Tsetlin patterns by first going through the intermediate step of initial monomials.
There is also a direct description of the result:
\begin{prop}
Let $\T$ be a reverse SSYT of shape $\lambda$ and let $g= \Gamma(\T)$.. Let $\lambda^k$ be the partition whose Young diagram is those boxes of $\T$ with entries $\geq k$. 
Then the $k$-th row, $(g_{k1}, g_{k2}, \ldots, g_{k(n+1-k)})$ of $g$ is $(\lambda^k_1, \lambda^k_2, \ldots, \lambda^k_{(n+1-k)})$.
\end{prop}

It is also easy to read off the content from a Gelfand-Tsetlin pattern:
\begin{prop} \label{prop:SSYTContent}
Let $\T$ be a reverse SSYT content $\alpha$ and let $g= \Gamma(\T)$. Then the row sums of $g$ are $(\alpha_1+\alpha_2+\cdots+\alpha_n, \alpha_2 + \cdots + \alpha_n, \cdots, \alpha_{n-1} + \alpha_n, \alpha_n)$. 
\end{prop}

For fixed $n$, the polyhedral cone of arrays $g_{ij}$ in $\RR_{\geq 0}^{\binom{n+1}{2}}$ obeying the Gelfand-Tsetlin inequalities is called the \newword{Gelfand-Tsetlin cone}. 
If we fix $\lambda_1 \geq \lambda_2 \geq \cdots \geq \lambda_n$ and consider the slice of the Gelfand-Tsetlin cone where $g_{1j} = \lambda_j$, we get a polytope in $\RR_{\geq 0}^{\binom{n}{2}}$ called the \newword{Gelfand-Tsetlin polytope}. 
We note that the Gelfand-Tsetlin polytope is bounded, since $\lambda_1 = g_{11} \geq g_{ij} \geq g_{1n} = \lambda_n$ for any Gelfand-Tsetlin pattern $g$.

\begin{eg} \label{eg:GTPolytope}
For $n=3$ and $\lambda_1 > \lambda_2 > \lambda_3$, the vertices of the Gelfand-Tsetlin polytope are 
\[ 
\begin{sbm}  \lambda_1 & \lambda_2 & \lambda_3  \\ \lambda_1 & \lambda_2 &  \\ \lambda_1 && \\ \end{sbm} \quad
\begin{sbm}  \lambda_1 & \lambda_2 & \lambda_3 \\ \lambda_1 & \lambda_3 & \\  \lambda_1 && \\ \end{sbm} \quad
\begin{sbm}  \lambda_1 & \lambda_2 & \lambda_3  \\ \lambda_1 & \lambda_2 & \\ \lambda_2 &&\\  \end{sbm} \quad
\begin{sbm}  \lambda_1 & \lambda_2 & \lambda_3 \\ \lambda_2 & \lambda_2 & \\ \lambda_2 && \\  \end{sbm} \quad
\begin{sbm}  \lambda_1 & \lambda_2 & \lambda_3 \\ \lambda_2 & \lambda_3 & \\  \lambda_2 && \\ \end{sbm} \quad
\begin{sbm}  \lambda_1 & \lambda_2 & \lambda_3 \\ \lambda_1 & \lambda_3 & \\  \lambda_3 &&  \\\end{sbm} \quad
\begin{sbm}  \lambda_1 & \lambda_2 & \lambda_3 \\ \lambda_2 & \lambda_3 & \\  \lambda_3 && \\  \end{sbm} .
\]
The figure below depicts this polytope for $(\lambda_1, \lambda_2, \lambda_3) = (4,2,1)$.  
We have drawn this polytope so that the coordinates on the page are the row sums of the Gelfand-Tsetlin patterns (so this is a ``view from infinity").
The polytope has $6$ faces in total, two triangles, two parallelograms and two trapezoids. 
Two of the trapezoids, at the back of the figure, meet along the dashed edge; the other four faces, with solid edges, are in the front of the figure.

\centerline{
\begin{tikzpicture}[scale=1.75]
  \node (v123) at (6,4) {$\begin{sbm} 4 &2 & 1  \\ 4 & 2 &  \\ 4&& \\ \end{sbm}$};
  \node (v132) at (5,4) {$\begin{sbm} 4 &2 & 1  \\ 4 & 1 &  \\ 4&& \\ \end{sbm}$};
  \node (v213) at (6,2) {$\begin{sbm} 4 &2 & 1  \\ 4 & 2 &  \\ 2&& \\ \end{sbm}$};  
  \node (v312) at (3,2) {$\begin{sbm} 4 &2 & 1  \\ 2 & 1 &  \\ 2&& \\ \end{sbm}$};  
  \node (v231) at (5,1) {$\begin{sbm} 4 &2 & 1  \\ 4 & 1 &  \\ 1&& \\ \end{sbm}$};  
  \node (v321) at (3,1) {$\begin{sbm} 4 &2 & 1  \\ 2 & 1 &  \\ 1&& \\ \end{sbm}$};  
  \node (w) at (4,2) {$\begin{sbm} 4 &2 & 1  \\ 2 & 2 &  \\ 2&& \\ \end{sbm}$};  
  \draw[dashed] (v132) -- (v231) ; 
  \draw[fill=white, white] (5,2) circle (0.2) ;
    \draw[fill=white, white] (5,3) circle (0.2) ;
  \draw (v123) -- (v132);
  \draw (v123) -- (v213);
  \draw (v213) -- (v231);
  \draw (v132) -- (v312);
  \draw (v231) -- (v321);
  \draw (v312) -- (v321);  
  \draw (v123) -- (w); \draw (w) -- (v312); 
  \draw (v321) -- (w); \draw (w) -- (v213); 
\end{tikzpicture}
}
\end{eg}

Theorem~\ref{PluckToGTAlg} is an algebraic statement, saying that the homogeneous coordinate ring $\Pluck$ has a flat degeneration to the semigroup ring of the Gelfand-Tsetlin cone.
The corresponding geometric statement is
\begin{Theorem} \label{PluckToGTGeom}
There is a flat proper family over $\AA^1$, whose fibers over $t \neq 0$ are all isomorphic to $\Fl_n$ and whose fiber over $t=0$ is the toric variety of the Gelfand-Tsetlin polytope.
\end{Theorem}
We will call this family the \newword{Gelfand-Tsetlin degeneration of $\Fl_n$}.
As is usual, each face $F$ of the Gelfand-Tsetlin polytope corresponds to a closed toric subvariety of the Gelfand-Tsetlin toric variety.

\subsection{Gelfand-Tsetlin degeneration of the coordinate ring of the Richardson variety} \label{sec:GTDegenRich}
We want to study, not $\Pluck$ but its quotient ring $\Pluck_u^w$. 
Define $\In(\Pluck_u^w)$ to be the quotient of $\In(\Pluck)$ by all monomials of the form $\In(f)$ where $f \in \Pluck$ is a polynomial which is zero in $\Pluck_u^w$. 
The analogs of Theorems~\ref{PluckToGTAlg} and~\ref{PluckToGTGeom} are
\begin{theorem}
There is a flat family of graded algebras over $\AA^1$, whose fibers over $t \neq 0$ are isomorphic to $\Pluck_u^w$ and whose fiber over $t=0$ is isomorphic to $\In(\Pluck_u^w)$.
There is a flat closed subfamily of the Gelfand-Tsetlin degeneration of $\Fl_n$, whose fiber over $t=0$ is the subscheme of the Gelfand-Tsetlin toric variety corresponding to $\In(\Pluck_u^w)$.
\end{theorem}

Before proceeding, we note a notational issue: $\In(\Pluck)$ is the semigroup ring of the semigroup of Gelfand-Tsetlin patterns.
However, Gelfand-Tsetlin patterns are written additively, and the semigroup operation in a semigroup ring is always written multiplicatively.
Therefore, for $g$ a Gelfand-Tsetlin pattern, we will write $\chi^g$ for the corresponding element of the semigroup ring.

\begin{eg} \label{eg:GTFlag3}
Let's consider the Schubert variety $X_{s_1}$ in $\Fl_3$.
The corresponding ideal in $\Pluck$ is $\Delta_1$.
In general, for any principal ideal $\langle f \rangle$, we have $\In(\langle f \rangle)  = \langle \In(f) \rangle$, so $\In(\Pluck_{213})$ is the quotient of $\In(\Pluck)$ by $\In(\Delta_1)$. 
The initial term $\In(\Delta_1)$ is $\begin{sbm} 1&0&0 \\ 0&0&0 \\ 0&0&0 \\ \end{sbm}$, and the corresponding Gelfand-Tsetlin pattern is $\Gamma_1 := \begin{sbm} 1&0&0 \\ 0&0& \\ 0&&\\ \end{sbm}$.
We want to take the semigroup ring of the Gelfand-Tsetlin semigroup and set $\chi^g$ equal to $0$ if $\chi^g$ is divisible by $\chi^{\Gamma_1}$ in the Gelfand-Tsetlin semigroup ring or, equivalently, if $g-\Gamma_1$ is a Gelfand-Tsetlin pattern.

Let $g$ be a Gelfand-Tsetlin pattern. Then the array
\[ g-\Gamma_1 = \begin{bmatrix} g_{11} - 1 &g_{12}&g_{13}\\ g_{21} & g_{22} & \\ g_{31} && \\ \end{bmatrix} \]
will be a Gelfand-Tsetlin pattern if and only if $g_{11} > g_{21}$.
So, in the quotient ring $\In(\Pluck_{s_1})$, we set $\chi^g=0$ for these $g$. In other words, a basis for $\In(\Pluck_{s_1})$ is $\chi^g$ for $g$ a Gelfand-Tsetlin pattern with $g_{11} = g_{21}$. 
We depict this pattern of equalities visually as 
\[ \left[   \!\begin{gathered} \xymatrix@R=0.1 in@C=0.1 in{ \bullet \ar@{=}[d] & \bullet & \bullet \\ \bullet & \bullet & \\ \bullet && \\ } \end{gathered} \right]. \]
The figure below shows the Gelfand-Tsetlin toric variety from example~\ref{eg:GTPolytope}; the face corresponding to this equality (one of the two trapezoids in the rear) is shaded:

\centerline{
\begin{tikzpicture}[scale=1.75]
  \node (v321) at (6,4) {$\begin{sbm} 4 &2 & 1  \\ 4 & 2 &  \\ 4&& \\ \end{sbm}$};
  \node (v231) at (5,4) {$\begin{sbm} 4 &2 & 1  \\ 4 & 1 &  \\ 4&& \\ \end{sbm}$};
  \node (v312) at (6,2) {$\begin{sbm} 4 &2 & 1  \\ 4 & 2 &  \\ 2&& \\ \end{sbm}$};  
  \node (v213) at (3,2) {$\begin{sbm} 4 &2 & 1  \\ 2 & 1 &  \\ 2&& \\ \end{sbm}$};  
  \node (v132) at (5,1) {$\begin{sbm} 4 &2 & 1  \\ 4 & 1 &  \\ 1&& \\ \end{sbm}$};  
  \node (v123) at (3,1) {$\begin{sbm} 4 &2 & 1  \\ 2 & 1 &  \\ 1&& \\ \end{sbm}$};  
  \node (w) at (4,2) {$\begin{sbm} 4 &2 & 1  \\ 2 & 2 &  \\ 2&& \\ \end{sbm}$};  
       \begin{pgfonlayer}{bg} 
  \filldraw[white, fill=gray] (v132.center) -- (v231.center) -- (v321.center) -- (v312.center) -- (v132.center) ;
  \end{pgfonlayer}
  \draw[thick] (v123) -- (v132);
  \draw[thick] (v123) -- (v213);
  \draw[thick] (v213) -- (v231);
  \draw[thick] (v132) -- (v312);
  \draw[thick] (v231) -- (v321);
  \draw[thick] (v312) -- (v321);  
  \draw[thick] (v123) -- (w); \draw[thick] (w) -- (v321);
  \draw[thick] (v312) -- (w); \draw[thick] (w) -- (v213);
\draw[dashed] (v132) -- (5,1.8) ; 
\draw[dashed] (5,2.2) -- (5,2.8) ; 
 \draw[dashed] (5,3.2) -- (v231) ;
\end{tikzpicture}
}
\end{eg}
%
%
In general, $X_u$ will degenerate to more than one face of the Gelfand-Tsetlin toric variety. 
We illustrate this with our next example.
\begin{eg}
Let's consider the Schubert variety $X_{s_2}$ in $\Fl_3$, also known as $X_{132}$. 
The corresponding ideal in $\Pluck$ is $\Delta_{12}$.
So $\In(\Pluck_{s_2}) = \In(\Pluck)/\In(\Delta_{12})$. 
The Gelfand-Tsetlin pattern corresponding to $\Delta_{12}$ is $\Gamma_{12} = \begin{sbm} 1&1& 0 \\ 1&0& \\ 0&& \\ \end{sbm}$. 
So, for a Gelfand-Tsetlin pattern $g$, the array $g - \Gamma_{12}$ will be a Gelfand-Tsetlin pattern if and only if $g_{21} > g_{31}$ \textbf{and} $g_{12} > g_{22}$.
In other words, a basis for $\In(\Pluck_{s_2})$ is $\chi^g$ for $g$ a Gelfand-Tsetlin pattern with $g_{21} = g_{31}$ \textbf{or} $g_{12} = g_{22}$.
We can depict this visually as
\[
 \left[   \!\begin{gathered} \xymatrix@R=0.1 in@C=0.1 in{ \bullet & \bullet  \ar@{=}[d]& \bullet \\ \bullet & \bullet & \\ \bullet && \\ } \end{gathered} \right] \cup
  \left[   \!\begin{gathered} \xymatrix@R=0.1 in@C=0.1 in{ \bullet & \bullet & \bullet \\ \bullet  \ar@{=}[d] & \bullet & \\ \bullet && \\ } \end{gathered} \right] .
  \]
  Again, we depict this as a union of faces of the Gelfand-Tsetlin polytope. This time, two faces are shaded; a triangle and a parallelogram:

\centerline{
\begin{tikzpicture}[scale=1.75]
  \node (v321) at (6,4) {$\begin{sbm} 4 &2 & 1  \\ 4 & 2 &  \\ 4&& \\ \end{sbm}$};
  \node (v231) at (5,4) {$\begin{sbm} 4 &2 & 1  \\ 4 & 1 &  \\ 4&& \\ \end{sbm}$};
  \node (v312) at (6,2) {$\begin{sbm} 4 &2 & 1  \\ 4 & 2 &  \\ 2&& \\ \end{sbm}$};  
  \node (v213) at (3,2) {$\begin{sbm} 4 &2 & 1  \\ 2 & 1 &  \\ 2&& \\ \end{sbm}$};  
  \node (v132) at (5,1) {$\begin{sbm} 4 &2 & 1  \\ 4 & 1 &  \\ 1&& \\ \end{sbm}$};  
  \node (v123) at (3,1) {$\begin{sbm} 4 &2 & 1  \\ 2 & 1 &  \\ 1&& \\ \end{sbm}$};  
  \node (w) at (4,2) {$\begin{sbm} 4 &2 & 1  \\ 2 & 2 &  \\ 2&& \\ \end{sbm}$};  
       \begin{pgfonlayer}{bg} 
  \filldraw[white, fill=gray] (v213.center) -- (v231.center) -- (v321.center) -- (v312.center) -- (v213.center) ;
  \end{pgfonlayer}
  \draw[thick] (v123) -- (v132);
  \draw[thick] (v123) -- (v213);
  \draw[thick] (v213) -- (v231);
  \draw[thick] (v132) -- (v312);
  \draw[thick] (v231) -- (v321);
  \draw[thick] (v312) -- (v321);  
  \draw[thick] (v123) -- (w); \draw[thick] (w) -- (v321);
  \draw[thick] (v312) -- (w); \draw[thick] (w) -- (v213);
\draw[dashed] (v132) -- (5,1.8) ; 
\draw[dashed] (5,2.2) -- (5,2.8) ; 
 \draw[dashed] (5,3.2) -- (v231) ;
\end{tikzpicture}
}
\end{eg}

We now state the general result:
\begin{theorem} \label{thm:PluckDegeneration}
Let $u \in S_n$ and let $g$ be a Gelfand-Tsetlin pattern. Let $\fP$ be $\{ (i,j)  : g_{ij} = g_{(i+1)j} \}$. Then the monomial $\chi^g$ is nonzero in $\In(\Pluck_u)$ if and only if $\bu(\fP) \succeq u$. In other words, $X_u$ degenerates to a reduced union of toric varieties, indexed by those subsets $\fP$ of $\ZZ_{>0}^2$ for which $\bu(\fP) \succeq u$, with $\fP$ corresponding to the toric subvariety corresponding to the face of the Gelfand-Tsetlin polytope where the equalities in $\fP$ hold.
\end{theorem}

We can describe this process visually. Each dashed line in the diagram below represents an equality which may or may not hold in $g$.
\[  \left[   \!\begin{gathered} \xymatrix@R=0.15 in@C=0.15 in{ s_1 \ar@{--}[d] & s_2  \ar@{--}[d]& s_3  \ar@{--}[d] & \bullet  \\ s_2 \ar@{--}[d] & s_3  \ar@{--}[d]  & \bullet & \\  s_3  \ar@{--}[d]  & \bullet&& \\  \bullet&&& \\ } \end{gathered} \right] \]
For each equality which holds, read the simple generators which are at the top of the dashed lines. Read from bottom to top, reading across each row from left to right.
The Demazure product of these generators is $\bu(\fP)$. 

The same combinatorics occurs in Theorem~\ref{thm:PluckDegeneration} and~\ref{thm:KnutsonMillerMore} and there is a geometric reason for this: 
Inside the Gelfand-Tsetlin degeneration of $\Fl_n$, there is a family whose every fiber is $\AA^{\binom{n}{2}}$: Inside $\Fl_n$, this fiber is $\Xo^{w_0} = \{ \Delta_n \Delta_{n(n-1)} \cdots \Delta_{n(n-1) \cdots 2} \neq 0 \}$; inside the Gelfand-Tsetlin toric variety, it is the locus where $\chi^{\Gamma_n} \chi^{\Gamma_{n(n-1)}} \cdots \chi^{\Gamma_{n(n-1) \cdots 2}} \neq 0$. 
In our figures, it is the neighborhood of the upper right vertex. 
Within this open neighborhood, each $X_u \cap \Xo^{w_0}$ degenerates to the Stanley-Reisner scheme given by pipe dreams for $u$. In our figures, this is a union of faces of the Gelfand-Tsetlin polytope incident on the top right vertex.

One would like to deduce Theorem~\ref{thm:PluckDegeneration} from Theorem~\ref{thm:KnutsonMiller}.
However, only one direction of the theorem is clear from this perspective. 
We explain the easy direction.
Let $f \in \Pluck \subset \kappa[z_{ij}]$. 
Let the initial monomial $\In(f)$ be $\prod z_{ij}^{A_{ij}}$.
Since $f \in \Pluck$, we know that $\gamma(A)$ is a Gelfand-Tsetlin pattern; put $\gamma(A) = g$.
It is straightforward to check that $g_{ij} = g_{i(j+1)}$ if and only if $A_{ij}=0$; let $\fP$ be the set of $(i,j)$ for which these equalities hold.
Then Theorem~\ref{thm:KnutsonMiller} tells us that, if $f$ vanishes on $X_u$, then $\bar{u}(\fP) \not\succeq u$.

However, given a Gelfand-Tsetlin pattern $g$ with $\bar{u}(g) \not\succeq u$, it is not clear that it is the initial term of an $f$ which is simultaneously both in $\Pluck$ and also $0$ on $X_u$. 
Kogan and Miller   claimed this result in~\cite{KoganMiller}, but their proof appears to be flawed; it deduces the statement in the open neighborhood of $w_0 B_+$ discussed above, and therefore only proves equality once we localize $\chi^{\Gamma_n} \chi^{\Gamma_{n(n-1)}} \cdots \chi^{\Gamma_{n(n-1) \cdots 2}}$.
For correct proofs, see~\cite{Chirivi}, \cite{Kiritchenko} and~\cite{Kim}; the more recent sources will require less translation to convert into the pipe dream notation which we have used here.

In matrix Schubert varieties, the Gr\"obner degeneration of $M_u$ (with respect to an antidiagonal term order) is reduced, and the Gr\"obner degeneration of $M^w$ (with respect to that same term order) need not be.
In the Pl\"ucker algebra, however, there is a symmetry which allows us to interchange the degeneration of $X_u$ and of $X^{u w_0}$, as we will now explain.

Recall that $\Pluck$ is multigraded by $\ZZ^n$ and is generated by the Pl\"ucker coordinates $\Delta_I$ for $\emptyset \subsetneq I \subseteq [n]$. 
Let $\SLPluck$ be the subring in multidegrees $\ZZ^{n-1} \times \{ 0 \}$; this is the subring generated by  $\Delta_I$ for $\emptyset \subsetneq I \subsetneq [n]$. 
We write $\SLPluck_u$, $\SLPluck^w$ and $\SLPluck_u^w$  for the corresponding quotient rings of $\SLPluck$.

There is a symmetry $\tau: \Delta_I \mapsto \pm \Delta_{[n] \setminus I}$ of $\SLPluck$, where the sign is $(-1)^{\# \{ (i,j) : i \in I,\ j \in [n] \setminus I,\ i<j \}}$.
Geometrically, this derives from the symmetry taking a flag $F_1 \subset F_2 \subset \ldots \subset F_{n-1}$ and sending it to the flag of orthogonal complements $F_{n-1}^{\perp} \subset F_{n-2}^{\perp} \subset \cdots \subset F_1^{\perp}$.
This symmetry takes $\SLPluck_u$ to $\SLPluck^{u w_0}$. 

\begin{eg}
$\SLPluck_{s_1}$ is the quotient of $\SLPluck$ by the ideal $\langle \Delta_1 \rangle$.
We have $\tau(\Delta_1) = \Delta_{23}$, and the quotient of $\SLPluck$ by $\langle \Delta_{23} \rangle$ is $\SLPluck^{s_2 s_1}$
As promised, we have $s_2 s_1 = s_1 w_0$.. 
To give a second example, 
$\SLPluck_{s_1 s_2}$ is $\SLPluck {\big/} \langle \Delta_1, \Delta_{12}, \Delta_{13} \rangle$.
We have $\tau{\big(}  \langle \Delta_1, \Delta_{12}, \Delta_{13} \rangle {\big)} =  \langle \Delta_{23}, \Delta_{3}, \Delta_{2} \rangle$; and $\SLPluck {\big /}   \langle \Delta_{23}, \Delta_{3}, \Delta_{2} \rangle = \SLPluck^{s_2}$. 
Again, we have $s_2 = (s_1 s_2) w_0$.
\end{eg}

For $f \in \SLPluck$, the monomial $\In(f)$ has no $z_{1n}$ factor, and the corresponding Gelfand-Tsetlin pattern has $g_{1n}=0$. 
Define an involution $\tau_0$ on the set of Gelfand-Tsetlin patterns with $g_{1n}=0$ by $\tau_0(g)_{ij} = g_{11} - g_{i(n+2-i-j)}$. 
For $f \in \SLPluck$, if $\In(f)$ corresponds to the Gelfand-Tsetlin pattern $g$, then $\In(\tau(f))$ corresponds to the Gelfand-Tsetlin pattern $\tau_0(g)$ and, thus,
$\In(\Pluck_u)$ is the image of $\In(\Pluck^{w_0 u})$ under the involution of the Gelfand-Tsetlin toric variety induced by $\tau_0$.
See~\cite{Kim} for details.

Therefore, Theorem~\ref{thm:PluckDegeneration} also gives a description of $\In(\Pluck^w)$. 
We now spell out the resulting description.

\begin{theorem}
Let $g$ be a Gelfand-Tsetlin pattern. Let $\fQ = \{ (i,j) : g_{i (n+2-i-j)} = g_{(i+1)(n+1-i-j)} \}$, which we will think of as a pipedream. 
Let $\bw = \bu(\fQ) w_0$.
 Then the monomial $\chi^g$ is nonzero in $\In(\Pluck^w)$ if and only if $\bw \preceq w$.
\end{theorem}

Again, we can depict this graphically. The dashed lines in the diagram below depict equalities which may or may not hold in $g$. 
\[  \left[   \!\begin{gathered} \xymatrix@R=0.15 in@C=0.15 in{ \bullet & s_3  \ar@{--}[dl] & s_2  \ar@{--}[dl] & s_1  \ar@{--}[dl]  \\   \bullet  & s_3  \ar@{--}[dl]  & s_2  \ar@{--}[dl] & \\  \bullet   & s_3  \ar@{--}[dl]&& \\  \bullet &&& \\ } \end{gathered} \right] \]
For each equality which holds, read the simple generators which are at the top of the dashed lines. Read from bottom to top, reading across each row from right to left.
The corresponding Demazure product is $\bu(\fQ)$ and  $\bw = \bu(\fQ) w_0$.

Thus, to each Gelfand-Tsetlin pattern $g$, we have associated two permutations $\bu$ and $\bw$, such that $\chi^g$ is nonzero in the ring $\In(\Pluck_u)$ if and only if $u \preceq \bu$, and such that $\chi^g$ is nonzero in  the ring $\In(\Pluck^w)$ if and only if $\bw \preceq w$.
It is tempting to guess that $\chi^g$ is nonzero in the ring $\In(\Pluck_u^w)$ if and only if $u \preceq \bu$ and  $\bw \preceq w$.
Kim~\cite{Kim} asserts this in a final conclusion, but he miscites \cite{GeomMonTheory} by referencing a result which holds only for $\lambda_1 > \lambda_2 > \cdots > \lambda_n$, not for $\lambda_1 \geq \lambda_2 \geq \cdots \geq \lambda_n$.
However, with that corrected, we do have
\begin{theorem}\label{thm:RichPluckDegeneration}
Let $g$ be a Gelfand-Tsetlin pattern with $g_{11} > g_{12} > \cdots > g_{1n}$. Let $\fP$, $\fQ$, $\bu$ and $\bw$ be defined as above. 
Then $\chi^g$ is nonzero in $\In(\Pluck_u^w)$ if and only if $u \preceq \bu$ and $\bw \preceq w$.
\end{theorem}

\begin{proof}
See \cite{Kim} and adjust the reference to \cite{GeomMonTheory} to use the correct hypothesis.
\end{proof}

If we are in the case that  $\lambda_1 \geq \lambda_2 \geq \cdots \geq \lambda_n$, then the conditions that  $u \preceq \bu$ and $\bw \preceq w$ are not sufficient, as we check with our now familiar example in $R_{1324}^{4231}$, which we began discussing in Example~\ref{eg:s2SquaredInFlags4}.
\begin{eg}
Let $g$ be a Gelfand-Tsetlin pattern. We first work out what it means to say that $1324 \preceq \bu$ and $\bw \preceq 4231$. 
The permutation $\bu(\fP)$ is defined as a Demazure product. Since $1324=s_2$, that Demazure product is $\succeq 1324$ if and only if one of the factors in the Demazure product is $s_2$: In other words, if and only if either $g_{12} = g_{22}$ or $g_{21} = g_{31}$ (or both). 
Similarly, $\bw \preceq 4231$ if and only if either $g_{22} = g_{13}$ or $g_{32} = g_{23}$  (or both).
Thus, in the figure below, we will have $\bu \succeq u$ and $\bw \preceq w$ if and only if at least one of the vertical dashed lines and at least one of the diagonal dashed lines is equality.
\[ \xymatrix@R=0.2in@C=0.2in{
g_{11} & g_{12}  \ar@{--}[d]  & g_{13}  \ar@{--}[dl] & g_{14} \\
g_{21} \ar@{--}[d] & g_{22} & g_{23}  \ar@{--}[dl] & \\
g_{31} & g_{32} && \\
g_{41} &&& \\ } \]

We first consider the $(2,1,2,1)$ weight space in $H^0(R_{1324}^{4231}, \cO(3,2,1,0))$. As $\lambda_1 > \lambda_2 > \lambda_3 > \lambda_4$, this is an example where the criterion of Theorem~\ref{thm:RichPluckDegeneration} applies. The Gelfand-Tsetlin patterns of weight $(2,1,2,1)$  for $\lambda = (3,2,1,0)$ are:
\[ 
\begin{sbm} 3&2&1&0 \\ 2&1&1& \\ 2&1&& \\ 1&&& \\ \end{sbm} \quad
\begin{sbm} 3&2&1&0 \\ 3&1&0& \\ 3&0&& \\ 1&&& \\ \end{sbm} \quad
\begin{sbm} 3&2&1&0 \\ 3&1&0& \\ 2&1&& \\ 1&&& \\ \end{sbm}  \quad. 
\begin{sbm} 3&2&1&0 \\ 2&2&0& \\ 2&1&& \\ 1&&& \\ \end{sbm}\]
The first two are nonzero in $\In(\Pluck_{1324}^{4231})$; the third is zero in $\In(\Pluck_{1324})$ and the fourth is zero in $\In(\Pluck^{4231})$.
So this weight space is two dimensional, as we have computed in Example~\ref{eg:s2SquaredInFlags4Part2}.

Now, however, we repeat this exercise for the $(1,1,1,1)$ weight space in  $H^0(R_{1324}^{4231}, \cO(2,1,1,0))$.
There are now three Gelfand-Tsetlin patterns:
\[  \begin{sbm} 2&1&1&0 \\ 2&1&0& \\ 2&0&& \\ 1&&& \\ \end{sbm} \quad
\begin{sbm} 2&1&1&0 \\ 2&1&0& \\ 1&1&& \\ 1&&& \\ \end{sbm} \quad
\begin{sbm} 2&1&1&0 \\ 1&1&1& \\ 1&1&& \\ 1&&& \\ \end{sbm} . \]
Since all of them have $g_{12} = g_{22} = g_{13}$, all of them have $u \preceq \bu$ and $\bw \preceq w$. 
However, only two of these are nonzero in $\In(\Pluck_{1324}^{4231})$.
Indeed, we have $\Delta_1 \Delta_{234} - \Delta_{2} \Delta_{134}=0$ in the quotient $\Pluck_{1324}^{4231}$.
The leading term of this binomial is $\Delta_{2} \Delta_{134}$, and the corresponding Gelfand-Tsetlin pattern is the middle one, so the middle Gelfand-Tsetlin pattern is $0$ in $\In(\Pluck_{1324}^{4231})$.
\end{eg}

The author is not aware of a source that  answers the following problem, although the answer should be extractable from the discussion of Standard Monomial Theory in Section~\ref{sec:SMT}.

\begin{prob}
Give a combinatorial criterion, solely in terms of the Gelfand-Tsetlin pattern $g$, for whether or not $\chi^g$ is zero in the ring $\In(\Pluck_u^w)$. Is the ring $\In(\Pluck_u^w)$ reduced?
\end{prob}

\subsection{Frobenius splitting and its consequences}  \label{sec:FrobeniusSplitting}

In this section, we will discuss the method of ``Frobenius splitting" and what it implies about the Pl\"ucker algebra and Richardson varieties.
In this section, let our ground field $\kappa$ have characteristic $p>0$.
The standard reference on Frobenius splitting is~\cite{BrionKumar}, which the author recommends for its clarity and thoroughness.
For a shorter introduction from a perspective very close to this chapter, the author also recommends~\cite{FrobSplitKnutson}.

Although we work over a field of characteristic $p$ throughout this section, there are standard ways to transport our major results back to fields of characteristic $0$.
We will point out examples of this as we go along; see Section~1.6 of~\cite{BrionKumar} for more.

Let $A$ be a commutative ring of characteristic $p$.
A map $\phi : A \to A$ is called a \newword{Frobenius splitting} if
\[ 
\phi(x+y) = \phi(x) + \phi(y), \
\phi(x^p y) = x \phi(y) \ \text{and} \
\phi(1) = 1. 
\]
Combining the second and third conditions, we see that $\phi$ is a left inverse to the $p$-th power map, so we can think of $\phi$ as a $p$-th root.

If $S^{-1} A$ is a localization of $A$ and $\phi : A \to A$ is a Frobenius splitting, then $\phi$ naturally extends to a splitting on $S^{-1} A$ by $\phi(s^{-1} a) = s^{-1} \phi(s^{p-1} a)$. 
Thus, a Frobenius splitting on $A$ gives a map of sheaves $\cO_{\Spec A} \to \cO_{\Spec A}$. 
More generally, if $X$ is any scheme over $\ZZ/p \ZZ$ then a Frobenius splitting on $X$ is a map of sheaves $\cO_X \to \cO_X$ obeying the above conditions.

Let $X$ be a  $\ZZ/p \ZZ$-scheme and let $\phi : \cO_X \to \cO_X$ be a Frobenius splitting.
Let $\cI \subset \cO_X$ be an ideal sheaf. Then $\cI$ is called \newword{compatibly split} if $\phi(\cI) \subseteq \cI$; in this case, $\phi$ descend to a Frobenius splitting on the zero scheme of $\cI$.
If $Y$ is a subscheme of $X$, we will say that $Y$ is compatibly split if the ideal sheaf of $Y$ is compatibly split.

Frobenius splittings is very powerful for proving that schemes are reduced.
\begin{prop}
Let $X$ be a scheme with a Frobenius splitting. Then $X$ is reduced. Thus, if $Y$ is a compatibly split subscheme of $X$, then $Y$ is reduced.
\end{prop}

\begin{proof}
Let $U$ be an open set of $X$ and let $z$ be a nilpotent element in $\cO_X(U)$. Choosing $n$ large enough, we have $z^{p^n}=0$. Then $\phi^n(z^{p^n}) = \phi^n(0) = 0$.
But, by the axioms of a Frobenius splitting, $\phi^n(z^{p^n}) = z \phi^n(1) = z$. So $z=0$. We have shown that the only nilpotent element of the coordinate ring is $0$, so $X$ is reduced.
\end{proof}

We will therefore speak, from now on, about compatibly split subvarieties. (We do not take the word ``variety" to imply irreducible.) The collection of compatibly split subvarieties is closed under many basic operations. 

\begin{prop} \label{prop:ComponentSplit}
Let $X$ be a scheme with a Frobenius splitting. Then any irreducible component of $X$ is compatibly split (with its reduced subscheme structure).
\end{prop}

\begin{prop} \label{prop:UnionIntersectionSplit}
Let $X$ be a scheme with a Frobenius splitting and let $Y$ and $Z$ be compatibly split subvarieties. Then $Y \cap Z$ is compatibly split (and, in particular, reduced), and $Y \cup Z$ is compatibly split (with the reduced subscheme structure).
\end{prop}

\begin{proof}[Proofs of Propositions~\ref{prop:ComponentSplit} and~\ref{prop:UnionIntersectionSplit}]
See \cite[Proposition 1.2.1]{BrionKumar} for both these statements.
\end{proof}

These conditions strongly restrict which subvarieties can be compatibly split:

\begin{eg}
There can \textbf{not} be a Frobenius splitting on $\AA^2$ where the lines $x=0$, $y=0$ and $x=y$ are compatibly split.
If they were, then the union of the coordinate axes, with corresponding ideal $\langle xy \rangle$ would be compatibly split.
Then the scheme theoretic intersection of this union with the diagonal would be compatibly split and hence reduced.
This scheme theoretic intersection has ideal $\langle xy \rangle + \langle x-y \rangle = \langle xy, x-y \rangle$.
But $\langle xy, x-y \rangle$ is not a radical ideal; it contains $x^2$ and not $x$.
\end{eg}
In particular, a Frobenius split scheme only can have finitely many split subschemes:
\begin{prop}
Let $X$ be a finite type scheme over a field of characteristic $p$, equipped with a Frobenius splitting. 
Then $X$ only has finitely many compatibly split subschemes.
\end{prop}

\begin{proof}
See~\cite{KumarMehta} or~\cite{Schwede}.
\end{proof}

Frobenius splittings also imply strong consequences for the cohomology of ample line bundles.

\begin{prop} \label{CohomVanishingSplit}
Let $X$ be a scheme with a Frobenius splitting and let $L$ be an ample line bundle on $X$. Then $H^j(X, L)=0$ for all $j > 0$.
If $Y$ is any compatibly split subvariety of $X$, then $H^0(X, L) \to H^0(Y, L)$ is surjective.
\end{prop}

\begin{proof}
See \cite[Theorem 1.2.8]{BrionKumar}.
\end{proof}

The following result is inspired by~\cite[Theorem 1]{GeomMonTheory} and is morally already in that source.

\begin{prop} \label{prop:SplittingBasis}
Let $X$ be a scheme with a Frobenius splitting and let $L$ be an ample line bundle on $X$. 
Then there is a basis $\cB$ for $H^0(X, L)$ such that, for any compatibly split subscheme $Y$ of $X$, the set $\{ b \in \cB : b \in \Ker(H^0(X, L) \to H^0(Y, L)) \}$ is a basis for $\Ker(H^0(X, L) \to H^0(Y, L))$.
In other words, if we take the image of $\cB$ in $H^0(Y, L)$ and delete the zero vectors, we get a basis for $H^0(Y, L)$.
\end{prop}

\begin{proof}
Fix $X$ and $L$ as above; put $V = H^0(X,L)$. For any compatibly split subset $Y$ of $X$, let $K(Y) = \Ker(H^0(X, L) \to H^0(Y, L)) \subseteq H^0(X, L)$. 

We claim that, for compatibly split $Y$ and $Z$, we have $K(Y \cup Z) = K(Y) \cap K(Z)$ and $K(Y \cap Z) = K(Y) + K(Z)$.
The first formula is true by definition.
For the second formula, we clearly have $K(Y \cap Z) \supseteq K(Y) + K(Z)$; we must show that, if $f \in K(Y \cap Z)$, then we can write $f = g+h$ for $g \in K(Y)$ and $h \in K(z)$.
Define a section $\overline{g}$ in $H^0(Y \cup Z, L)$ by $\overline{g}|_Y = 0$, $\overline{g}|_Z = f|_Z$; this is possible since $(f|_Z)_Y = f|_{Y \cap Z} = 0$. 
Using the surjectivity in Proposition~\ref{CohomVanishingSplit}, we lift $\overline{g}$ to $g \in V$. Then $g \in K(Y)$ and $f-g \in K(Z)$. 

The formulas $K(Y \cup Z) = K(Y) \cap K(Z)$ and $K(Y \cap Z) = K(Y) + K(Z)$ show that  $Z \mapsto K(Z)$ is an anti-homomorphism of lattices, from the lattice of compatibly split subvarieties of $X$ to the lattice of subspaces of $V$.
But the operations $\cup$ and $\cap$ distribute over each other, so this means that the subspaces $K(Z)$ form a distributive sublattice of the lattice of subspaces of $V$.
For any distributive  sublattice of the lattice of subspaces, there is such a basis; see~\cite[Chapter 1, Proposition 7.1]{PP}.
\end{proof}

So far, we have describe the general theory of Frobenius splitting without any reference to Richardson varieties.
The connection to Richardson varieties comes through the following theorems:
\begin{theorem} \label{FLnSplit}
There is a Frobenius splitting on $\Fl_n$ such that the compatibly split irreducible subvarieties are precisely the Richardson subvarieties.
\end{theorem}

\begin{proof}
Theorem~2.3.1 of~\cite{BrionKumar} shows that there is a splitting of $\Fl_n$ which splits all Schubert varieties $X_u$ and opposite Schubert varieties $X^w$; therefore, it also splits the Richardson varieties $R_u^w = X_u \cap X^w$. 
For the fact that they are the only split subvarieties, see~\cite{Hague} or~\cite{KLS1}. 
\end{proof}

\begin{theorem} \label{GknSplit}
Let $\Fl_n(k_1, k_2, \ldots, k_r)$ be a partial flag manifold. There is a splitting on $\Fl_n(k_1, k_2, \ldots, k_r)$ such that the compatibly split  irreducible subvarieties are precisely the projected Richardson subvarieties.
\end{theorem}

\begin{proof}
See~\cite{KLS1}.
\end{proof}

Since we have obtained these splittings for any $p$ (in particular, for sufficiently large $p$), standard semi-continuity properties let us deduce the following consequences in characteristic zero as well:
\begin{theorem}
Over any field, any scheme-theoretic intersection of Richardson varieties is reduced, and is a union of Richardson varieties.
Over any field, for any ample line bundle $\cL(\lambda)$ on $\Fl_n$ and any positive dimensional Richardson $R_u^w$, we have $H^i(R_u^w, \cL(\lambda))=0$ for $i>0$, also, the map $H^0(\Fl_n, \cL(\lambda)) \to H^0(R_u^w, \cL(\lambda))$ is surjective.
\end{theorem}

Thus, we can take Richardson subvarieties in $\Fl_n$ and apply the operations of (scheme-theoretic) intersection, of union, and of taking irreducible components, in any order, and never leave the world of Richardson varieties and never obtain a non-radical ideal sheaf. 
And the same holds for projected Richardson varieties in a partial flag manifold.

The consequences of Proposition~\ref{prop:SplittingBasis} are slightly better than we might expect, as we now explain:
\begin{theorem} \label{thm:GeomSMBasis}
Let $\lambda_1 \geq \lambda_2 \geq \cdots \geq \lambda_n$. Then, over any field, there is a basis $\cB$ for $H^0(\Fl_n, \cL(\lambda))$ such that, for each Richardson $R_u^w$, the vectors in $H^0(R_u^w, \cL(\lambda))$ which are nonzero images of vectors from $\cB$ form a basis for $H^0(R_u^w, \cL(\lambda))$.
\end{theorem}

\begin{proof}
We describe the proof over a field of characteristic $p$, and omit the semi-continuity argument.

If $\lambda_1 > \lambda_2 > \cdots > \lambda_n$, then $\cL(\lambda)$ is ample on $\Fl_n$, so this follows immediately from Proposition~\ref{prop:SplittingBasis}.

Now, suppose that we only have $\lambda_1 \geq \lambda_2 \geq \cdots \geq \lambda_n$. Let $k_1$, $k_2$, \dots, $k_r$ be the indices for which $\lambda_{k_i} > \lambda_{k_i +1}$.
Then $\cL(\lambda)$ is pulled back from an ample line bundle on $\Fl_n(k_1, k_2, \ldots, k_r)$.
We abbreviate $\Fl_n(k_1, k_2, \ldots, k_r)$ to $Y$, the map $\Fl_n \to Y$ to $\pi$ and write $\cL$ for the ample line bundle on $Y$ such that $\cL(\lambda) \cong \pi^{\ast}(\cL)$.

Then $\pi^{\ast} : H^0(Y, \cL) \to H^0(\Fl_n, \cL(\lambda))$ and $\pi^{\ast} : H^0(\Pi_u^w, \cL) \to H^0(R_u^w, \cL(\lambda))$ are isomorphisms for each $(u,w)$, fitting into a commuting diagram
\[ \xymatrix{
H^0(\Fl_n, \cL(\lambda)) \ar@{->>}[r]  & H^0(R_u^w, \cL(\lambda))  \\
H^0(Y, \cL) \ar@{->>}[r] \ar[u]^{\cong} & H^0(\Pi_u^w, \cL) \ar[u]^{\cong} \\
} .\]
Since the line bundle $\cL$ is ample on $Y$, Proposition~\ref{prop:SplittingBasis} implies that the desired basis exists in $H^0(Y, \cL)$, and we can then use these diagrams to transport it to a basis in $H^0(\Fl_n, \cL(\lambda))$.
\end{proof}

Theorem~\ref{thm:GeomSMBasis} is a slight improvement on Theorem~1 of~\cite{GeomMonTheory}.
That reference states this result for Richardsons in the partial flag manifold, not for projected Richardsons.
Thus, if we simply want to cite Theorem~1 of~\cite{GeomMonTheory}, we can only deduce Theorem~\ref{thm:GeomSMBasis} for $\lambda_1 > \lambda_2 > \cdots > \lambda_n$.

We discussed splittings on $\Fl_n$ which also split each Richardson $X_u^w$; the reader who prefers commutative algebra to algebraic geometry might prefer to split the ring $\Pluck$ and its quotients $\Pluck_u^w$. 
This is possible:
\begin{prop}
There a Frobenius splitting $\Pluck \to \Pluck$ which descends to each quotient $\Pluck_u^w$.
\end{prop}

\begin{proof}
For any $\lambda$, Lemma~1.1.4 of~\cite{BrionKumar} gives such a splitting on the subring of $\Pluck$ in degrees $\{ k \lambda: k \in \ZZ_{\geq 0} \}$, and it is not hard to see that this splittings are compatible on the whole ring $\Pluck$.
\end{proof}

We describe, without proof, the splittings on $\Fl_n$ and $G(k,n)$ from Theorems~\ref{FLnSplit} and~\ref{GknSplit}.
Let $X$ be a scheme over a field of characteristic $p$. 
A \newword{near splitting} of $X$ is a map $\phi : \cO_X \to \cO_X$ obeying 
\[ \phi(x+y)  = \phi(x)  + \phi(y) \ \text{and} \ \phi(x^p y) = x \phi(y) \]
but not necessarily $\phi(1) = 1$.
The near splittings on $X$ form a $\kappa$-vector space.
If $X$ is regular, then the vector space of near splittings is isomorphic to $H^0(X, \omega_X^{-(p-1)})$ (see~\cite[Section 1.3]{BrionKumar}). 
So we can specify a near-splitting on $X$ by specifying a section $\eta$ of $\omega_X^{-(p-1)}$. In fact, we will specify a section $\sigma$ of $\omega_X^{-1}$ and take $\eta = \sigma^{p-1}$.
Since we will be working on spaces $X$ where we have $H^0(X, \cO^{\times}) = \kappa^{\times}$, a section of $\omega_X^{-1}$ will be determined up to scalar multiple by its vanishing locus. We'll take sections of $\omega_X^{-1}$ which are defined over $\FF_p$; since $a^{p-1} = 1$ for $a \in \FF_p^{\times}$, the resulting section $\sigma^{p-1}$ of $\omega^{-(p-1)}$ will be completely determined by its vanishing locus. 
\begin{theorem}
The splittings on $\Fl_n$ and $G(k,n)$ which compatibly split, respectively, all Richardson varieties and all positroid varieties, correspond to sections of $\omega_X^{-1}$ which vanish on, respectively, the $2(n-1)$ Schubert divisors and the $n$ positroid divisors.
\end{theorem}

This is part of a heuristic for finding splittings in general.
If we want to split a regular variety $X$ compatibly with distinct divisors $D_1$, $D_2$, \dots, $D_r$, then we should choose a section $\sigma$ of $\omega_X^{-1}$ which vanishes on $\bigcup D_i$. 
When we are lucky, as we are in the cases of partial flag varieties, $\sum [D_i]$ will be anticanonical, so there will be a unique such $\sigma$ up to scalar multiple.
One must then check whether $\sigma$ gives a splitting, or only a near splitting. 
See~\cite[Section 1.4]{FrobSplitKnutson} for concrete criteria which can address this question.

\section{The Bott-Samelson varieties and brick varieties}

\subsection{Bott-Samelson varieties} \label{sec:BSVarieties}
Let $s_{i_1}$, $s_{i_2}$, \dots, $s_{i_{a}}$ be a sequence of simple generators for $S_n$. 
We define the \newword{open Bott-Samelson variety} $\BS^{\circ}(i_1, i_2, \ldots, i_{a})$ to be the subvariety of $\Fl_n^{a+1}$ consisting of sequences of flags $(F^0, F^1, F^2, \ldots, F^{a})$ such that $F^0 = e B_+$ and $F^{j-1} \xrightarrow{s_{i_j}} F^j$. (Recall the notation $E_{\bullet} \xrightarrow{\ w\ } F_{\bullet}$ from the end of Section~\ref{BruhatSchubert}.)
Concretely, this means that $F^{j-1}_i = F^j_i$ for $i \neq i_j$ and $F^{j-1}_{i_j} \neq F^j_{i_j}$.
We define the \newword{Bott-Samelson variety} $\BS(i_1, i_2, \ldots, i_{a})$ to be the subvariety of $\Fl_n^{a+1}$ where $F^0 = e B_+$ and, for $1 \leq j \leq a$, we either have  $F^{j-1} \xrightarrow{s_{i_j}} F^j$ or $F^{j-1} = F^j$. In other words, we impose that $F^{j-1}_i = F^j_i$ for $i \neq i_j$  but impose no condition on the relation between $F^{j-1}_{i_j}$ and  $F^j_{i_j}$.

\begin{remark}
Bott-Samelson varieties were introduced by Bott and Samelson \cite{BottSamelson55} as smooth manifolds.
The first papers to consider Bott-Samelson varieties as algebraic varieties were Hansen~\cite{Hansen73} and Demazure~\cite{Demazure74}; Demazure introduced the term ``Bott-Samelson variety".
\end{remark}

There is an elegant way to represent points of a Bott-Samelson variety using wiring diagrams, which was introduced by Magyar~\cite{Magyar98} and which was further developed by Escobar, Pechenik, Tenner and Yong~\cite{EPTY}.
A \newword{wiring diagram} is a collection of $n$ paths $\sigma_1$, $\sigma_2$, \dots, $\sigma_n$ (called \newword{wires}) in $\RR^2$ obeying the following topological conditions:
Each path $\sigma_i$ is the graph of a continuous function $\RR \to \RR$, and we also denote the function by $\sigma_i$.
For each real number $x$, either all the values $(\sigma_1(x), \sigma_2(x), \ldots, \sigma_n(x))$ are distinct, or else precisely two of them are equal, and there are finitely many values $x_1 < x_2 < \cdots < x_a$ where two of the $\sigma_h(x)$ become equal.
If $\sigma_{h_1}(x_j) = \sigma_{h_2}(x_j)$, then we impose that, after possibly switching $h_1$ and $h_2$, we have $\sigma_{h_1}(x) < \sigma_{h_2}(x)$ for $x<x_j$ and $\sigma_{h_1}(x) > \sigma_{h_2}(x)$ for $x>x_j$.

To such a wiring diagram, we associate a word $(s_{i_1}, s_{i_2}, \ldots, s_{i_a})$ in the simple generators of $S_n$. 
Specifically, $i_j$ is the index such that, if $\sigma_{h_1}$ and $\sigma_{h_2}$ cross at $x_j$, then $\sigma_{h_1}(x_j) = \sigma_{h_2}(x_j)$ are the $i_j$-th largest value of the $n-1$ numbers $\{ \sigma_1(x_j), \sigma_2(x_j), \ldots, \sigma_n(x_j) \}$. 
The word $(s_{i_1}, s_{i_2}, \ldots, s_{i_a})$ is reduced if and only if each pair of wires crosses at most once.
We number the wires such that $\sigma_1(x) < \sigma_2(x) < \cdots < \sigma_n(x)$ for $x \ll 0$.
Then, putting $w = s_{i_1} s_{i_2} \cdots s_{i_a}$, we have $\sigma_{w(1)}(x) < \sigma_{w(2)}(x) < \cdots < \sigma_{w(n)}(x)$ for $x \gg 0$.

\begin{eg} \label{eg:wiring}
Here is a wiring diagram for the word $s_2 s_1 s_2 s_3$ in $S_4$: \\ 
\smallskip

\centerline{
\begin{tikzpicture}
\draw (-1,1) -- (1.75,1) -- (2.25,2) -- (2.75,2) -- (3.25,3) -- (3.75,3) -- (4.25,4) -- (6,4) ;
\draw (-1,2) -- (0.75,2) -- (1.25,3) -- (2.75,3) -- (3.25,2) -- (6,2) ;
\draw (-1,3) -- (0.75,3) -- (1.25,2) -- (1.75,2) -- (2.25,1)  -- (6,1) ;
\draw (-1,4) -- (3.75,4) -- (4.25,3)  -- (6,3) ;
\end{tikzpicture}
}
\end{eg}

The regions in the complement of the wiring digaram are called \newword{chambers}.
We will say that a chamber lies at \newword{height} $i$ if $i$ wires pass below it.
For $x \in (x_j, x_{j+1})$, the vertical line $x \times \RR$ crosses through $n-1$ chambers, at heights  $1$, $2$, \dots, $n-1$.
Given a sequence of flags $(F^0, F^1, \dots, F^{a})$ of $\BS(i_1, i_2, \dots, i_{a})$, we write the subspaces $(F^j_1, F^j_2, \dots, F^j_{n-1})$ of the flag $F^j$ in the chambers meeting $(x_j, x_{j+1}) \times \RR$, with $F^j_i$ in the chamber at height $i$. 

Thus, $\BS(i_1, i_2, \dots, i_{a})$ is equated with the space of ways to label the chambers of the wiring diagram with subspaces, so that a chamber at height $i$ is labeled by an $i$-dimensional space and such that, when one chamber is above another, the space in the top chamber contains the point in the bottom chamber.
In the open Bott-Samelson $\BS^{\circ}(i_1, i_2, \dots, i_{a})$, we impose further that the chambers on the two sides of a crossing are not labeled by the same subspace.

\begin{eg} \label{eg:WiringChambers}
We return to the wiring diagram in Example~\ref{eg:wiring} and fill the chambers with subspaces:\\
\smallskip

\centerline{
\begin{tikzpicture}[xscale=1.5, yscale=1.5]
\draw (-1,1) -- (1.75,1) -- (2.25,2) -- (2.75,2) -- (3.25,3) -- (3.75,3) -- (4.25,4) -- (6,4) ;
\draw (-1,2) -- (0.75,2) -- (1.25,3) -- (2.75,3) -- (3.25,2) -- (6,2) ;
\draw (-1,3) -- (0.75,3) -- (1.25,2) -- (1.75,2) -- (2.25,1)  -- (6,1) ;
\draw (-1,4) -- (3.75,4) -- (4.25,3)  -- (6,3) ;
\node (L0) at (0.5,1.5) {$\Span(e_1)$};
\node (L1) at (4,1.5) {$L$};
\node (P0) at (0,2.5) {$\Span(e_1, e_2)$};
\node (P1) at (2,2.5) {$P_1$};
\node (P2) at (4.5,2.5) {$P_2$};
\node (S0) at (1.5,3.5) {$\Span(e_1, e_2, e_3)$};
\node (S1) at (5,3.5) {$H$};
\draw[dashed] (L0) -- (P0); \draw[dashed] (P0) -- (S0);
\draw[dashed] (L0) -- (P1); \draw[dashed] (P1) -- (S0);
\draw[dashed] (L1) -- (P1);
\draw[dashed] (L1) -- (P2); \draw[dashed] (P2) -- (S0); \draw[dashed] (P2) -- (S1);
\end{tikzpicture}
}
\smallskip

The Bott-Samelson $\BS(2,1,2,3)$ corresponds to collections of subspaces $(L, P_1, P_2, H)$, of dimensions $1$, $2$, $2$ and $3$ respectively, obeying the containments shown by the dashed lines.
In the open Bott-Samelson, $\BSo(2,1,2,3)$, we impose additionally that $\Span(e_1) \neq L$, $\Span(e_1, e_2) \neq P_1 \neq P_2$ and $\Span(e_1, e_2, e_3) \neq H$. (However, we do not impose that $\Span(e_1, e_2) \neq P_2$.)
The sequence of flags is
\[ \begin{array}{c@{}c@{\ }c@{\ }c@{}}
F^0 =& {\big(} \Span(e_1) ,& \Span(e_1, e_2) ,& \Span(e_1, e_2, e_3) {\big)} \\
F^1 =& {\big(} \Span(e_1) ,& P_1 ,& \Span(e_1, e_2, e_3) {\big)}  \\
F^2 =&{\big(} L ,& P_1 ,& \Span(e_1, e_2, e_3) {\big)}  \\
F^3 =& {\big(} L ,& P_2 ,& \Span(e_1, e_2, e_3) {\big)}  \\
F^4 =&{\big(} L ,& P_2 ,& H  {\big)}.  \\
\end{array} \]
\end{eg}

We close by describing the geometry of the Bott-Samelson varieties:

\begin{lemma}
The Bott-Samelson variety $\BS(i_1, i_2, \ldots, i_{a})$ is a smooth $a$-dimensional variety and is a repeated $\PP^1$ bundle over a point. The open Bott-Samelson variety $\BS^{\circ}(i_1, i_2, \ldots, i_{a})$ is isomorphic to  $\AA^{a}$.
\end{lemma}

\begin{proof}[Proof sketch]
Consider the map $\BS(i_1, i_2, \ldots,i_{a-1},  i_{a}) \longrightarrow \BS(i_1, i_2, \ldots, i_{a-1})$ which forgets the last flag $F^{a}$. 
All the components of $F^{a}$ are the same as those of $F^{a-1}$ except for $F^{a}_{i_{a}}$. 
The subspace $F^{a}_{i_{a}}$ lies between $F^{a}_{i_{a}-1} = F^{a-1}_{i_{a-1}}$ and $F^{a}_{i_{a}+1} = F^{a-1}_{i_{a+1}}$. The set of subspaces lying between $F^{a-1}_{i_{a-1}}$ and $F^{a-1}_{i_{a+1}}$  is in bijection with the projective line $\PP(F^{a-1}_{i_{a+1}}/F^{a-1}_{i_{a-1}})$. So (ignoring some details),  $\BS(i_1, i_2, \ldots,i_{a-1},  i_{a})$ is a $\PP^1$-bundle over  $\BS(i_1, i_2, \ldots,i_{a-1})$. 
If we work with the open spaces, then the condition $F^{a}_{i_{a}} \neq F^{a-1}_{i_{a}}$ deletes one point from each $\PP^1$ fiber, so $\BS^{\circ}(i_1, i_2, \ldots,i_{a-1},  i_{a})$ is an $\AA^1$-bundle over  $\BS^{\circ}(i_1, i_2, \ldots,i_{a-1})$. 

This shows that $\BS^{\circ}(i_1, i_2, \ldots, i_{a})$  is a repeated $\AA^1$-bundle, but we are supposed to show more strongly that it is isomorphic to $\AA^{a}$.
This follows from Lemma~\ref{lem:TrivialBundle} below.
\end{proof}

\begin{lemma} \label{lem:TrivialBundle}
Let $X_{a} \to X_{a-1} \to \cdots \to X_2 \to X_1 \to X_0$ be a sequence of spaces and maps where $X_0$ is a point and each $X_j \to X_{j-1}$ is either an $\AA^1$ or $\GG_m$-bundle. 
Then all the bundles are trivial, so $X_{a} \cong \AA^k \times \GG_m^{a-k}$, where there are $k$ $\AA^1$-bundles and $a-k$ $\GG_m$-bundles in the sequence..
\end{lemma}

\begin{proof}
Our proof is by induction on $a$, with the base case $a=0$ being obvious. Thus, we need to prove that any $\AA^1$ or $\GG_m$-bundle over $\AA^p \times \GG_m^{q}$ is trivial.

$\GG_m$-bundles over a variety $X$ are classified by $H^1(X, \cO^{\times})$. 
If $X = \Spec A$ for a ring $A$, then this is isomorphic to the Cartier class group of $A$; in particular, if $A$ is a UFD then every $\GG_m$-bundle over $\Spec A$ is trivial.
The space $\AA^p \times \GG_m^{q}$ is the spectrum of a polynomial ring in $p$ ordinary variables and $q$ Laurent variables, which is a UFD, so all $\GG_m$-bundles over $\AA^p \times \GG_m^{q}$ are trivial.

The case of $\AA^1$ is similar but slightly harder. 
$\AA^1$-bundles over $X$ are classified by the non-Abelian cohomology group $H^1(X, \cO^{\times} \ltimes \cO)$, which lies in a long exact sequence between $H^1(X, \cO)$ and $H^1(X, \cO^{\times})$.
If $X$ is affine then $H^1(X, \cO)$ is trivial. If $X$ is the spectrum of a UFD then, as discussed above, $H^1(X, \cO^{\times})$ is trivial.
For the case $X = \AA^p \times \GG_m^{q}$, both hypotheses apply so $H^1(X, \cO^{\times} \ltimes \cO)$ is trivial and the result follows.
 \end{proof}
 
 \subsection{Matrix products formulas for open Bott-Samelson varieties} \label{BSMatrix}

We have shown that the open Bott-Samelson variety $\BSo(i_1, i_2, \ldots, i_a)$ is isomorphic to $\AA^a$, but we will want an explicit isomorphism for computations. 
This computation will use a trick that will also occur many other times, so we set up notation to describe it.
For $1 \leq i \leq n-1$, let $\rho_i : \GL_2 \to \GL_n$ be the group homomorphism 
\[ \rho_i \left( \begin{sbm} a&b \\ c&d \\ \end{sbm} \right) = 
\begin{sbm} 
1&&& && &&\\
&{\tiny \ddots}&& && &&\\
&&1& &&&&\\
&&& a&b &&&\\
&&& c&d &&&\\
&&&&& 1&& \\
&&&&&& {\tiny \ddots}& \\
&&&&&&& 1 \\ 
\end{sbm} \]
where the $2 \times 2$ block is in rows and columns $i$ and $i+1$.

\begin{lemma}
Let $g \in \GL_n$ and let $h \in \GL_2$. Write  $B_+(2)$ for the upper Borel in $\GL_2$. If $h \in B_+(2)$, then $gB_+ = g \rho_i(h) B_+$; if $h \not\in B_+(2)$ then $gB_+ \xrightarrow{s_i} g \rho_i(h) B_+$.
The map $hB_+(2) \mapsto g \rho_i(h) B_+$ is an isomorphism from $\GL_2/B_+(2) \cong \PP^1$ to the set of flags $\{ F : gB_+=F \ \text{or}\ gB_+  \xrightarrow{s_i}  F \}$.
\end{lemma}

\begin{proof}[Proof sketch]
Recall that the $k$-th space in the flag $g B_+$ is the span of the $k$ leftmost columns of $g$. 
Multiplying $g$ on the right by $\rho_i(h)$ acts on the $i$-th and $(i+1)$-th columns of $g$ by column operations, so it preserves all spaces in the flag $gB_+$ except possibly the $i$-th flag; moreover, it preserves the $i$-th flag if and only if $\rho_i(h)$ is a rightward column operation, which happens if and only if $h \in B_+(2)$. We leave the rest to the reader.
\end{proof}

We set $z_i(t) = \rho_i\left( \begin{sbm} t&1 \\ 1&0 \end{sbm} \right)$.

\begin{cor}
Let $g \in \GL_n$. The map $t \mapsto g z_i(t) B_+$ is an isomorphism from $\AA^1$ to the set of flags $F$ with $gB_+  \xrightarrow{s_i}  F$.
\end{cor}

\begin{proof}
$t \mapsto  \begin{sbm} t&1 \\ 1&0 \end{sbm} B_+(2)$ is an isomorphism from $\AA^1$ to the space of non-identity cosets in $\GL_2/B_+(2)$. 
\end{proof}

By induction, we immediately obtain:
\begin{prop}
Let $(s_{i_1}, s_{i_2}, \ldots, s_{i_a})$ be any word in the simple generators of $S_n$. Map $\AA^a$ to $\Fl_n^{a+1}$ by
\[  (t_1, t_2, \ldots, t_a) \mapsto
 {\big(} B_+, z_{i_1}(t_1) B_+,  z_{i_1}(t_1) z_{i_2}(t_2) B_+, \ldots, z_{i_1}(t_1) z_{i_2}(t_2) \cdots z_{i_a}(t_a) B_+ {\big)} . \]
This is an isomorphism $\AA^a \to \BSo(i_1, i_2, \ldots, i_a)$.
\end{prop}


\begin{eg}
Let $n=2$. The Bott-Samelson variety $\BS(1,1,1, \ldots, 1)$ is $(\PP^1)^{a}$. The open Bott-Samelson variety is the space of $(z_1, z_2, \ldots, z_{a})$ in $(\PP^1)^{a}$ where $\Span \begin{sbm} 1\\0 \end{sbm} \neq z_1 \neq z_2 \neq \cdots \neq z_{a}$. 
We can give an explicit isomorphism $\AA^{a} \to \BS^{\circ}(1,1,\ldots,1)$ by sending $(t_1, t_2, \ldots, t_{a})$ to the sequence of points where
\[ z_i = \Span \begin{sbm} t_1 & 1 \\ 1 & 0 \\ \end{sbm} \begin{sbm} t_2 & 1 \\ 1 & 0 \\ \end{sbm}  \cdots \begin{sbm} t_i & 1 \\ 1 & 0 \\ \end{sbm} \begin{sbm} 1\\0 \end{sbm} .\]
\end{eg}

\begin{eg} \label{eg:BS121}
Let $n=3$. The Bott-Samelson variety $\BS(1,2,1)$ is the set of sequences of flags of the form
\[\begin{array}{lcl}
 \Span(e_1) &\subset& \Span(e_1, e_2) \\
   L_1 &\subset& \Span(e_1, e_2) \\
    L_1 &\subset& P_1 \\
     L_2 &\subset& P_1. \\
     \end{array} \]
The open subvariety $\BS^{\circ}(1,2,1)$ imposes that $\Span(e_1) \neq L_1 \neq L_2$ and $\Span(e_1, e_2) \neq P_1$. 
We observe that $L_1$ can be recovered from the flag $(L_2, P_1)$ by $L_1 = P_1 \cap \Span(e_1, e_2)$. 

We can coordinatize $\BS^{\circ}(1,2,1)$ by $\AA^3$ by sending $(t_1, t_2, t_3)$ to
 \[\begin{array}{lcl}
 \Span(e_1) &\subset& \Span(e_1, e_2)\\  \Span(t_1 e_1+ e_2) &\subset& \Span(e_1, e_2)\\ \Span(t_1 e_1+ e_2)  &\subset& \Span(t_1 e_1+ e_2, t_2 e_1 + e_3) \\ \Span((t_1 t_3 + t_2)e_1+ t_3 e_2+e_3)  &\subset& \Span(t_1 e_1+ e_2, t_2 e_1 + e_3) . 
 \end{array}\]
 We rewrite this in terms of matrices; our chain of flags is
 \[
  \begin{bmatrix} 1&0&0 \\ 0&1&0 \\ 0&0&1 \end{bmatrix} B_+,\ \
   \begin{bmatrix} t_1&1&0 \\ 1&0&0 \\ 0&0&1 \end{bmatrix} B_+,\ \
      \begin{bmatrix} t_1&t_2&1 \\ 1&0&0 \\ 0&1&0 \end{bmatrix} B_+,\ \
            \begin{bmatrix} t_1t_3 + t_2&t_1&1 \\ t_3&1&0 \\ 1&0&0 \end{bmatrix} B_+ .
      \]
\end{eg}

In the following section, we will often want to use coset representatives for $\GL_2/B_+(2)$ other than $\begin{sbm} t&1 \\ 1&0 \end{sbm}$.
We therefore adopt the more general convention: For a word $(s_{i_1}, s_{i_2}, \ldots, s_{i_a})$ and a sequence of elements $(h_1, h_2, \ldots, h_a)$ in $\GL_2$ not lying in $B_+(2)$, let $\mu_{i_1 i_2 \ldots i_a}(h_1, h_2, \ldots, h_a)$ be the sequence of flags
\[ \mu_{i_1 i_2 \ldots i_a}(h_1, h_2, \ldots, h_a):=  {\big(} B_+, \rho_{i_1}(h_1) B_+,  \rho_{i_1}(h_1) \rho_{i_2}(h_2) B_+, \ldots, \rho_{i_1}(h_1) \rho_{i_2}(h_2) \cdots \rho_{i_a}(h_a) B_+ {\big)} \]
in $\BSo(i_1, i_2, \ldots, i_a)$.
We close by remarking on some other choices we could use for $h_i$:

\begin{remark} \label{rem:ZvsZdot}
Instead of $\begin{sbm} t&1 \\ 1&0 \end{sbm}$, we could use $\begin{sbm} t&-1 \\ 1&0 \end{sbm}$. This is a bit more natural, because our matrices then lie in $\SL_2$, and we will switch to this convention in Section~\ref{DeodharMatrix}.
For now, however, we make computations easier by omitting the sign.
\end{remark}

\begin{remark}
Instead of $\begin{sbm} t&1 \\ 1&0 \end{sbm}$, we could use $\begin{sbm} 1&0 \\ t&1 \end{sbm}$ for $t \neq 0$. This would parametrize a copy of $\GG_m$ inside the $\AA^1$ of flags $F$ with $g B_+ \xrightarrow{s_i} F$; which $\GG_m$ we obtain will depend on the specific matrix $g$, not only on the flag $g B_+$. We would thus obtain a dense torus $\GG_m^a$ inside $\BSo(i_1, \ldots, i_a)$; this is an example of a ``Deodhar torus", which we will discuss more generally in Section~\ref{sec:Deodhar}.
\end{remark}

\subsection{Maps from Bott-Samelsons to Schubert cells}

We recall the notion of the Demazure product, denoted $\ast$, from Section~\ref{sec:CombinatorialNotation}.
We can use the Demazure product to describe the image of the Bott-Samelson map in $\Fl_n$.

\begin{theorem}
Let $(s_{i_1}, s_{i_2}, \ldots, s_{i_{a}})$ be any word in the simple generators of $S_n$ and let $w = s_{i_1} \ast s_{i_2} \ast \cdots \ast s_{i_{a}}$. Let $\pi : \BS(i_1, i_2, \ldots, i_{a}) \longrightarrow \Fl_n$ be the projection onto the last flag. 
Then the image of $\pi$ is $X^w$. If $(s_{i_1}, s_{i_2}, \ldots, s_{i_{a}})$ is a reduced word for $w$, then the map $\BS^{\circ}(i_1, i_2, \ldots, i_{a}) \to \Xo^w$ is an isomorphism, and $\BS(i_1, i_2, \ldots, i_{a}) \to X^w$ is birational.
\end{theorem}

\begin{proof}
This is due to Demazure. For the statements in the reduced case, see \cite[Th\'eor\`eme 1, Section 3.11]{Demazure74}. The first statement is easily extracted from  \cite[Proposition 4, Section 3.10]{Demazure74}.
\end{proof}

Thus, if we choose any reduced word $(s_{i_1}, s_{i_2}, \ldots, s_{i_{\ell}})$ for $w$, the map $\BS(i_1, i_2, \ldots, i_{\ell}) \to X^w$  is a resolution of singularities of $X^w$.
This map has good cohomological properties:
\begin{theorem}
Let $(s_{i_1}, s_{i_2}, \ldots, s_{i_{\ell}})$ be a reduced word for $w$. We abbreviate $Y = \BS(i_1, i_2, \ldots, i_{\ell})$ and $X= X^w$, so $\pi : Y \to X$ is a resolution of singularieties.
Then $\pi_{\ast} \cO_Y = \cO_X$,  $\pi_{\ast} \omega_Y = \omega_X$ and $(R^j \pi)_{\ast} \cO_Y = (R^j \pi)_{\ast} \omega_Y =0$ for $j>0$.
\end{theorem}

\begin{proof}
See~\cite[Section~3.1]{Andersen} or \cite[Theorem~4]{Ramanathan}.
\end{proof}

\begin{remark}
If  $(s_{i_1}, s_{i_2}, \ldots, s_{i_{a}})$ is a reduced word for $w$ then the composition $\AA^a \to \BSo(i_1, i_2, \ldots, i_a) \to \Xo^w$ is an isomorphism. 
See~\cite{KLR} for a formula to invert this isomorphism.
\end{remark}

We could say much more about Bott-Samelson parametrizations of Schubert varieties, but our goal is to describe Richardsons, not Schuberts, so we move on.

\subsection{Brick varieties and Richardsons}
Let $u \leq w$ be elements of $S_n$. 
In this section, we will describe how to adapt the ideas in the previous section to give parametrizations of the Richardson $R_u^w$ and the open Richarson $\Ro_u^w$.

Let $(s_{i_1}, s_{i_2}, \dots, s_{i_a})$ be a reduced word for $w$ and let $(s_{j_1}, s_{j_2}, \dots, s_{j_b})$ be a reduced word for $u^{-1} w_0$. 
\begin{lemma} \label{lem:ProductW0}
The Demazure product $s_{i_1} \ast s_{i_2} \ast \cdots \ast s_{i_a} \ast s_{j_1} \ast s_{j_2} \ast \cdots \ast s_{j_b}$ is equal to $w_0$.
\end{lemma}

\begin{proof}
Since $u \preceq w$, there must be some subword $(s_{i_{k_1}}, s_{i_{k_2}}, \ldots, s_{i_{k_c}})$ of $(s_{i_1}, s_{i_2}, \dots, s_{i_a})$ which is a reduced word for $u$. Then 
\[ s_{i_1} \ast s_{i_2} \ast \cdots \ast s_{i_a} \ast s_{j_1} \ast s_{j_2} \ast \cdots \ast s_{j_b} \succeq
s_{i_{k_1}} \ast  \cdots \ast s_{i_{k_c}} \ast s_{j_1} \ast s_{j_2} \ast \cdots \ast s_{j_b} = u \ast (u^{-1} w_0) = w_0 . \]
The only element of $S_n$ which is $\succeq w_0$ is $w_0$.
\end{proof}

We will abbreviate the Bott-Samelson $\BS(i_1, i_2, \ldots, i_a, j_1, j_2, \ldots, j_b)$ to $Y$.
From Lemma~\ref{lem:ProductW0}, the map $\pi : Y \to X^{w_0} = \Fl_n$ is surjective. 
We define the \newword{brick variety} $Z$ to be $\pi^{-1}(w_0 B_+)$.
The terminology ``brick variety" was introduced by Escobar~\cite{EscobarBrick}, in reference to earlier related constructions in polyhedral combinatorics~\cite{PilaudSantos, PilaudStump}.

\begin{theorem}
The brick variety is projective and smooth of dimension $\ell(w) - \ell(u)$. Projection onto the $a$-th flag is a birational surjection $\pi : Z \to R_u^w$. 
\end{theorem}

\begin{proof}
The brick variety is obviously a closed subvariety of the Bott-Samelson variety, which is a closed subvariety of a product of flag varieties, so it is projective.
Birationality is also clear: If $F$ is a flag in $\Ro_u^w = \Xo_u \cap \Xo^w$, then there is a unique chain of flags $e B_+ \xrightarrow{s_{i_1}} \cdots \xrightarrow{s_{i_a}} F \xrightarrow{s_{j_1}} \cdots \xrightarrow{s_{j_b}} w_0 B_+$.
For smoothness in characteristic zero, see~\cite[Theorem 3]{EscobarBrick}; the characteristic zero hypothesis is not explicitly stated by Escobar, but Escobar works over $\CC$ throughout and the proof invokes Kleiman's transversality theorem, which only holds in characteristic zero.
For a proof in arbitrary characteristic, see~\cite{Balan} or~\cite[Lemma~A.2]{KLS1}.
\end{proof}

As in the case of Schubert varieties, the map $\pi : Z \to R_u^w$ has good cohomological properties:
\begin{theorem} \label{RichardsonRationalSingularities}
We abbreviate $R_u^w$ to $R$ and continue to use all the other notation above.
We write $\omega_R$ for the dualizing sheaf of $R$ and write $\omega_Z$ for the dualizing sheaf of $Z$; since $Z$ is smooth, the latter is the same as the top wedge power of the cotangent bundle.
With these notations, we have  $\pi_{\ast} \cO_Z = \cO_R$,  $\pi_{\ast} \omega_Z= \omega_R$ and $(R^j \pi)_{\ast} \cO_Z = (R^j \pi)_{\ast} \omega_Z =0$ for $j>0$.
\end{theorem}

\begin{remark}
Let $\kappa$ have characteristic $0$ and let $X$ be a finite type variety over $\kappa$. We say that $X$ has \newword{rational singularities} if there is a proper birational map $f: Y \to X$ from a smooth $Y$ such that $f_{\ast} \cO_Y = \cO_X$ and $R^j f_{\ast} \cO_Y=0$ for $j>0$. If this condition holds for one resolution $f: Y \to X$, then it holds for all resolutions; moreover, in this case, it follows that $X$ is Cohen-Macaulay, that $f_{\ast} \omega_Y = \omega_X$ and that $R^j f_{\ast} \omega_Y=0$. So Theorem~\ref{RichardsonRationalSingularities} states that Richardson varieties have rational singularities in characteristic zero. In finite characteristic, there are several different definitions of ``rational singularities";
Theorem~\ref{RichardsonRationalSingularities} states that Richardson varieties have ``rational resolutions", in the sense of \cite[Definition 3.4.1]{BrionKumar}.
\end{remark}

\begin{proof}
For a full proof in all characteristics, see \cite[Appendix A]{KLS1}.
In characteristic zero, Brion~\cite{BrionPositivity} shows that Richardson varieties have rational singularities (see Lemmas~2 and~3 and the discussion thereafter); tracing through Brion's proof shows that the map from the brick manifold is the particular resolution being constructed.
\end{proof}

The case of open Richardsons is even better. 
We define an open subset $\Zo$ of $Z$ as follows: 
Recall that $Z$ is a closed subset of $\BS(i_1, i_2, \ldots, i_a, j_1, j_2, \ldots, j_b)$ and hence a closed subset of $\Fl^{a+b+1}$; let the flags from this embedding be
$F^0$, $F^1$, \dots, $F^{a+b}$. 
Let $\Zo$ be the open subvariety of $Z$ defined by the conditions $F^0 \xrightarrow{\ w\ } F^a \xrightarrow{u^{-1} w_0} F^{a+b}$. 
Since  $(s_{i_1}, s_{i_2}, \dots, s_{i_a})$ and $(s_{j_1}, s_{j_2}, \dots, s_{j_b})$ are reduced words, these are open conditions.
Moreover, they force $(F^0, F^1, \ldots, F^a, \ldots, F^{a+b})$ to lie in the open Bott-Samelson $\BSo(i_1, i_2, \ldots, i_a, j_1, j_2, \ldots, j_b)$

\begin{theorem} \label{ZoIsRo}
The projection map onto flag $F^a$ gives an isomorphism $\Zo \to \Ro_u^w$. 
\end{theorem}

\begin{proof}[Proof sketch]
For any flag $F \in \Ro_u^w$, we must show that there is a unique chain $(F^0, \ldots, F^a, \ldots, F^{a+b})$ in $\Zo$ with $F = F^a$. 
The uniqueness of $(F^0, F^1, \ldots, F^a)$ follows from $F^0 \xrightarrow{\ w\ } F^a$ and the assumption that $(s_{i_1}, \ldots, s_{i_a})$ is reduced; similarly, the uniqueness of $(F^a, F^{a+1}, \ldots, F^{a+b})$ follows from $F^a \xrightarrow{u^{-1} w_0} F^{a+b}$ and the assumption that $(s_{j_1}, s_{j_2}, \dots, s_{j_b})$  is reduced.
\end{proof}

Let $W$ be the locus in $\BS^{\circ}(i_1, i_2, \ldots, i_a, j_1, j_2, \ldots, j_b)$ where $F^0 \xrightarrow{w_0} F^{a+b}$. 

\begin{prop} \label{AffineOpen}
$W$ is a distinguished open subvariety of $\AA^{a+b}$ and we have $W \cong \AA^{\binom{n}{2}} \times \Ro_u^w$.
\end{prop}

\begin{proof}[Sketch of proof]
We know that $\BS^{\circ}(i_1, i_2, \ldots, i_a, j_1, j_2, \ldots, j_b) \cong \AA^{a+b}$.
The locus $W$ is open in the Bott-Samelson $\BS^{\circ}(i_1, i_2, \ldots, i_a, j_1, j_2, \ldots, j_b)$, since we are imposing an open condition on $(F^0, F^{a+b})$.
Let us see that it is a distinguished open: Write the flag $F^{a+b}$ as $g B_+$, then $(F^0, \ldots, F^{a+b})$ is in $W$ if and only if the bottom left minors of $g$ are nonzero.
So $W$ is given by inverting these minors. We note that, using the explicit matrix parametrizations of $\BS^{\circ}(i_1, i_2, \ldots, i_a, j_1, j_2, \ldots, j_b)$ as a product of $z_i(t)$'s, we can write down the polynomial on $\AA^{a+b}$ which is inverted.

Let $N_+$ be the unipotent radical of $B_+$; in other words, $N_+$ is upper triangular matrices with ones on the diagonal. 
Let $N_+$ act on $\BS^{\circ}(i_1, i_2, \ldots, i_a, j_1, j_2, \ldots, j_b)$ by $g \cdot (F^0, F^1, \ldots, F^{a+b}) = (gF^0, gF^1, \ldots, gF^{a+b})$. 
To check that this action maps the Bott-Samelson variety to itself, note that $g (e B_+) = B_+$, and, if $E \xrightarrow{\ u\ } F$, then $gE \xrightarrow{\ u\ } gF$.
The $N_+$ action takes $W$ to itself.

Note that $\Zo$ is the closed subvariety $F^{a+b} = w_0 B_+$ in $W$.
The group $N_+$ acts freely and transitively on the Schubert cell $X^{w_0}$, so each orbit of $N_+$ on $W$ contains exactly one point of $\Zo$. 
So the map $N_+ \times \Zo \longrightarrow W$ sending $(g, (F^0, F^1, \ldots, F^{a+b}))$ to $(g F^0, g F^1, \ldots, g F^{a+b})$ is an isomorphism.
We deduce that $W \cong N_+ \times \Zo$, by Theorem~\ref{ZoIsRo}, this is isomorphic to $\AA^{\binom{n}{2}} \times \Ro_u^w$.
\end{proof}

\begin{eg}
This example takes place in $\Fl_3$.
Let $u = s_1$ and $w = s_2 s_1$, so $u^{-1} w_0 = s_2 s_1$. We will parametrize $\BSo(2,1,2,1)$ by the product of the matrices $z_i(t)$.
We compute that the final flag in our parametrization is
\[ z_2(t_1) z_1(t_2) z_2(t_3) z_1(t_4)B_+ = 
\begin{bmatrix} 
t_3 + t_2 t_4 & t_2 & 1 \\ 1 + t_1 t_4 & t_1 & 0 \\ t_4 & 1 & 0
\end{bmatrix} B_+. \]
The condition $F^0 \xrightarrow{w_0} F^4$ is then that $t_4$ and $\det \begin{sbm} 1 + t_1 t_4 & t_1 \\  t_4 & 1 \end{sbm}$ are nonzero; the latter is automatic as $\det \begin{sbm} 1 + t_1 t_4 & t_1 \\  t_4 & 1 \end{sbm}=1$. So $\Ro_{213}^{312} \times \AA^3 \cong \{ (t_1, t_2, t_3, t_4) : t_4 \neq 0 \} \cong \GG_m \times \AA^3$. As one might guess, $\Ro_{213}^{312}$ is isomorphic to $\GG_m$ in this case. 

If we want to use this method to describe $\Zo \cong \Ro_u^w$ itself, we have to impose that the final flag is $w_0 B_+$.
This gives the equations $t_3+t_2 t_4 = 1+t_1 t_4 = t_2 = 0$ or, equivalently, $t_2=t_3=0$, $t_4 = -t_1^{-1}$. 
Then the parametrization of $\Ro_u^w$ is $z_2(t_1) z_1(t_2) B_+ = z_2(t_1) z_1(0) B_+$ which gives the parametrization
\[z_2(t_1) z_1(0) B_+ =  \begin{bmatrix}
 0 & 1 & 0 \\
t_1  & 0 & 1 \\
 1 & 0 & 0 \\
 \end{bmatrix} B_+ \quad \text{for}\  t_1 \neq 0 . \]
\end{eg}

Proposition~\ref{AffineOpen} is useful because there are standard algorithms for computing the cohomology and mixed Hodge structure on open sets of affine space, see~\cite{OT, Walther}, and the \texttt{Macaulay 2} command \texttt{deRham} in the \texttt{Dmodules} package.
It also provides a concise proof of the following results:

\begin{cor} \label{OpenRichardsonSmoothIrreducible}
The affine variety $\Ro_u^w$ is smooth and irreducible of dimension $\ell(w) - \ell(u)$, and is the spectrum of a UFD.
\end{cor}

\begin{proof}
From the proposition, $\Ro_u^w \times \AA^{\binom{n}{2}}$ is an distinguished open subvariety of $\AA^{a+b}$. 
Thus, $\Ro_u^w \times \AA^{\binom{n}{2}}$ is smooth of dimension $a+b = \ell(w) + \ell(w_0 u^{-1})$ and is the spectrum of a UFD.
We thus deduce that $\Ro_u^w$ is smooth of dimension $\ell(w) + \ell(w_0 u^{-1}) - \binom{n}{2} = \ell(w) - \ell(u)$ and is the spectrum of a UFD.
(Note that, if $R$ is not a UFD, then $R[t]$ is also not a UFD, so we can cancel the $\AA^{\binom{n}{2}}$ factor.)
\end{proof}

The smoothness was obtained earlier by Richardson~\cite{Richardson92}, who directly checked that the diagonal $E=F$ in $\Fl_n^2$ is transverse to $\Xo^w \times \Xo_u$.
Unique factorization  in the Grassmannian was obtained earlier by Levinson and Purbhoo~\cite{LevinsonPurbhoo}; the flag variety case may be original to this Handbook.

This raises a natural question:
\begin{problem}
Let $u \preceq w$. Is the open Richardson $\Ro_u^w$ always isomorphic to a distinguished open subvariety of $\AA^{\ell(w) - \ell(u)}$?
\end{problem}

\section{The Deodhar decomposition} \label{sec:Deodhar}
We now describe a decomposition of $\Ro_u^w$ into pieces of the form $\GG_m^i \times \GG_a^j$, due to Deodhar~\cite{Deodhar85}.
This decomposition will depend on a choice of reduced word $s_{i_1} s_{i_2} \cdots s_{i_a}$ for $w$.
More generally, for any word $(s_{i_1}, s_{i_2}, \ldots, s_{i_a})$, not necessarily reduced, we will describe a decomposition of the open Bott-Samelson $\BSo(i_1, i_2, \ldots, i_a)$.
When $(s_{i_1}, s_{i_2}, \ldots, s_{i_a})$ is a reduced word for $w$, we have $\BSo(i_1, i_2, \ldots, i_a) \cong \Xo^w$, and our decomposition will refine the Richardson decomposition $\bigsqcup_{u \preceq w} \Ro_u^w$ of $\Xo^w$.

\subsection{The Deodhar pieces}

So, let $(s_{i_1}, s_{i_2}, \ldots, s_{i_a})$ be any word in the simple generators of $S_n$.
Given any point  $(F^0, F^1, \ldots, F^a)$ in $\BSo(i_1, i_2, \ldots, i_a)$, let $v^j$ be the permutation such that $F^j  \in \Xo_{v_j}$. 
For a sequence $(v^0, v^1, \ldots, v^a)$, let $\cD_{\seq}(v^0, v^1, \ldots, v^a)$ be the locally closed subvariety of $\BSo(i_1, i_2, \ldots, i_a)$ corresponding to sequences of flags in these Schubert strata.
(The subscript $\seq$ is in anticipation of a different indexing which is more common.)
If $(s_{i_1}, s_{i_2}, \ldots, s_{i_a})$ is a reduced word for $w$ then, under the identification of $\BSo(i_1, \ldots, i_a)$ with $\Xo^w$, the open Richardson $\Ro_u^w$ is the union of those Deodhar pieces where $v^a = u$. 
Note that we always have $v^0 = e$.

We will call the $\cD_{\seq}(v^0, v^1, \ldots, v^a)$ \newword{Deodhar pieces}.
They are often called ``Deodhar strata", but we are avoiding this term because they do not always form a stratification, see Section~\ref{EG:Dudas}.

Our first task is to identify those sequences $(v^0, v^1, \cdots, v^a)$ which occur for some point in $\BSo(i_1, \ldots, i_a)$.

\begin{lemma} \label{lem:LemDistinguished}
Let $(F^0, F^1, \ldots, F^a)$ be a sequence of flags in $\BSo(i_1, i_2, \ldots, i_a)$ and let $F^j \in \Xo_{v^j}$. 
For each index $1 \leq j \leq a$, we have exactly one of the following three conditions:
\begin{enumerate}
\item $v^j = v^{j-1}$ and $v_j s_{i_j}  \succ v_j$. 
\item $v^j = v^{j-1} s_{i_j} $ and $v_{j-1} \prec v_j$. 
\item $v^j = v^{j-1}  s_{i_j}$ and $v_{j-1} \succ v_j$. 
\end{enumerate}
\end{lemma}

\begin{proof}
We abbreviate $i_j$ to $i$.

The flags $F^{j-1}$ and $F^j$ agree in all subspaces except for the $i$-th subspace; write them as $F^{j-1} = (V_1, V_2, \ldots, V_{i-1}, V, V_{i+1}, \ldots, V_{n-1})$ and $F^j = (V_1, V_2, \ldots, V_{i-1}, V', V_{i+1} \ldots, V_{n-1})$. 
Thus, we either have $v^j = v^{j-1} s_{i} $ or $v^j = v^{j-1}$. What we need to do is to show that, if $v^j = v^{j-1}$, then $v^j s_{i_j} \succ v^j$. 

The space of flags of the form $(V_1, V_2, \ldots, V_{i-1}, X, V_{i+1} \ldots, V_{n-1})$ is a $\PP^1$.
The generic point of this $\PP^1$ is in $\Xo_v$ for some $v$ with $\ell(v s_i)> \ell(v)$, and there is exactly one point which is in $\Xo_{v s_i}$.
Since $F^j \neq F^{j-1}$, it is impossible that $F^j$ and $F^{j-1}$ are both the unique point which is in $\Xo_{v s_i}$ so, if $v^j = v^{j-1}$, it must be the case that $v^j = v^{j-1} = v$.
We then have $v^j s_i = v s_i \succ v = v^{j-1}$, as desired. 
\end{proof}

\begin{remark}
If we stratified the closed Bott-Samelson variety $\BS(i_1, i_2, \ldots, i_a)$ in the analogous manner, the case $v^j = v^{j-1}$ with $v^j s_{i_j} \prec v^j$ could occur as well.
\end{remark}

We define a sequence $(v^0, v^1, \dots, v^j)$ to be a \newword{distinguished sequence of $(s_{i_1}, s_{i_2}, \cdots, s_{i_a})$ for $u$} if $v^0 = e$, $v^a = u$ and, for each index $j$, one of the three conditions in Lemma~\ref{lem:LemDistinguished} applies.
We note that we can encode a distinguished sequence as an subword of  $s_{i_1} s_{i_2} \cdots s_{i_a}$ with product $u$, by recording the positions where $v^j = v^{j-1} s_{i_j}$.
We will replace the omitted letters by the symbol $\bullet$. We will call the resulting word in the alphabet $\{ s_1, s_2, \ldots, s_{n-1}, \bullet \}$ a \newword{distinguished subword}.
The notion of distinguished subwords was introduced by Deodhar~\cite{Deodhar85}.

For a distinguished subword $x$, we will denote the corresponding Deodhar piece by $\cD(x)$. Note the absence of the $\seq$; this is the more common indexing which we anticipated before.
We note that both the notations $\cD$ and $\cD_{\seq}$ are only meaningful in the presence of an understood ambient word. 
We have tried to use the words ``sequence" and ``subword" to help orient the reader as to which indexing convention we are using at any given point.

\begin{eg}
Take $n=3$; let $u=s_1$ and $w = w_0$; we use the word $(s_1, s_2, s_1)$ for $w$. The word $(s_1, s_2, s_1)$ has two subwords with product $s_1$, namely $(\bullet, \bullet, s_1)$ and $(s_1, \bullet, \bullet)$, which correspond to the sequences $(e,e,e,s_1)$ and $(e,s_1, s_1, s_1)$. The first sequence is distinguished but the second is not; since $s_1 s_2 \succ s_1$, we violate the distinguished condition at $j=3$.
So $\Ro_u^w$ has a single Deodhar piece, $\cD_{\seq}(e,e,e,s_1) = \cD(\bullet, \bullet, s_1)$.
\end{eg} 

\begin{eg} \label{eg:Deodhar121}
Take $n=3$; let $u=e$ and $w = w_0$; we use the word $(s_1, s_2, s_1)$ for $w$. There are two subwords of $(s_1, s_2, s_1)$ with product $e$, namely $(\bullet, \bullet, \bullet)$ and $(s_1, \bullet, s_1)$, and both are distinguished.
So $\Ro_e^{w_0}$ is the union of two Deodhar pieces: $\cD(\bullet, \bullet, \bullet) = \cD_{\seq}(e,e,e,e)$ and $\cD(s_1, \bullet, s_1) = \cD_{\seq}(e,s_1, s_1, e)$.

We describe points of $\BS^{\circ}(1,2,1)$ using two lines, $L_1$, $L_2$, and a plane $P_1$, as in Example~\ref{eg:BS121}. 
The Richardson $\Ro_e^{w_0}$ is the open subvariety of $\BS^{\circ}(1,2,1)$ where $L_2$ is transverse to $\Span(e_2, e_3)$ and $P_1$ is transverse to $\Span(e_3)$.
In the coordinates of Example~\ref{eg:BS121}, the open Richardson $\Ro_e^{w_0}$ is the open locus $t_2 (t_1 t_3 + t_2) \neq 0$.

The piece $\cD_{\seq}(e,e,e,e)$ is the piece where $L_1$ is transverse to $\Span(e_2, e_3)$; the piece $\cD_{\seq}(e,s_1,s_1,e)$ is the piece where $L_1 \subset \Span(e_2, e_3)$.
We reinterpret this condition in terms of the flag $(L_2, P_1)$ and in terms of the coordinates $(t_1, t_2, t_3)$.
Since $L_1 = P_1 \cap \Span(e_1, e_2)$, the first piece is the piece where $P_1 \cap \Span(e_1, e_2)$ is transverse to $\Span(e_2, e_3)$; equivalently, the first piece is the piece where $\Span(e_2) \not\subset P_1$ and the second piece is the piece where $\Span(e_2) \subset P_1$.
In terms of the $(t_1, t_2, t_3)$ coordinates, these pieces are $t_1 \neq 0$ and $t_1 = 0$. 
\end{eg}

Given a distinguished sequence $(v^0, v^1, \ldots, v^a)$, we define $J_{=}$, $J_{\uparrow}$ and $J_{\downarrow}$ to be the sets of indices $j$ with $v^j = v^{j-1}$, $v^{j-1} \prec v^j$ and $v^{j-1} \succ v^j$ respectively.
Let  $m_=$, $m_{\uparrow}$ and $m_{\downarrow}$ be the respective cardinalities of these sets.

%

\begin{lemma} \label{lem:PieceStructure}
Let  $(v^0, v^1, \ldots, v^a)$ be a distinguished sequence of $(s_{i_1}, s_{i_2}, \cdots, s_{i_a})$ for $u$.
With the notations $m_{=}$ and $m_{\downarrow}$ as above, we have
$\cD_{\seq}(v^0, v^1, \ldots, v^a) \cong (\GG_m)^{m_{=}} \times \AA^{m_{\downarrow}}$.
\end{lemma}

\begin{proof}[Proof sketch]
Note that, if $(v^0, v^1, \ldots, v^{a-1}, v^a)$ is a distinguished sequence of $(s_{i_1}, s_{i_2}, \cdots, s_{i_{a-1}}, s_{i_a})$, then $(v^0, v^1, \ldots, v^{a-1})$ is a distinguished sequence of $(s_{i_1}, s_{i_2}, \cdots, s_{i_{a-1}})$. 
So deleting the last flag $F^a$ gives a map  $\cD_{\seq}(v^0, v^1, \ldots, v^{a-1}, v^a) \longrightarrow  \cD_{\seq}(v^0, v^1, \ldots, v^{a-1})$.
We abbreviate $\cD_{\seq}(v^0, v^1, \ldots, v^a)$ to $\cD$ and abbreviate $\cD_{\seq}(v^0, v^1, \ldots, v^{a-1})$ to $\cD'$. 
We will sketch a proof that this map is a $\GG_m$-bundle if $j \in J_=$, this map is an isomorphism if $j \in J_{\uparrow}$  and this map is an $\AA^1$-bundle if $j \in J_{\downarrow}$. 
Then, by Lemma~\ref{lem:TrivialBundle}, these bundles are all trivial.

So, let us consider the fibers of the map $\cD \to \cD'$. 
Abbreviate $i_a$ to $i$ and abbreviate $v^{j-1}$ to $v$.
Put $v_+ = v \ast s_i$ and put $v_- = (v \ast s_i) s_i$, so $v_+ = v_- s_i$ and $v_+ \succ v_-$.

Let $(F^0, F^1, \ldots, F^{a-1})$ be a point of $\cD'$ and let $F^{a-1} = (V_1, V_2, \ldots, V_{i-1}, V_i, V_{i+1}, \dots, V_{n-1})$. 
Then $F^a$ is of the form $(V_1, V_2, \ldots, V_{i-1}, X, V_{i+1}, \dots, V_{n-1})$. 
As in the proof of Lemma~\ref{lem:LemDistinguished}, the set of flags of this form is a $\PP^1$.
There is exactly one flag in this pencil, call it $ (V_1, V_2, \ldots, V_{i-1}, W, V_{i+1}, \dots, V_{n-1})$, which is in $\Xo_{v_+}$, and all the other points are in $\Xo_{v_-}$.

If $j \in J_=$, then we must have $v^j = v^{j-1} = v_-$. Then $F^a$ can be any flag where $X \neq V_i$, $W$. So $F^a$ must be chosen from a $\PP^1$ with two points deleted or, in other words, from $\GG_m$.

If $j \in J_{\uparrow}$ then $F^a$ is the unique flag where $X=W$, so the fiber of $\cD \to \cD'$ is a single point.

If $j \in J_{\downarrow}$ then we must have $V_i = W$, and  can be any flag where $X \neq V_i$. 
So $F^a$ must be chosen from a $\PP^1$ with one point deleted or, in other words, from $\AA^1$.
\end{proof}

\begin{remark}
If we were stratifying the closed Bott-Samelson variety instead, the analogous formula would be $\AA^{m_=^+ +m_{\downarrow}}$, where $m_=^+$ would count only the cases with  $v^j = v^{j-1}$ and $v^j s_{i_j} \succ v^j$, not the cases with  $v^j = v^{j-1}$ and $v^j s_{i_j} \prec v^j$.
\end{remark}

We now prove several corollaries of Lemma~\ref{lem:PieceStructure}.
\begin{cor}
Let $\FF_q$ be the finite field with $q$ elements. Then $\#(\Ro_u^w(\FF_q))$ is a polynomial in $q$.
\end{cor}

\begin{proof}
Since $\Ro_u^w$ is the disjoint union of its Deodhar pieces, the number of $\FF_q$ points of $\Ro_u^w$ is the sum, over all distinguished subsequences, of the number of $\FF_q$ points in each Deodhar piece. 
The number of $\FF_q$ points in $\GG_m^{m_{=}} \times \AA^{m_{\downarrow}}$ is $(q-1)^{m_{=}} q^{m_{\downarrow}}$. 
\end{proof}

We can prove something stronger than this: 
\begin{prop}
Continue the notation of Lemma~\ref{lem:PieceStructure}. We have $\dim \Ro_u^w = m_{=} + 2 m_{\downarrow}$.
\end{prop}

\begin{proof}
We have $m_{=} + m_{\downarrow} + m_{\uparrow} = a = \ell(w)$. 
Since $v^0 = e$ and $v^a = u$, we have $m_{\uparrow} - m_{\downarrow} = \ell(u)$.
Subtracting one equation from the other, we obtain $m_{=} + 2 m_{\downarrow} = \ell(w) - \ell(u) = \dim \Ro_u^w$.
\end{proof}

\begin{cor}
Let $u \leq w$ and let $d = \dim \Ro_u^w = \ell(w) - \ell(u)$. Let $\cR_u^w(q)$ be the number of $\FF_q$ points of $\Ro_u^w$. 
Then the polynomial $\cR_u^w$ is palindromic with a sign twist: $\cR_u^w(q) = (-q)^d \cR_u^w(q^{-1})$.
\end{cor}

\begin{proof}
By the previous Proposition, $\cR_u^w(q)$ is a sum of polynomials of the form $(q-1)^{d-2m} q^{m}$, where we have abbreviated $m_{\downarrow}$ to $m$.
We have $(-q)^d (q^{-1}-1)^{d-2m} (q^{-1})^{m} = (q-1)^{2d-m} q^m$, so the sum $f_u^w(q)$ obeys the same relation.
\end{proof}

\begin{remark}
The polynomials $\cR_u^w(q)$ occur in Kazhdan-Lusztig theory, where they are called the \newword{Kazhdan-Lusztig $R$-polynomials}. See~\cite{KL1} and~\cite[Chapter 5]{BB2}.
(These are related to, but much simpler than, the more famous ``Kazhdan-Lusztig polynomials".) 
The relation $\cR_u^w(q) = (-q)^d \cR_u^w(q^{-1})$ can be alternately stated by saying that there is a polynomial $\widetilde{\cR}_u^w$ such that $\cR_u^w(q) = q^{(\ell(w) - \ell(u))/2} \widetilde{\cR}_u^w(q^{1/2} - q^{-1/2})$; the polynomials $\widetilde{\cR}$ are also common in Kazhdan-Lusztig theory.
Deodhar~\cite{Deodhar85} originally introduced his decomposition in order to study the $R$-polynomials.
\end{remark}

\begin{remark}
Work of Shende, Treumann and Zaslow~\cite{STZ} has uncovered relations between the number of $\FF_q$ points in $\Ro_u^w$ and the HOMFLYPT polynomials of knot theory, and more generally between the mixed Hodge structure of $\Xo_u^w$ and Khovanov 
homology.
See Galashin and Lam~\cite{GalashinLamKnot}, Shen and Weng~\cite{ShenWeng} and Trinh~\cite{Trinh} for further developments. \end{remark}

Since the open Richardson $\Ro_u^w$ is irreducible and is the disjoint union of its Deodhar pieces, there must be exactly one Zariski dense Deodhar piece.
This piece must have dimension $\ell(w) - \ell(u)$, so it must have $m_{\downarrow} = 0$. 
The corresponding Deodhar piece is then a torus $\GG_m^{\ell(w) - \ell(u)}$, which we will call the \newword{Deodhar torus}.

We conclude the section by describing the combinatorics of the distinguished sequence, $(v^0, v^1, \cdots, v^a)$, corresponding to this Deodhar torus.
Since $m_{\downarrow}=0$, for every $j$, we must have either
\begin{enumerate}
\item $v^j = v^{j-1}$ and $v^{j} s_{i_j} \succ v^{j}$ \\ or
\item $v^j  = v^{j-1} s_{i_j}$ and $v^{j-1} \prec v^j$.
\end{enumerate}
We call a distinguished sequence of this form \newword{positive}, and we will also use the term \newword{positive} for the corresponding distinguished subword of $(s_{i_1}, s_{i_2}, \ldots, s_{i_a})$. 
So we have shown that, for each $u \leq w$, there is a unique positive subword of $(s_{i_1}, s_{i_2}, \ldots, s_{i_a})$ with product $u$.

We can rewrite the defining conditions of a positive sequence as
\[ v^{j-1} = \begin{cases} v^j & v^j s_{i_j} \succ v^j \\ v^j s_{i_j} & v^j s_{i_j} \prec v^j \\ \end{cases} .\]
We can compute the positive sequence for $u$ by putting $v^a = u$ and using the above equation as a recursion for $v^{j-1}$ in terms of $v^j$.
The positive subword for $u$ can also be described as the rightmost subword of  $(s_{i_1}, s_{i_2}, \ldots, s_{i_a})$ which is a reduced word for $u$.

 \subsection{Matrix product formulas for Deodhar pieces} \label{DeodharMatrix}

We now turn to the problem of parametrizing the Deodhar pieces.
Our primary source is Marsh and Rietsch~\cite{MR}.
We recall the notation $\rho_i : \GL_2 \to \GL_n$ from Section~\ref{BSMatrix}. Define:
\[
\dot{s}_i = \rho_i\left( \begin{sbm} 0&-1 \\ 1&0 \end{sbm}  \right) \qquad
 \dot{z}_i(t) = \rho_i \left( \begin{sbm} t&1 \\ -1&0 \end{sbm} \right) \qquad
  \ddot{z}_i(t) = \rho_i \left( \begin{sbm} t&-1 \\ 1&0 \end{sbm} \right) \qquad
 y_i(t) = \rho_i \left( \begin{sbm} 1&0 \\ t&1 \end{sbm} \right) .
\]
The value of the $-1$'s is in the theory of total positivity, and in the relationship to other Lie types. 
We won't see these advantages here but include them for compatibility with~\cite{MR}; see Theorem~\ref{thm:DeodharPositivity} for a result where the signs matter.

Let $(s_{i_1}, s_{i_2}, \ldots, s_{i_a})$ be a word with Demazure product $w$. 
Let $(v^0, v^1, \ldots, v^a)$ be a distinguished sequence and continue to use the notations $J_{=}$, $J_{\uparrow}$ and $J_{\downarrow}$ as before.
Let $(t_1, t_2, \ldots, t_a)$ be a point in $\AA^a$ with $t_j \in \GG_m$ if $j \in J_{=}$, $t_j =0$ if $j \in J_{\uparrow}$ and $t_j \in \AA^1$ if $j \in J_{\downarrow}$.
So the set of possible values of $t$ forms the space $\GG_m^{m_=} \times \AA^{m_{\downarrow}}$, which we know is isomorphic to the Deodhar piece $\cD_{\seq}(v^0, v^1, \ldots, v^a)$.
Our goal is to give an explicit isomorphism from the set of values of $t$ to  $\cD_{\seq}(v^0, v^1, \ldots, v^a)$.
We will sometimes use the convention of labeling the terms which lie in $\GG_m$ as $t_j$ (for ``torus") and those lying in $\AA^1$ as $u_j$ (for ``unipotent").

Define 
\[ h_j  = \begin{cases}
y_{i_j}(t_j) & j \in J_{=} \\
\dot{s}_{i_j} & j \in J_{\uparrow} \\
\dot{z}_{i_j}(u_j) & j \in J_{\downarrow} \\
\end{cases} . \]
Recall the map $\mu_{i_1 i_2 \cdots i_a}$ from $\left( \GL_2 - B_+(2)\right)^a$ to $\BSo(i_1, i_2, \ldots, i_a)$ introduced in Section~\ref{BSMatrix}.
We will write $g^k$ for the partial product $h_1 h_2 \cdots h_k$, so the image of $\mu_{i_1 \cdots i_a}$ is $(B_+, g^1 B_+, \ldots, g^a B_+)$.

\begin{theorem} \label{thm:MRParametrization}
With the above notation, the map sending $(t_1, t_2, \ldots, t_a)$ to  $(B_+, g^1 B_+, \ldots, g^a B_+)$ in $\BSo(i_1, i_2, \ldots, i_a)$ is an isomorphism from $\GG_m^{m_=} \times \AA^{m_{\downarrow}}$ to the Deodhar piece $\cD_{\seq}(v^0, v^1, \ldots, v^a)$.
\end{theorem}

\begin{eg} \label{eg:121MRFormula}
In $\Fl_3$, consider the word $(s_1, s_2, s_1)$ from Example~\ref{eg:BS121}.
There are two distinguished subexpressions ending in $e$: $(e, e, e, e)$ and $(e, s_1, s_1, e)$; they correspond to $(J_=, J_{\uparrow}, J_{\downarrow}) = ( \{ 1,2,3 \}, \emptyset, \emptyset)$ and $(\{ 2 \}, \{ 1 \}, \{ 3 \})$ respectively. 
Since $(s_1, s_2, s_1)$  is reduced, the projection of $\BSo(1,2,1)$ onto the last flag $\Fl_3$ is an isomorphism with its image, so we focus on describing the final flag $F^3$.

The corresponding matrix products are
\[ \begin{array}{c@{}c@{}ccc}
\begin{bmatrix} 1&0&0 \\ t_1&1&0 \\ 0&0&1 \\ \end{bmatrix} &
\begin{bmatrix} 1&0&0 \\ 0&1&0 \\ 0&t_2&1 \\ \end{bmatrix} &
\begin{bmatrix} 1&0&0 \\ t_3&1&0 \\ 0&0&1 \\ \end{bmatrix} &=& 
\begin{bmatrix} 1 & 0 & 0 \\ t_1+t_3 & 1 & 0 \\ t_2 t_3 & t_2 & 1 \\ \end{bmatrix} \\[0.5 cm]
  \begin{bmatrix} 0&-1&0 \\ 1&0&0 \\ 0&0&1 \\ \end{bmatrix} &
\begin{bmatrix} 1&0&0 \\ 0&1&0 \\ 0&t_2&1 \\ \end{bmatrix} &
\begin{bmatrix} u_3&1&0 \\ -1&0&0 \\ 0&0&1 \\ \end{bmatrix} &=&
\begin{bmatrix}  1 & 0 & 0 \\  u_3 & 1 & 0 \\ -t_2 & 0 & 1 \\ \end{bmatrix} \\
\end{array}
\]
with $t_1$, $t_2$, $t_3 \in \GG_m$ and $u_3 \in \AA^1$.
These two pieces disjointly cover the open Richardson $\Ro^{321}_{123}$, which is $\Delta_1 \Delta_{12} \Delta_{3} \Delta_{23} \neq 0$.
The first piece is the open set $\Delta_{13} \neq 0$, and the second piece is the closed set $\Delta_{13}=0$. 
The reader is invited to compute these minors and see that they are zero or nonzero as appropriate.
\end{eg}

\begin{eg} \label{eg:11111MRFormula}
We give an example with a nonreduced word. We work in $\Fl_2$ with the word $(s_1, s_1, s_1, s_1, s_1)$. 
Then a distinguished sequence is a sequence of six $e$'s and $s_1$'s which starts with $e$ and has no consecutive pair of $s_1$'s.
As a concrete example, we will take the sequence $(e,e,e,s_1,e,s_1)$, corresponding to the subword $(\bullet, \bullet, s_1, s_1, s_1)$.
The matrices $\phi_j$ are
\[ 
\begin{bmatrix} 1&0 \\ t_1 & 1 \end{bmatrix}, \
\begin{bmatrix} 1&0 \\ t_2 & 1 \end{bmatrix} ,\
 \begin{bmatrix} 0&-1 \\ 1 & 0 \end{bmatrix}, \
 \begin{bmatrix} u_4 & 1 \\ -1 & 0 \end{bmatrix}, \ 
 \begin{bmatrix} 0 & -1 \\ 1 & 0 \end{bmatrix} 
 \qquad t_1 \in \GG_m,\ u_3 \in \AA^1
 \]
 The successive partial products $g^j$ are
 \[
  \begin{bmatrix} 1&0 \\ 0&1 \end{bmatrix}, \ 
  \begin{bmatrix} 1&0 \\ t_1&1 \end{bmatrix}, \ 
    \begin{bmatrix} 1&0 \\ t_1+t_2&1 \end{bmatrix}, \ 
  \begin{bmatrix} 0 & -1 \\ 1 & -t_1-t_2 \end{bmatrix}, \
  \begin{bmatrix} 1&0 \\ t_1+t_2+u_4 &1 \end{bmatrix}, \ 
  \begin{bmatrix} 0 & -1 \\ 1 & -t_1-t_2-u_4 \end{bmatrix} . \]
A flag in $\Fl_2$ is simply a point on the projective line $\PP^1$; the sequence of flags $(F^0, F^1, \ldots, F^5)$ in this case is
\[
 \begin{bmatrix} 1 \\ 0 \end{bmatrix}, \ 
  \begin{bmatrix} 1 \\ t_1 \end{bmatrix}, \ 
   \begin{bmatrix} 1 \\ t_1+t_2 \end{bmatrix}, \ 
    \begin{bmatrix} 0 \\ 1 \end{bmatrix}, \ 
       \begin{bmatrix} 1 \\ t_1+t_2+u_4 \end{bmatrix}, \ 
    \begin{bmatrix} 0 \\ 1 \end{bmatrix} . \]       
    Note that consecutive elements of this sequence are always distinct points of $\PP^1$; this is the Bott-Samelson condition.
    Note also that $F^0$, $F^1$, $F^2$ and $F^4$ are in the Schubert cell $\Xo_e = \{ \Delta_1 \neq 0 \}$ where as $F^3$ and $F^5$ are in the Schubert cell $\Xo_{s_1} = \{ \Delta_1 = 0 \}$; this is the additional Deodhar condition.
\end{eg}

To prove Theorem~\ref{thm:MRParametrization}, we want the following lemma:
\begin{lemma} \label{lem:KeyDeodharLemma}
Let $(s_{i_1}, s_{i_2}, \ldots, s_{i_a})$ be a word in the simple generators of $S_n$, and let $(v^0, v^1, \ldots, v^a)$ be a distinguished sequence.
Let $(g^0, g^1, \ldots, g^a)$ be as above. Then $g^j \in B_- v^j$.

More precisely, let $\dot{v}^j$ be the signed permutation matrix whose nonzero entries are in positions $(v(k), k)$ and where the entries are $\pm 1$ with the signs chosen such that the nonzero left-justified minors of $\dot{v}^j$ are equal to $1$. Then $g^j \in N_- \dot{v}^j$.
\end{lemma}

The signed permutation matrix $\dot{v}$ has the important property that $\dot{(v s_i)} = \dot{v} \dot{s}_i$ if $v s_i \succ v$. 
We will not need to get the signs right if our only goal is to prove Theorem~\ref{thm:MRParametrization}, but it is nicer to get them right now than to put them in later.
The signs can be seen concretely in Examples~\ref{eg:121MRFormula} and~\ref{eg:11111MRFormula} by noting that in each product, if we take the topmost nonzero minor in columns $1$ through $k$, that minor is equal to $1$.

\begin{proof}
Our proof is by induction on $j$. The base case, $j=0$, is clear, since $g^0 = v^0 = e$, and $e \in N_- e$.

For $j \geq 1$, we break into cases according to whether $j$ is in $J_=$, $J_{\uparrow}$ or $J_{\downarrow}$.
We abbreviate $i_j$ to $i$ and $v^{j-1}$ to $v$.

\textbf{Case 1: $j \in J_{=}$:} Since $j \in J_{=}$, we have $v^j = v^{j-1} = v$ and $v s_i \succ v$.
Assume, inductively, that $g^{j-1} = b \dot{v}$ for some $b \in N_-$. 
Then $g^j = b \dot{v} y_i(t) = b (\dot{v} y_i(t) \dot{v}^{-1}) \dot{v}$. The condition that  $v s_i \succ v$ is equivalent to $v(i) < v(i+1)$, which implies that $\dot{v} y_i(t) \dot{v}^{-1} \in N_-$. 
So $b (\dot{v} y_i(t) \dot{v}^{-1}) \in N_-$ as desired.

\textbf{Case 2: $j \in J_{\uparrow}$:} Since $j \in J_{\uparrow}$, we have $v^j = v s_i \succ v$,  and so $\dot{v}^j = \dot{v}^{j-1} \dot{s}_i$.
Assume, inductively, that $g^{j-1} = b \dot{v}$ for some $b \in N_-$. 
Then $g^j = b \dot{v} \dot{s}_i = b \dot{v}^j$, as desired.

\textbf{Case 3: $j \in J_{\downarrow}$:} Since $j \in J_{\downarrow}$, we have $v^j = v s_i \prec v$, and so $\dot{v}^j = \dot{v}^{j-1} \dot{s}_i^{-1}$.
Assume, inductively, that $g^{j-1} = b \dot{v}$ for some $b \in N_-$. 
Then $g^j = b \dot{v} \dot{z}_i(u) = b (\dot{v} \dot{z}_i(u) \dot{s}_i \dot{v}^{-1}) (\dot{v} s_i^{-1})$. The matrix $\dot{z}_i(u) \dot{s}_i$ has a $2 \times 2$ block of the form $\begin{sbm} 1 & u \\ 0 &1 \end{sbm}$ in rows and columns $i$ and $i+1$, and is the otherwise the identity.
The condition that  $v s_i \succ v$ is equivalent to $v(i) > v(i+1)$, which implies that $\dot{v} (\dot{z}_i(u) \dot{s}_i) \dot{v}^{-1} \in N_-$. So $b (\dot{v} \dot{z}_i(u) \dot{s}_i \dot{v}^{-1}) \in N_-$ as required.
\end{proof}

\begin{proof}[Sketch of proof of Theorem~\ref{thm:MRParametrization}]
Since $g^j$ and $g^{j-1}$ differ by right multiplication by a matrix  which is the identity except in rows and columns $i_j$ and $i_{j+1}$, the flags $F^{j-1}$ and $F^j$ differ only in the $i_j$-th subspace, and the two flags do differ in that subspace. So $F^{j-1} \xrightarrow{s_{i_j}} F^j$ and thus our sequence of flags lies in $\BSo(i_1, i_2, \ldots, i_a)$. 

Moreover, let $gB_+$ be any flag. 
The space of flags $F$ with $g B_+  \xrightarrow{s_i} F$ is an affine line.
The map $u \mapsto g \dot{z}_i(u) B_+$ is an isomorphism from $\AA^1$ to that affine line.
The flag $g \dot{s}_i B_+$ is one point on that line, and  $t \mapsto g y_i(t) B_+$ is an isomorphism from $\GG_m$ to the complement of that point. 
So $\BSo(i_1, i_2, \ldots, i_a)$ is stratified into the images of the maps corresponding to the various distinguished subsequences, and each map is an isomorphism onto its image. 

What remains to be checked is that the image of the map corresponding to distinguished subsequence $(v^0, v^1, \ldots, v^a)$ is the Deodhar stratum $\cD(v^0, v^1, v^2, \ldots, v^a)$.
Let's spell out what this means.
Truncating the distinguished subsequence to $(v^0, v^1, \ldots, v^{j-1})$ gives a distinguished subsequence for $(s_{i_1}, s_{i_2}, \ldots, s_{i_{j-1}})$; assume inductively that the image of our map is the Deodhar stratum for this truncated sequence.

If $j \in J_{\downarrow}$, then  $\cD(v^0, v^1, \ldots, v^{j-1}, v^j)$ is simply the set of all sequences of flags $(F^0, F^1, \ldots, F^{j-1}, F^j)$ where $(F^0, F^1, \ldots, F^{j-1}) \in \cD(v^0, v^1, \ldots, v^{j-1})$ and $F^{j-1} \xrightarrow{s_{i_j}} F^j$, and the image of our map is also this set of sequences of flags.

The more interesting case is where $j \in J_=$ or $J_{\uparrow}$. In this case, we can consider the set of sequences of flags $(F^0, F^1, \ldots, F^{j-1}, F^j)$ where $(F^0, F^1, \ldots, F^{j-1}) \in \cD(v^0, v^1, \ldots, v^{j-1})$ and $F^{j-1} \xrightarrow{s_{i_j}} F^j$; this is an $\AA^1$ bundle over $\cD(v^0, v^1, \ldots, v^{j-1})$. 
In each $\AA^1$ fiber, there is one point which is in $\cD(v^0, v^1, \ldots, v^{j-1}, v^{j-1} s_{i_j})$, which is distinguished by the property that $F^j \in \Xo_{v_{j-1} s_{i_j}}$, and all the other points of the $\AA^1$ fiber are in $\cD(v^0, v^1, \ldots, v^{j-1}, v^{j-1})$. There is also one point of the $\AA^1$ fiber which is in the image of our map for $(v^0, v^1, \ldots, v^{j-1}, v^{j-1} s_{i_j})$ and the other points are in the image of our map for $(v^0, v^1, \ldots, v^{j-1}, v^{j-1})$. 
In other words, we need to check that our $g^j B_+ \in \Xo_{v_j}$ and then the rest will follow.

We have $\Xo_{v_j}  = (B_- v^j B_+)/B_+$. 
By Lemma~\ref{lem:KeyDeodharLemma}, we have $g^j \in N_- \dot{v}^j \subset B_- v^j$, so we have $g^j B_+ \in (B_- v^j B_+)/B_+$ as desired.
\end{proof}

Inverting the isomorphism $\GG_m^{m_=} \times \AA^{m_{\downarrow}} \longrightarrow \cD(v^0, v^1, \ldots, v^a)$ is quite complex; see~\cite{MR} for the general formula.
We will describe the result for the Deodhar torus in the case where $(s_{i_1}, \ldots, s_{i_a})$ is reduced (which is the case which is relevant to Richardsons).
We first set up some auxilliary functions, called \newword{chamber minors}.

Let $(v^0, v^1, \ldots, v^a)$ be the positive sequence for $u$; it will also be convenient to put $w^j = s_{i_1} s_{i_2} \cdots s_{i_j}$. 
Because $(s_{i_1}, \ldots, s_{i_a})$ is reduced, at any point $(F^0, \ldots, F^a)$ of $\BSo(i_1, \ldots, i_a)$, the flag $F^j$ is in $\Xo^{w^j}$. By the definition of the Deodhar piece, if $(F^0, \ldots, F^a)$ is in $\cD_{\seq}(v^0, \ldots, v^a)$, then $F^j \in \Xo_{v^j}$. So, combining these, $F^j \in \Ro_{v^j}^{w^j}$. 
This means that the $k$-th subspace, $F^j_k$, is in the Grasmmannian Richardson variety $\Ro_{v^j[k]}^{w^j[k]}$. In particular, the Pl\"ucker coordinates $\Delta_{v^j[k]}(F^j_k)$ and $\Delta_{w^j[k]}(F^j_k)$ are nonzero.
We define the ratio 
\[ \Phi^j_k = \frac{\Delta_{v^j[k]}(F^j_k)}{\Delta_{w^j[k]}(F^j_k)} \]
to be the $(j,k)$-\newword{chamber minor}.
Since this is a ratio of two Pl\"ucker coordinates, it is a well defined invariant of the subspace $F^j_k$.
We can visualize the chamber minors as written in the chambers of the wiring diagram; the chamber minor in a given chamber is the ratio of the first and last nonzero Pl\"ucker coordinate of the subspace in that chamber.

Then Marsh and Rietsch's formula is the following:
\begin{theorem} \label{MRFormula}
With the above notation, let $j \in J_=$, so that $v^{j-1} = v^j$, and set $k = i_j$. Then
\[ t_k = \frac{\Phi^j_{k+1} \Phi^j_{k-1}}{\Phi^{j-1}_k \Phi^j_k}. \]
\end{theorem}
We remark that $F^{j-1}_{k \pm 1} = F^j_{k \pm 1}$, so we could switch the superscripts in the numerator to $j-1$ without effecting the formula.
Visually, these are the minors for the four chambers surrounding the $j$-th crossing of the wiring diagram.

\begin{remark}
Suppose we parametrize the open Bott-Samelson by $\mu_{i_1 \cdots i_a}\left( \begin{sbm} t_1&-1 \\ 1&0 \end{sbm},  \begin{sbm} t_2&-1 \\ 1&0 \end{sbm}, \cdots,  \begin{sbm} t_a&-1 \\ 1&0 \end{sbm} \right)$, and write $g^j := \ddot{z}_{i_1}(t_1) \cdots \ddot{z}_{i_j}(t_j)$. Then the minor $\Delta_{w^j[k]}(g^j)$ will always be $1$, so the formula for the chamber minor simplifies to  $\Delta_{v^j[k]}(g^j)$. This may appease the reader who wonders why we are using the term ``chamber minor" for a ratio of two minors.
\end{remark}

\begin{remark}
On the other hand, we can parametrize the Deodhar piece using using $\dot{s}_i$ and $y_i$, as we have described in this section. In that case, Lemma~\ref{lem:KeyDeodharLemma}, shows that the minor $\Delta_{v^j[k]}(g^j)$ is $1$ so, in this case, the formula for the chamber minor simplifies to $\Delta_{w^j[k]}(g^j)^{-1}$. This may irritate the reader who doesn't think we should use the term ``chamber minor" for the reciprocal of a minor.
\end{remark}

Formula~\ref{MRFormula} tells us how to compute the $t_k$ parameters given the full list of flags $(F^0, F^1, \ldots, F^a)$, but we often want a formula in terms of just $F^a$. 
For a particular word known as the unipeak word, the flag $F^j$ is easily recovered from $F^a$; see Section~\ref{sec:unipeak}. 
For the general case, see~\cite{MR}.

Recently, Galashin, Lam, Sherman-Bennett and Speyer~\cite{GLSS} constructed a cluster structure on $\Ro_u^w$. The cluster variables are certain monomials in the chamber minors.

\subsection{The Deodhar pieces do not form a stratification} \label{EG:Dudas}

In this section, we will verify that the Deodhar pieces do not, in general, form a stratification, meaning that the closure of a Deodhar piece is not, in general, a union of Deodhar pieces.
Dudas~\cite{Dudas} demonstrated this earlier in Lie type $B$.
Our examples were found by unfolding Dudas's example from type $B_n$ to type $A_{2n-1}$, and then discovering a number of simplifications.
As far as we know, this Handbook is the first source to present such examples in Lie type $A$.

Recall that Bott-Samelson varieties for $\Fl_n$ are indexed by words in the alphabet $\{1,2,\ldots, n-1 \}$, and that Deodhar pieces are indexed by certain words in the alphabet $\{ 1,2,\ldots, n-1, \bullet \}$.
Given a word $x$, we denote its reversal by $\rev{x}$, and we denote the concatenation of words $x$ and $y$ by $xy$.

It is easier to give a counterexample in an open Bott-Samelson variety, without the assumption that $(s_{i_1}, s_{i_2}, \ldots, s_{i_a})$ is reduced. 
 Let $x$ be a word in $\{1,2,\ldots, n-1 \}$ and let $a$ and $b$ be distinct distinguished subwords of $x$. 
Then $a \rev{a}$ and $b \rev{b}$ will be distinguished subwords of $x \rev{x}$ (for the identity).
Suppose that 
\begin{enumerate}
\item $\overline{\cD(a)} \supset \cD(b)$ in $\BSo(x)$ but
\item $\dim \cD(a \rev{a}) \leq \dim \cD(b \rev{b})$ in $\BSo(x \rev{x})$.
\end{enumerate}

A concrete example is $x = (1,2,1)$, $a = (1,\bullet, 1)$, $b = (1,2,\bullet)$.
To check the first condition, we describe points of $\BS^{\circ}(1,2,1)$ using two lines, $L_1$, $L_2$, and a plane $P_1$, as in Example~\ref{eg:BS121}. Then $\cD(1,\bullet,1)$ is the subvariety where $L_1 = \Span(e_2)$ and $P_1 \neq \Span(e_2, e_3)$, and that $\cD(1,2,\bullet)$ is the subvariety where $L_1 = \Span(e_2)$, $P_1=\Span(e_2, e_3)$ and $L_2 \neq e_3$. It is easy to see from this description that $\overline{\cD(1,\bullet,1)}$ contains $\cD(1,2,\bullet)$. (Concretely, $\overline{\cD(1,\bullet,1)}$ is the locus where $L_1 = \Span(e_2)$.)
To check the second condition, note that $\cD(1, \bullet, 1, 1 \bullet, 1) \cong \cD(1,2,\bullet, \bullet, 2,1) \cong \GG_m^2 \times \GG_a^2$. 

\begin{prop}
For $x$, $a$ and $b$ as above, we have $\overline{\cD(a \rev{a})} \cap \cD(b \rev{b}) \neq \emptyset$, but $\overline{\cD(a \rev{a})} \not\supseteq \cD(b \rev{b})$.
Thus, the Deodhar pieces do not form a stratification of $\BSo(x \rev{x})$.
\end{prop}

\begin{proof}
Let $r$ be the length of the word $x$. Inside $\BSo(x \rev{x})$, let $M$ be the subvariety of sequences of flags of the form $(F^0, F^1, \ldots, F^{r-1}, F^r, F^{r-1}, \ldots, F^1, F^0)$. Then projection onto the first $r+1$ flags is an isomorphism $M \cong \BSo(x)$, and this isomorphism takes $\cD(a \rev{a}) \cap M$ to $\cD(a)$, and $\cD(b \rev{b}) \cap M$ to $\cD(b)$. Thus, the first condition makes $\overline{\cD(a \rev{a}) \cap M} \supset \cD(b \rev{b}) \cap M$ and, in particular, $\overline{\cD(a \rev{a})} \cap \cD(b \rev{b})   \supset  \cD(b \rev{b}) \cap M \neq \emptyset$.

Since $\cD(a \rev{a})$ is irreducible, the second condition then implies that $\overline{\cD(a \rev{a})} \not\supseteq \cD(b \rev{b})$.
\end{proof}

We note for future reference that we can parametrize $M$ as
\[ \mu_{x \rev{x}} \left( \begin{sbm} t_1&1 \\ 1&0 \end{sbm},  \begin{sbm} t_2&1 \\ 1&0 \end{sbm}, \cdots,  \begin{sbm} t_r&1 \\ 1&0 \end{sbm},  \begin{sbm} t_r&1 \\ 1&0 \end{sbm}^{-1},  \cdots,  \begin{sbm} t_2&1 \\ 1&0 \end{sbm}^{-1},   \begin{sbm} t_1&1 \\ 1&0 \end{sbm}^{-1} \right) . \]

Unfortunately, nonempty words of the form $x \rev{x}$ are never reduced, so we need to work harder to give an example with a reduced word, which will thus give an example in a Richardson.
To this end, find a word $x$ in the letters $\{2,3,\ldots, n-2\}$, and distinct distinguished subwords $a$ and $b$ of $x$ such that 
\begin{enumerate}
\item $x s_1 s_{n-1} \rev{x}$ is reduced
\item $\overline{\cD(a)} \supset \cD(b)$ in $\BSo(x)$
\item $\dim \cD(a \bullet \bullet \rev{a}) \leq \dim \cD(b  \bullet \bullet  \rev{b})$ in $\BSo(x  \bullet \bullet  \rev{x})$.
\end{enumerate}
An example is $x = (2,3,2)$, $a = (2,\bullet,2)$ and $b = (2,3, \bullet)$.
The reader can check that $s_2 s_3 s_2 s_1 s_4 s_2 s_3 s_2$ is reduced. 
To see that  $\overline{\cD(2, \bullet, 2)} \supset \cD(2,3, \bullet)$ in $\BSo(2,3,2)$, note that, although $\BSo(x)$ is a space of sequences of $5$-dimensional flags, the first and fourth subspaces in these flags never change, so it is morally a sequence of $3$-dimensional flags and, thought of in this way, we already did the computation above.
Finally, we compute that $\cD(2, \bullet, 2, \bullet, \bullet, 2, \bullet, 2) \cong \cD(2,3,\bullet, \bullet, \bullet, \bullet, 3,2) \cong \GG_m^4 \times \GG_a^2$.

\begin{prop}
For $x$, $a$ and $b$ as above, we have $\overline{\cD(a \bullet \bullet \rev{a})} \cap \cD(b  \bullet \bullet \rev{b}) \neq \emptyset$, but $\overline{\cD(a \bullet \bullet \rev{a})} \not\supseteq \cD(b \bullet \bullet \rev{b})$.
Thus, the Deodhar pieces do not form a stratification of $\BSo(x s_1 s_n \rev{x})$.
\end{prop}

\begin{proof}
We put $x = (i_1, i_2, \ldots, i_r)$ and we let $w$ be the product $s_{i_1} s_{i_2} \cdots s_{i_r}$.

Once again, we introduce a subvariety $M$ of $\BSo(x s_1 s_n \rev{x})$:
Let $M$ be the subvariety of $\BSo(x s_1 s_n \rev{x})$ which can be parametrized as
\[  \mu_{x s_1 s_n \rev{x}}\left( \begin{sbm} t_1&1\\1&0 \end{sbm}, \cdots,  \begin{sbm} t_r&1\\1&0 \end{sbm}, \begin{sbm} 1&0 \\ u&1 \end{sbm}, \begin{sbm} 1&0 \\ v&1 \end{sbm}, \begin{sbm} t_r&1 \\ 1&0 \end{sbm}^{-1}, \dots, \begin{sbm} t_1&1 \\ 1&0 \end{sbm}^{-1} \right) .\]
Projection onto the first $r$ flags makes $M$ into a $\GG_m^2$ bundle over $\BSo(x)$, and this bundle is trivial by Lemma~\ref{lem:TrivialBundle}, so $M \cong \BSo(x) \times \GG_m^2$.
We claim that this identifies $M \cap \cD(a \bullet \bullet \rev{a})$ with $\cD(a) \times \GG_m^2$ (and respectively for $b$). Once we check this, we will have $\overline{\cD(a \bullet \bullet \rev{a})} \cap \cD(b \bullet \bullet \rev{b}) \supset M \cap \cD(b \bullet \bullet \rev{b})$ as before, and the dimensionality argument is exactly as before. 

However, the claim that the projection identifies $M \cap \cD(a \bullet \bullet \rev{a})$ with $\cD(a) \times \GG_m^2$  now requires computation.

Let $(v^0, v^1, \ldots, v^{r-1}, v^r, v^r, v^r,  v^{r-1}, \ldots, v^1, v^0)$ be the distinguished sequence corresponding to $(a, \bullet, \bullet, a^R)$. 
Let $(F^0, F^1, \ldots, F^{2r+2})$ be a list of flags in $M$ which corresponds to $\cD(a) \times \GG_m^2$.
By definition, $F^j$ is in the Schubert cell $\Xo_{v^j}$ for $0 \leq j \leq r$; we must check that this also holds for $r+1 \leq j \leq 2r+2$.

We first do the case of $F^{r+1}$.
 For brevity, abbreviate $z_{i_k}(t_k)$ to $q_k$  and put $g_k = q_1 q_2 \cdots q_k$.
 Then $F^r = g_r B_+$ and $F^{r+1} = g_r y_1(u) B_+$. 
 Note that $z_i(t) = s_i y_i(t)$, so $g_r = s_{i_1} y_{i_1}(t_1) s_{i_2} y_{i_2}(t_2)  \cdots s_{i_r} y_{i_r}(t_r)$.
  Commuting the $s$'s to the left, and using that $x$ is a reduced word for $w$, we can write $g_r = s_{i_1} s_{i_2} \cdots s_{i_r} Y=wY$ for some $Y$ in $N_-$. 
  Then $F^{r+1} =  g_r y_1(u) = (g_r y_1(u) g_r^{-1}) g_r B_+ = (w Y y_1(u) Y^{-1} w^{-1}) g_r B_+ =  (w Y y_1(u) Y^{-1} w^{-1}) F^r$. The conjugate $Yy_1(u)Y^{-1}$ is a matrix in $N_-$ which is $0$ in position $(i,j)$ for $2 \leq j < i \leq n$.
  Since the word $x$ contains no $1$'s, the conjugate $w Yy_1(u)Y^{-1} w^{-1}$ is then also in $N_-$. So $(w Y y_1(u) Y^{-1} w^{-1}) B_+$ must be in the same Schubert cell as $F^r$, as required by the first of the two central bullet points in the word $a \bullet \bullet \rev{a}$.
  
  We now do the case of $F^{2r+2-k}$ for $0 \leq k \leq r$. We have $F^{2r+2-k} = g_r y_1(u) y_n(v) q_r^{-1} \cdots q_{k+1}^{-1} B_+ = (g_r y_1(u) y_n(v) g_r^{-1}) g_r q_r^{-1} \cdots q_{k+1}^{-1} B_+ =  (g_r y_1(u) y_n(v) g_r^{-1}) g_k B_+ =  (g_r y_1(u) y_n(v) g_r^{-1}) F^k$.
  As before, we have $g_r y_1(u) y_n(v) g_r^{-1} = w Y y_1(u) y_n(v) Y^{-1} w^{-1}$ for some $Y \in N_-$. The conjugate $Y y_1(u) y_n(v) Y^{-1}$ is an element of $N_-$ which is $0$ in position $(i,j)$ for $2 \leq j < i \leq n-1$. Since $x$ has not $1$'s or $n$'s, the conjugate $w Yy_1(u) y_n(v) Y^{-1} w^{-1}$ is again in $N_-$. So $F^{2r+2-k}$ and $F^k$ are in the same Schubert cell, as desired.
\end{proof}

\begin{remark}
The reader might wonder if we could reduce clutter by looking at a word $x$ in $\{1,2,\ldots, n-2 \}$ such that $x s_{n-1} \rev{x}$ is reduced instead. Unfortunately, the only such $x$'s are of the form $s_k s_{k+1} \cdots s_{n-2}$, and these do not have distinguished subwords $a$ and $b$ with the requisite properties.
\end{remark}

\subsection{Unipeak and univalley words} \label{sec:unipeak}
For each $w$, there are two particularly nice reduced words for $w$, for which the Deodhar pieces are unusually simple.
Let $(s_{i_1}, s_{i_2}, \cdots, s_{i_a})$ be a reduced word for $w$, and consider the wiring diagram for this word.
See Section~\ref{sec:BSVarieties} for notations associated with wiring diagrams.

We will say that $(s_{i_1}, s_{i_2}, \cdots, s_{i_a})$ is \newword{unipeak} if, for each wire $\sigma_j$, the $x$-coordinates of the crossings where $\sigma_j$ crosses upward are all less than the $x$-coordinates of the crossings where $\sigma_j$ crosses downward.
We note that a unipeak word is automatically reduced. To see this, suppose that $(s_{i_1}, s_{i_2}, \cdots, s_{i_a})$  is not reduced, so there are some $\sigma_{h_1}$ and $\sigma_{h_2}$ which cross at both $x_{j_1}$ and $x_{j_2}$, and not crossing between these points. Without loss of generality, suppose that $x_{j_1} < x_{j_2}$, with $\sigma_{h_1}$ crossing downward at $x_{j_1}$. Then $\sigma_{h_1}$ crosses upward at $x_{j_2}$, making our word non-unipeak.

We will similarly say that  $(s_{i_1}, s_{i_2}, \cdots, s_{i_a})$ is \newword{univalley} if, for each wire $\sigma_j$, the values of $x$ at which $\sigma_j$ crosses downward are all less than the values of $x$ at which $\sigma_j$ crosses upward.

\begin{eg}
Consider $w = 321$ in $S_3$. The unique unipeak word is $s_2 s_1 s_2$ (shown on the left); the unique univalley word is $s_1 s_2 s_1$ (shown on the right).\\
\centerline{
\begin{tikzpicture}
\draw (-0.5,1) -- (1.75,1) -- (2.25,2) -- (2.75,2) -- (3.25,3) -- (4.5, 3);
\draw (-0.5,2) -- (0.75,2) -- (1.25,3) -- (2.75,3) -- (3.25,2) -- (4.5, 2);
\draw (-0.5,3) -- (0.75,3) -- (1.25,2) -- (1.75,2) -- (2.25,1) -- (4.5, 1);
\end{tikzpicture} \qquad
\begin{tikzpicture}
\draw (-0.5,3) -- (1.75,3) -- (2.25,2) -- (2.75,2) -- (3.25,1) -- (4.5, 1);
\draw (-0.5,2) -- (0.75,2) -- (1.25,1) -- (2.75,1) -- (3.25,2) -- (4.5, 2);
\draw (-0.5,1) -- (0.75,1) -- (1.25,2) -- (1.75,2) -- (2.25,3) -- (4.5, 3);
\end{tikzpicture}
}
\end{eg}

There is an elegant geometric way to construct a unipeak word for $w$. 
Draw two axes -- one, which we will call the $p$-axis, pointing northeast and one, which we will call the $q$-axis, pointing northwest.
Plot points at coordinates $(p,q) = (p, w(p))$ and draw line segments from this point to $(0, w(p))$ and $(p, 0)$ to form a ``peak", then extend this peak with horizontal rays to the left and right.
This gives a unipeak wiring diagram for $w$.
We omit the proof, but all unipeak wiring diagrams for $w$ are the same as this one, up to planar isotopy; in other words, there is a unique unipeak reduced word for $w$ up to commutation.

\begin{eg}
Let $w = 3241 = s_2 s_1 s_2 s_3$. 
In the figure below, we use this construction to find the unipeak wiring diagram for $w$ (shown in thick lines).
We observe that we obtain the same reduced word, and the same wiring diagram up to planar isotopy, as in Example~\ref{eg:wiring}.

\centerline{\begin{tikzpicture}
\draw[->] (0,0)--(5,5) node[right]{$p$};
\draw[->] (0,0)--(-5,5) node[left]{$q=w(p)$};
\draw (0.9,1.1) -- (1.1, 0.9) node[below]{$1$} ;
\draw (1.9,2.1) -- (2.1, 1.9) node[below]{$2$} ;
\draw (2.9,3.1) -- (3.1, 2.9) node[below]{$3$} ;
\draw (3.9,4.1) -- (4.1, 3.9) node[below]{$4$} ;
\draw (-0.9,1.1) -- (-1.1, 0.9) node[below]{$1$} ;
\draw (-1.9,2.1) -- (-2.1, 1.9) node[below]{$2$} ;
\draw (-2.9,3.1) -- (-3.1, 2.9) node[below]{$3$} ;
\draw (-3.9,4.1) -- (-4.1, 3.9) node[below]{$4$} ;
\draw[ultra thick] (5,1) -- (1,1) -- (-2,4) -- (-3,3) -- (-5,3) ;
\draw[ultra thick] (5,2) -- (2,2) -- (0,4) -- (-2,2) -- (-5,2) ;
\draw[ultra thick] (5,3) -- (3,3) -- (-1,7) -- (-4,4) -- (-5,4) ;
\draw[ultra thick] (5,4) -- (4,4) -- (3,5) -- (-1,1) -- (-5,1) ;
\end{tikzpicture}}
\end{eg}

%

Let $C$ be a chamber of the unipeak diagram and consider the geometric construction of the unipeak word above.
The top of $C$ consists of two line segments, one with slope $1$ and one with slope $-1$.
Let the line segment with slope $1$ come from the wire $\sigma_i$ and let the line segment with slope $-1$ come from wire $\sigma_j$; we have $i \leq j$.
We will say that $(i,j)$ is the \newword{roof} of $C$.

Let $(s_{i_1}, s_{i_2}, \ldots, s_{i_a})$ be the unipeak word for $w$ and let $(F_1, F_2, \ldots, F_{n-1})$ be a flag in $\Xo^w$.
Since $(s_{i_1}, s_{i_2}, \ldots, s_{i_a})$ is reduced, there is a unique chain of flags $(F^0, F^1, \ldots, F^a)$ in $\BSo(i_1, i_2, \ldots, i_a)$ ending with $F^a = (F_1, F_2, \ldots, F_{n-1})$, and thus a unique labeling of the chambers by subspaces. 

\begin{prop}
In the above notation, the subspace in chamber $C$ is $\Span(e_1, e_2, \ldots, e_{i-1}) + F_{w^{-1}(j)-1}$.
\end{prop}

\begin{proof}[Proof sketch]
It is clear that $C$ is above the height $i-1$ chamber on the far left of the diagram and above the height $w^{-1}(j)-1$ chamber on the far right.
These open chambers are labeled with $\Span(e_1, e_2, \ldots, e_{i-1})$ and $F_{w(j)-1}$ respectively, so $V$ contains both of these spaces and therefore $V \supseteq \Span(e_1, e_2, \ldots, e_{i-1}) + F_{w(j)-1}$.
We will now show that $\dim V = \dim {\big(} \Span(e_1, e_2, \ldots, e_{i-1}) + F_{w(j)-1} {\big)}$, showing equality.
For notational convenience, put $a = i-1$ and $b = w(j)-1$.

The dimension of $V$ is the height of chamber $C$, which is the number of wires passing below any point of  $C$.
The wire $\sigma_k$ passes below $C$ if and only if $k<i$ or $w(k) < w(j)$ (or both). 
So the height of $C$ is $\# \left( \{ k : k \leq a \} \cup \{ k : w(k) \leq b \} \right) = \# \left( [a] \cup w^{-1}([b]) \right)$.
We want to show that $ \# \left( [a] \cup w^{-1}([b]) \right)$ is equal to $\dim \Span(e_1, e_2, \ldots, e_{a}) + F_b$.

We have $\dim {\big(} \Span(e_1, e_2, \ldots, e_{a}) + F_b {\big)} = a+b - \dim{\big(} \Span(e_1, e_2, \ldots, e_{a}) \cap F_b {\big)}$. 
Because the flag $F$ is in $\Xo^w$, we have $\dim{\big(} \Span(e_1, e_2, \ldots, e_{a}) \cap F_b {\big)} = \#(w[a] \cap [b])$. 
So $\dim {\big(} \Span(e_1, e_2, \ldots, e_{a}) + F_b {\big)} = a+b - \#(w[a] \cap [b]) = \#(w[a] \cup [b])$.
Applying the bijection $w$, we have  $\#(w[a] \cup [b]) =  \# \left( [a] \cup w^{-1}([b]) \right)$ as desired.
\end{proof}

\begin{eg}
We take the chambers of Example~\ref{eg:WiringChambers} and fill them as described here:

\centerline{
\begin{tikzpicture}[xscale=1.5, yscale=1.5]
\draw (-1,1) -- (1.75,1) -- (2.25,2) -- (2.75,2) -- (3.25,3) -- (3.75,3) -- (4.25,4) -- (6,4) ;
\draw (-1,2) -- (0.75,2) -- (1.25,3) -- (2.75,3) -- (3.25,2) -- (6,2) ;
\draw (-1,3) -- (0.75,3) -- (1.25,2) -- (1.75,2) -- (2.25,1)  -- (6,1) ;
\draw (-1,4) -- (3.75,4) -- (4.25,3)  -- (6,3) ;
\node (L0) at (0.5,1.5) {$\Span(e_1)$};
\node (L1) at (4,1.5) {$F_1$};
\node (P0) at (0,2.5) {$\Span(e_1, e_2)$};
\node (P1) at (2,2.5) {$\Span(e_1) + F_1$};
\node (P2) at (4.5,2.5) {$F_2$};
\node (S0) at (1.5,3.5) {$\Span(e_1, e_2, e_3)$};
\node (S1) at (5,3.5) {$F_3$};
\end{tikzpicture}
}
\end{eg}

We give the analogous formula for univalley wiring diagrams. 
If $C$ is a chamber of a univalley wiring diagram then there are two wires running along the bottom of $C$; let $\sigma_i$ be the decreasing wire and let $\sigma_j$ be the increasing wire.
Then the subspace in $C$ is $\Span(e_1, e_2, \ldots, e_i) \cap F_{w(j)}$. 

Using this, we can give an explicit description of the Deodhar strata for a unipeak wiring diagram.

\begin{prop}
Let $w$ be a permutation, let $(s_{i_1}, s_{i_2}, \ldots, s_{i_n})$ be the unipeak wiring diagram for $w$ and let $F$ be a flag in $\Xo^w$. Then the knowledge of which Deodhar piece $F$ is in is equivalent to the knowledge, for all $i \leq i'$ and all $j$, of $\dim {\big(} \Span(e_1, e_2, \ldots, e_{i-1}, e_{i'+1}, e_{i'+1}, \ldots, e_n) + F_j {\big)}$. 
If we let $F = gB_+$ then, equivalently, the knowledge of which Deodhar piece $F$ is in is equivalent to the knowledge, for all $i \leq i'$ and all $j$, of the rank of the submatrix of $g$ in rows $\{ i, i+1, \ldots, i' \}$ and columns $\{ 1,2,\ldots, j \}$.
\end{prop}

\begin{rem}
The author learned this result from Allen Knutson.
The author is not aware of a published source for this statement.
\end{rem}

\begin{remark}
We can think of $\Span(e_1, e_2, \ldots, e_{i-1}, e_{i'+1}, e_{i'+1}, \ldots, e_n)$ as $E_{i-1} + E^{\text{op}}_{n-i'}$ where $E_{\bullet}$ is the standard flag, and $E^{\text{op}}_{\bullet}$ is the opposite flag. More generally,  Curtis~\cite{Curtis}  and Shapiro, Shapiro and Vainshtein~\cite{SSV} consider an arbitrary pair of flags $E_{\bullet}$ and $E'_{\bullet}$ and decompose $\Fl_n$ according to the values of $\dim(E_i + E'_j + F_k)$ for $F_{\bullet} \in \Fl_n$ and all $(i,j,k)$; each stratum is of the form $\GG_a^r \times \GG_m^s$. It would be interesting to incorporate this more general decomposition into the discussion here.
\end{remark}

\begin{proof}[Proof sketch]
Knowing what Deodhar stratum $F$ is in is equivalent to knowing what Grassmannian Schubert cell the various subspaces in the various chambers are in.
We have now seen that the chamber labels are equal to $\Span(e_1, e_2, \ldots, e_i) \cap F_j$ for various $j$. 
Knowing which Grassmannian Schubert cell $V$ is in is equivalent to knowing $\dim {\big(} V + \Span(e_{i'+1}, e_{i'+2}, \ldots, e_n)$ for all $i'$.
Thus, knowing which Deodhar stratum $F$ is in is equivalent to knowing $\dim {\big(} \Span(e_1, e_2, \ldots, e_{i-1}, e_{i'+1}, e_{i'+1}, \ldots, e_n) + F_j {\big)}$ for all $1 \leq i \leq i' \leq n$ and all $1 \leq j \leq n$.

Since $F_j$ has as a basis the first $j$ columns of $g$, the dimension of $\Span(e_1, e_2, \ldots, e_{i-1}, e_{i'+1}, e_{i'+1}, \ldots, e_n) + F_j$ is $(i-1) + (n-i') + \text{rank}( g_{[i,i'] \times [j]})$, where $g_{[i,i'] \times [j]}$ is the submatrix of $g$ in rows $\{ i, i+1, \ldots, i' \}$ and columns $\{ 1,2,\ldots, j \}$.
\end{proof}

Analogously, for a univalley word, knowing which Deodhar stratum $F$ is in is equivalent to knowing $\dim {\big(} \Span(e_i, e_{i+1},\ldots, e_{i'}) \cap F_j {\big)}$  for all $1 \leq i \leq i' \leq n$ and all $1 \leq j \leq n$.

\section{Total positivity} \label{sec:Positivity}
A real matrix $A$ is called \newword{totally positive} if every minor $\Delta_{I,J}(A)$ is positive, and is called \newword{totally nonnegative} if every minor $\Delta_{I,J}(A)$ is nonnegative.
Totally positive matrices are common in algebraic statistics and representation theory and have been studied for over a century; we recommend~\cite{FZPositivity} and the sources therein for an overview.

The Cauchy-Binet identity states that 
\[ \Delta_{I,K}(AB) = \sum_J \Delta_{I,J}(A) \Delta_{J,K}(B) \]
where the matrix $A$ is $i \times j$, the matrix $B$ is $j \times k$, the sets $I$ and $K$ are $r$-element subsets of $[i]$ and $[k]$ respectively, and the sum is over all $r$-element subsets $J$ of $[j]$.
Thus, the Cauchy-Binet identity implies that the totally nonnegative matrices form a semigroup under multiplication.
A Loewner-Whitney theorem (proved in~\cite{Loewner55}, making key use of results from \cite{Whitney52}) states that the semigroup of totally nonnegative matrices in $\SL_n(\RR)$ is generated by the Chevalley generators $\Id + t e_{i, i+1}$ and $\Id + t e_{i+1, i}$ for $t>0$. (Here $e_{ij}$ is the matrix which is $1$ in position $(i,j)$ and $0$ everywhere else.)
There is an enormous literature on the structure of the semigroup of totally nonnegative matrices, but we need to move on to totally nonnegative points of the flag manifold.

\subsection{Totally nonnegative subspaces and flags}
A subspace $V \subset \RR^n$ is called \newword{totally nonnegative} if all of its Pl\"ucker coordinates are nonnegative. 
More precisely, since the Pl\"ucker coordinates are homogeneous coordinates, we impose that they are all nonnegative up to a global sign flip.
We denote the set of totally nonnegative subspaces of dimension $k$ in $\RR^n$ by $G(k,n)_{\geq 0}$ and call it the \newword{totally nonnegative Grassmannian}.
A flag $F_1 \subset F_2 \subset \cdots \subset F_{n-1}$ is called \newword{totally nonnegative} if all of the subspaces $F_k$ are totally nonnegative.
We will denote the set of totally nonnegative points of the flag manifold as $\Fl_n^{\geq 0}$.

\begin{eg} \label{FL3Positive}
An element of $\Fl_3(\RR)$ is a pair of a point $p=(\Delta_1 \! : \! \Delta_2 \! : \! \Delta_3)$ in $\PP^2(\RR)$ and a line $\ell = \{ (x_1\! : \!x_2\! : \!x_3) \! : \! \Delta_{23} x_1 - \Delta_{13} x_2 + \Delta_{12} x_3 = 0 \}$ in $\PP^2(\RR)$ with $p \in \ell$.
The condition that $\Delta_1$, $\Delta_2$ and $\Delta_3 \geq 0$ says that $p$ lies in a projective triangle with vertices $(1\! : \!0\! : \!0)$, $(0\! : \!1\! : \!0)$ and $(0\! : \!0\! : \!1)$.
The condition that $\Delta_{12}$, $\Delta_{13}$ and $\Delta_{23} \geq 0$ means that the line $\ell$ crosses the boundary of this triangle on the edge from $(1\! : \!0\! : \!0)$ to $(0\! : \!1\! : \!0)$, and again on the edge from $(0\! : \!1\! : \!0)$ to $(0\! : \!0\! : \!1)$.
See the figure below.

\centerline{
\begin{tikzpicture}
\node  [label={[xshift=-1.0cm, yshift=-0.6cm]${[1:0:0]}$}] (v1) at (0,0) {};
\node  [label={[xshift=0.0cm, yshift=0.0cm]${[0:1:0]}$}] (v2) at (1,1.73) {};
\node  [label={[xshift=1.0cm, yshift=-0.6cm]${[0:0:1]}$}] (v3) at (2,0) {};
\draw (v1.center) -- (v2.center) -- (v3.center) -- (v1.center) ;
\draw[ultra thick] (-1,0.5) -- (3, 1.5) node [right,above] {$\ell$};;
\node  [label={[xshift=0.3cm, yshift=-0.6cm]$p$}] (p) at (1, 1) {$\bullet$};
\end{tikzpicture}
}
\end{eg}

\begin{lemma}
If $V$ is a totally nonnegative subspace and $g \in \GL_n$ is a totally nonnegative matrix, then the subspace $gV$ is totally nonnegative.
If $F_{\bullet}$ is a totally nonnegative flag and $g \in \GL_n$ is a totally nonnegative matrix, then the flag $gF_{\bullet}$ is totally nonnegative.
\end{lemma}

\begin{proof}
The first sentence follows from the Cauchy-Binet identity.
The first sentence implies the second, by considering the application of $g$ to each subspace $F_k$.
\end{proof}

It turns out that every totally nonnegative subspace is included in a totally nonnegative flag. 
The following result appears as~\cite[Theorem 3.8]{Post}, but the proof there is only a sketch, so we provide some details.

\begin{prop}\label{prop:SubspaceInFlag}
 If $V$ is a totally nonnegative subspace with dimension $k$, then there is a totally nonnegative flag with $F_k = V$.
 \end{prop}

 \begin{proof}
 We need to construct a sequence of totally nonnegative subspaces $F_1 \subset F_2 \subset \cdots \subset F_{k-1} \subset V \subset F_{k+1} \subset \cdots \subset F_{n-1}$.
 We will show how to construct $F_{k-1}$ from $V$; iterating in this manner, we get a chain $F_1 \subset F_2 \subset \cdots \subset F_{k-1} \subset V$.
 A similar construction makes the chain $V \subset F_{k+1} \subset \cdots \subset F_{n-1}$.
 
We write $x_1$, $x_2$, \dots, $x_n$ for the coordinates on $\RR^n$.
Let $j$ be the index such that $x_{j+1} = x_{j+2} = \cdots = x_n=0$ on $V$ and $x_j$ is not identically $0$ on $V$ and let $V' = V \cap \{ x_j = 0 \}$. 
Then we can choose a basis $\vec{v}_1$, $\vec{v}_2$, \dots, $\vec{v}_k$ for $V$ where $\vec{v}_1$, $\vec{v}_2$, \dots, $\vec{v}_{k-1}$ is a basis for $V'$ and $\vec{v}_k$ is of the form $e_j + \sum_{i<j} c_i e_i$. 
Computing Pl\"ucker coordinates using the matrix whose columns are the $\vec{v}_i$, we have 
\[ \Delta_{J}(V') = \begin{cases}\Delta_{J \cup \{ j \}}(V) & J \subseteq [j-1] \\ 0 &  J \not\subseteq [j-1] \end{cases} . \]
So $V'$ is totally nonnegative, and we take $F_{k-1} = V'$.
\end{proof}

Lusztig~\cite{Lusztig94} originally defined the totally nonnegative flag variety differently, in a way which we now know to be equivalent.
Let $B_-^{\geq 0}$ be the semigroup of totally nonnegative matrices in the lower Borel $B_-$. 
Then $\Fl_n^{\geq 0}$ was defined to be the closure, in $\Fl_n^{\geq 0}(\RR)$, of the set of flags $gB_+$ for $g \in B_-^{\geq 0}$.
Lusztig also showed that the nonnegative flag variety can be defined by inequalitites using canonical basis elements; see~\cite[Proposition 8.17]{Lusztig94} and~\cite[Theorem 3.4]{LusztigPartial}.
However, this is not the same as defining the nonnegative flag variety using Pl\"ucker coordinates, as we have done.
Fortunately, there are now short proofs that these definitions are equivalent; see~\cite[Theorem 0.8]{Lusztig23arXiv} or~\cite[Theorem 1.1]{BK1}.
So we will cite old results which use either definition without concerning our self about which definition they are using.

We will now describe several results of Rietsch~\cite{RietschCellular} and Marsh and Rietsch~\cite{MR}, many of which were first conjectured in~\cite{Lusztig94}, saying that the totally nonnegative part of the flag variety is extremely well behaved, and interacts very nicely with the Richardson and Deodhar decompositions.
We write $\Ro_u^{w, \geq 0} := \Ro_u^w(\RR) \cap \Fl_n^{\geq 0}$.


\begin{theorem} \label{thm:DeodharPositivity}
Let $u \preceq w$. Then $\Ro_u^{w, \geq 0} \cong \RR_{>0}^{\ell(w) - \ell(u)}$. More precisely, let  $(i_1, i_2, \ldots, i_a)$ be any reduced word for $w$. 
Let $\cD$ be the corresponding Deodhar torus in $\Ro_u^w$, so $\cD(\RR) \cong (\RR_{>0} \sqcup \RR_{<0})^{\ell(w) - \ell(u)}$. 
Then $\Ro_u^{w, \geq 0}$ is one of the connected components of $\cD(\RR)$.
To be precise, if we use the parametrization of $\cD(\RR)$ where we  use $\dot{s}_{i_j}$ and $y_{i_j}(t_j)$, then $\Ro_u^{w, \geq 0}$ is the connected component where the $t_j$ are positive for $j \in J_{=}$. 
\end{theorem}

\begin{proof}
The full statement is~\cite[Theorem 1.3]{MR}. Parts of this result appear in~\cite[Theorem 2.8]{RietschCellular}.
\end{proof}

\begin{remark}
Since $\RR$ and $\RR_{>0}$ are isomorphic as smooth manifolds, $\Ro_u^w(\RR)$ is isomorphic to both $\RR^{\ell(w) - \ell(u)}$ and $\RR_{>0}^{\ell(w) - \ell(u)}$ as a smooth manifold.  
We will always write $\RR_{>0}^{\ell(w) - \ell(u)}$, since the natural coordinates on $\Ro_u^w$ are multiplicative coordinates on a Deodhar torus.
\end{remark}

\begin{eg}
We consider the largest Richardson $\Ro_{123}^{321}$ in $\Fl_3$. 
The totally positive points are of the form $(\Delta_1: \Delta_2 : \Delta_3) \times (\Delta_{12} : \Delta_{13} : \Delta_{23})$ with all the Pl\"ucker coordinates positive. 
There are two reduced words for $321$, namely, $s_1 s_2 s_1$ and $s_2 s_1 s_2$, corresponding to the two parametrizations
\[ 
\begin{array}{c@{}c@{}l@{}c@{}lc}
(t_1, t_2, t_3) &\mapsto&
\begin{bmatrix} 1&0&0 \\ t_1&1&0 \\ 0&0&1 \\ \end{bmatrix}
\begin{bmatrix} 1&0&0 \\ 0&1&0 \\ 0&t_2&1 \\ \end{bmatrix}
\begin{bmatrix} 1&0&0 \\ t_3&1&0 \\ 0&0&1 \\ \end{bmatrix} &=& 
\begin{bmatrix} 1&0&0 \\ t_1+t_3 & 1 & 0 \\ t_2 t_3 & t_2 & 1 \\ \end{bmatrix} 
& \text{and} \\[0.25 in]
(u_1, u_2, u_3) &\mapsto&
\begin{bmatrix} 1&0&0 \\ 0&1&0 \\ 0&u_1&1 \\ \end{bmatrix}
\begin{bmatrix} 1&0&0 \\ u_2&1&0 \\ 0&0&1 \\ \end{bmatrix}
\begin{bmatrix} 1&0&0 \\ 0&1&0 \\ 0&u_3&1 \\ \end{bmatrix} &=& 
\begin{bmatrix} 1&0&0 \\ u_2 & 1 & 0 \\ u_1 u_2 & u_1 + u_3 & 1 \\ \end{bmatrix}. & \\
\end{array}
\]
If $(t_1, t_2, t_3)$ and $(u_1, u_2, u_3)$ each range over $\RR_{\neq 0}^3$, then these maps have different images; the image of the first map is $\Delta_1 \Delta_3 \Delta_{12} \Delta_{23} \Delta_{13} \neq 0$ and the image of the second map is $\Delta_1 \Delta_3 \Delta_{12} \Delta_{23} \Delta_2 \neq 0$. 
(The first four factors are the same and the fifth factor is different.)
However, if $(t_1, t_2, t_3)$ and $(u_1, u_2, u_3)$ each range over $\RR_{> 0}^3$, then both maps have the same image; the points where all six of the Pl\"ucker coordinates are positive.
This example illustrates how $\Ro_u^{w, \geq 0}$ is independent of the choice of reduced word for $w$, even through the Deodhar pieces depend on the choice of reduced word.

In terms of the geometric depiction in Example~\ref{FL3Positive}, the condition $\Delta_2=0$ corresponds to $p$ lying on the bottom edge of the triangle. As the reader can see, one can have a flag such $p$ lies on this bottom edge but all Pl\"ucker coordinates of $\ell$ are nonzero; such a flag would be in the image of the first parametrization and not the second. However, in such a flag, $\ell$ will either meet the boundary of the triangle in the left and bottom edges, or else in the right and bottom edges, so such an $\ell$ will not be positive.
\end{eg}

\begin{remark} \label{rem:Unipotent}
We describe special features of the case $u=e$. 
The Schubert cell $\Xo_e$ is the largest Schubert cell, isomorphic to $\AA^{\binom{n}{2}}$, and can be concretely be identified with the lower unipotent group $U_-$; the map $u \mapsto u B_+$ is an isomorphism $U_- \longrightarrow \Xo_e$. 
The totally nonnegative flags in $\Xo_e^{\geq 0}$ thus decompose into a union of these Richardson pieces:
\[ \Xo_e^{\geq 0} = \bigsqcup_{w \in S_n} \Ro_e^{w, \geq 0} . \]
Let $(s_{i_1}, s_{i_2}, \ldots, s_{i_a})$ be a reduced word for $w$.
Since $u=e$, we have $J_{\uparrow} = \emptyset$, so we are only multiplying together $y_i(t)$'s, not $\dot{s}_i$'s.
Theorem~\ref{thm:DeodharPositivity} tells us that  $(t_1, t_2, \ldots, t_a) \mapsto y_{i_1}(t_1) y_{i_2}(t_2) \cdots y_{i_a}(t_a) B_+$ is an isomorphism $\RR_{>0}^a \longrightarrow \Ro_e^{w, \geq 0}$. The product $y_{i_1}(t_1) y_{i_2}(t_2) \cdots y_{i_a}(t_a)$ is in $U^{\geq 0}_-$ since each factor is in $U^{\geq 0}_-$ and the totally nonnegative matrices form a semigroup and $U_-$ is a group.
It turns out that the image of the map   $(t_1, t_2, \ldots, t_a) \mapsto y_{i_1}(t_1) y_{i_2}(t_2) \cdots y_{i_a}(t_a)$ from $\RR_{>0}^a$ to $U_-$ is independent of the choice of reduced word; the image of this map is called a \newword{unipotent cell}, and we will denote it by $U_-^{w, \geq 0}$.
So $U_-^{\geq 0}$ decomposes into the disjoint union $\bigcup_{w \in S_n} U_-^{w, \geq 0}$, with $U_-^{w, \geq 0} \cong \RR_{>0}^{\ell(w)}$.
\end{remark}

\begin{remark}
When $u \neq 0$, we could likewise consider the set of matrices in $\GL_n$ that we get by multiplying together the matrices $y_{i_j}(t_j)$ and $\dot{s}_{i_j}$. 
In this case, the image in $\GL_n$ does depend on the choice of reduced word.
As an example, in $n=3$, we take $u=213$ and $w =321$. There are two reduced words for $w$:  $(s_1, s_2, s_1)$ and $(s_2, s_1, s_2)$. The corresponding matrix products are
\[ \begin{array}{lcl}
\begin{bmatrix} 1&0&0 \\ t_1 & 1 & 0 \\ 0&0&1 \\ \end{bmatrix} \begin{bmatrix} 1&0&0 \\ 0& 1 & 0 \\ 0&t_2&1 \\ \end{bmatrix} \begin{bmatrix} 0&-1&0 \\ 1 & 0 & 0 \\ 0&0&1 \\ \end{bmatrix} 
&=& \begin{bmatrix} 0&-1&0 \\ 1 & -t_1 & 0 \\ t_2&0&1 \\ \end{bmatrix} \\[0.3 in]
\begin{bmatrix} 1&0&0 \\ 0 & 1 & 0 \\ 0&u_1&1 \\ \end{bmatrix} \begin{bmatrix} 0&-1&0 \\ 1& 0 & 0 \\ 0&0&1 \\ \end{bmatrix} \begin{bmatrix}1 & 0 & 0 \\ 0&1 &0\\ 0&u_3&1 \\ \end{bmatrix} &=& \begin{bmatrix} 0&-1&0 \\ 1& 0 & 0 \\ u_1&u_3&1 \\ \end{bmatrix} . \\
\end{array}\]
If we take $u_1 = t_2$ and $u_3 = t_1 t_2$, then these are the same flag, but they are clearly not the same matrix.
\end{remark}

We have seen that Deodhar parametrizations are best understood in terms of Bott-Samelson varieties. 
Here is the main result about Bott-Samelson varieties and total positivity:
\begin{theorem} \label{thm:BSPositivity}
Let $F$ be a totally nonnegative flag in $\Ro_u^{w, \geq 0}$. Let $(s_{i_1}, s_{i_2}, \ldots, s_{i_a})$ be a reduced word for $w$, and let $(F^0, F^1, \ldots, F^a)$ be the unique sequence of flags in $\BSo(i_1, i_2, \ldots, i_a)$ with $F^a = F$.
Then all the flags $F^j$ are totally nonnegative.
\end{theorem}

\begin{proof}[Proof sketch]
This result is essentially due to Bethany Marsh and Konstanze Rietsch~\cite{MR}; we explain where to find this result in their work.

The variable $u$ doesn't play a role in the theorem; we simply need that $F$ is in $\Xo^{w, \geq 0} = \bigsqcup_{u \preceq w} \Ro_u^{w, \geq 0}$.
Set $w_1 = s_{i_1} s_{i_2} \cdots s_{i_j}$ and $w_2 = s_{i_{j+1}} s_{i_{j+2}} \cdots s_{i_a}$. The map $\Xo^{w} \to \Xo^{w_1}$ sending $F$ to $F^j$ is the map that Marsh and Rietsch call $\pi^w_{w_1}$, and this result is~\cite[Lemma 11.5]{MR}. The proof of that lemma relies on results of Rietsch in~\cite{RietschCellular}.
\end{proof}

We also remark on the behavior of Pl\"ucker coordinates on $\Ro_u^{w, \geq 0}$. 
We first recall the situation without imposing positivity.
Let $J$ be a $k$-element subset of $[n]$. If $J = u[k]$ or $J = w[k]$, then $\Delta_J$ is nonzero on $\Ro_u^w$.
If there is no $v$ with $u \preceq  v \preceq w$ and $v[k]=J$, then $\Delta_J$ is identically zero on $\Ro_u^w$.
In the last case, where there is a $v$ with $u \prec v \prec w$ and $v[k]  = J$, but where $J \neq u[k]$, $w[k]$, then the function $\Delta_J$ is not identically zero on $\Ro_u^w$, but may vanish at some points of $\Ro_u^w$.

In the situation of positivity, things are more elegant:
\begin{theorem} \label{thm:PositivePluckers}
Let $u \preceq w$ and let $J$ be a $k$-element on $[n]$. If there is a $v$ with $u \preceq  v \preceq w$ and $v[k]=J$ then $\Delta_J > 0$ everywhere on $\Ro_u^{w, \geq 0}$. If there is no such $v$, then $\Delta_J = 0$ everywhere on $\Ro_u^{w, \geq 0}$.
\end{theorem}

Theorem~\ref{thm:PositivePluckers} appears in other language in~\cite[Section 7]{BruhatInterval} and in~\cite[Theorem 1.2]{BK1}.
We provide a direct combinatorial proof.
We first must set up notation.

Fix a reduced word $(s_{i_1}, s_{i_2}, \ldots, s_{i_a})$ for $w$ and a permutation $u \preceq w$.
Put $w^j  = s_{i_1} s_{i_2} \cdots s_{i_j}$. 
Let $v^j$ be the positive subsequence, which can be computed recursively by $v^a = u$ and
\[ v^{j-1} = \begin{cases} v^j s_{i_j}  & v^j  s_{i_j} \prec v^j \\ v^j &   v^j s_{i_j} \succ v^j \end{cases} .\]
It is convenient to define $r_j$ to be $e$ if $j \in J_=$ and to be $s_{i_j}$ if $j \in J_{\uparrow}$, so that $v^j = r_1 r_2 \cdots r_j$. 

For each $j \in J_=$, we take a variable $t_j$. For $j \in [a]$, we define
\[ h_j := \begin{cases} y_{i_j}(t_j) & j \in J_= \\ \dot{s}_{i_j} & j \in J_{\uparrow} \end{cases} . \]
So Theorem~\ref{thm:DeodharPositivity} says that $\Ro_u^{w, \geq 0}$ is parametrized by sending  the vector $(t_j)_{j \in J_=}$ in $\RR_{>0}^{J_=}$ to the flag $h_1 h_2 \cdots h_a B_+$.

\begin{lemma} \label{lem:CoeffsPositive}
For any $1 \leq j \leq a+1$, and any $J \subseteq [n]$, the minor $\Delta_J(h_j h_{j+1} \cdots h_a)$ is a polynomial in the $t$-variables with nonnegative coefficients.
\end{lemma}

\begin{prob}
Give a direct combinatorial rule for computing this polynomial. When $e=1$, so that $J_{\uparrow} = \emptyset$, this polynomial can be computed via flows through wiring diagrams and the Gessel-Viennot formula; see~\cite{FZPositivity}.
\end{prob}

\begin{proof}[Proof of Lemma~\ref{lem:CoeffsPositive}]
The proof is by reverse induction on $j$. When $j=a+1$, we are taking the minors of $\Id$, which are all $1$ or $0$.

So, now assume that we know the result for $j+1$ and want to prove the result for $j$. We define $g' := h_{j+1} h_{j+2} \cdots h_a$ and $i:=i_j$.
We break into two cases: $j \in J_=$ or $j \in J_{\uparrow}$.

\textbf{Case 1:} $j \in J_=$. Then $h_j h_{j+1} \cdots h_a = y_{i}(t_j) g'$. Then the Caucy-Binet formula tells us that 
\[ \Delta_J(y_i(t_j) g') = \begin{cases} t_j \Delta_{s_i J}(g') + \Delta_J(g') & J \cap \{ i, i+1 \} = \{ i+1 \} \\ \Delta_J(g') & \text{otherwise} \\ \end{cases} . \]
So the result follows by induction.

\textbf{Case 2:} $j \in J_{\uparrow}$. Then $h_j h_{j+1} \cdots h_a = \dot{s}_i g$. 
In this case, the Cauchy-Binet formula tells us that
\[ \Delta_J(\dot{s}_i  g) = 
\begin{cases} \phantom{-} \Delta_{s_i J}(g') & J \cap \{ i, i+1 \} = \{ i+1 \} \\ - \Delta_{s_i J}(g') & J \cap \{ i, i+1 \} = \{ i \} \\  \phantom{-} \Delta_J(g') & \text{otherwise} \\ \end{cases} . \]
In the first and third cases, the claim follows by induction, but it appears that the second case will cause a problem.
We claim that, in this case, $\Delta_{s_i J}(g')$ will be identically zero.

Let $k=|J|$.
Put $w' = s_{i_{j+1}} s_{i_{j+2}} \cdots s_{i_a}$ and $u' = r_{j+1} r_{j+2} \cdots r_a$. 
So, inductively,  $g' B_+$ is in $\Ro_{u'}^{w'}$ and all Pl\"ucker coordinates of  $g' B_+$ are polynomials in the $t$'s with nonnegative coefficients.

Suppose, for the sake of contradiction that $\Delta_{s_i J}(g')$ is not identically zero as a polynomial in the $t$'s for some $J$ with $i \in J$ and $i+1 \not\in J$.
Then there is some $v'$ with $u' \preceq v' \preceq w'$ and $v'[k] = s_i J$.
So $i \in s_i v'[k]$ and $i+1 \not\in s_i v'[k]$.
This implies that $s_i v' \prec v'$ (whereas we have $s_i u' \succ u'$ and $s_i w' \succ w'$).

Now, we have $s_i u' \preceq s_i v' \preceq v' \preceq w'$, so there is a reduced subword of $(s_{i_{j+1}}, s_{i_{j+2}}, \cdots, s_{i_n})$ with product $s_i u'$. 
But then $(r_j, r_{j+1}, \cdots, r_n)$ cannot be the rightmost subword with product $s_i u'$, contradicting that $(r_1, r_{2}, \cdots, r_n)$ was supposed to be the positive subword with product $u$.
\end{proof}

Lemma~\ref{lem:CoeffsPositive} immediately implies Theorem~\ref{thm:PositivePluckers}:
\begin{proof}[Proof of Theorem~\ref{thm:PositivePluckers}]
Since $\Delta_J$ is a polynomial in the $t_j$ with nonnegative coefficients, either $\Delta_J$ is identically $0$, or it is positive for all $(t_j) \in \RR_{>0}^{J_=}$. 
\end{proof}

\subsection{Cell complexes and total positivity}
We have decomposed $\Fl_n^{\geq 0}$ into the pieces $\Ro_u^{w, \geq 0}$ where $\Ro_u^{w, \geq 0} \cong \RR_{>0}^{\ell(w) - \ell(u)}$. 
This suggests that $\Fl_n^{\geq 0}$ is a CW complex. 
This is true and was proved by Rietsch and Williams~\cite{RW1}. Galashin, Karp  and Lam~\cite{GKL} showed a much stronger statement, that this CW complex is regular:
\begin{theorem} \label{thm:CellComplex}
$\Fl_n^{\geq 0}$ is a regular CW complex, where the cells are the $\Ro_u^{w, \geq 0}$. This means that the closure $\overline{\Ro_u^{w, \geq 0}}$ of each $\Ro_u^{w, \geq 0}$ is a union of various $\Ro_{u'}^{w', \geq 0}$ and the pair $(\overline{\Ro_u^{w, \geq 0}}, \Ro_u^{w, \geq 0})$ is homeomorphic to the pair $(\text{closed ball},\ \text{interior of ball})$. More precisely, we have $\overline{\Ro_u^{w, \geq 0}} = \bigsqcup_{u \preceq u' \preceq w' \preceq w} \Ro_{u'}^{w', \geq 0}$.
\end{theorem}

\begin{eg}
The cell complex $\Fl_3^{\geq 0}$ has $1$ three-dimensional cell, $4$ two-dimensional cells, $8$ one-dimensional cells and $6$ zero-dimensional cells. 
Topologically, $\Fl_3^{\geq 0}$ is a closed $3$-dimensional ball and its boundary is a $2$-dimensional sphere.
We would draw it, but we have already done so in Example~\ref{Flag3Example}: The $2$-dimensional sphere is formed by gluing the two hexagons along their common hexagonal boundary, and the $3$-dimensional cell $\Ro_{123}^{312, \geq 0}$ fills the interior of the sphere.
\end{eg}


\begin{remark}
Theorem~\ref{thm:CellComplex} was foreshadowed by a similar result for unipotent cells, as discussed in Remark~\ref{rem:Unipotent}.
Recall that this remark decomposed $U_-^{\geq 0}$ into unipotent cells $U_-^{w, \geq 0}$, indexed by $w \in S_n$, with $U_-^{w, \geq 0} \cap \RR_{>0}^{\ell(w)}$. We have  $U_-^{v, \geq 0}$ in the closure of  $U_-^{w, \geq 0}$ if and only if $v \preceq w$.
The topological properties of this decomposition were considered by Fomin and Shapiro~\cite{FominShapiro}.
In order to work in a compact setting, they intersect with the hyperplane $x_{21}+x_{32} + \cdots + x_{n(n-1)}=1$; we'll call this hyperplane $H$.
For $w \neq e$, the intersection $U_-^{w, \geq 0} \cap H$ is isomorphic to $\RR_{>0}^{\ell(w)-1}$. 
(When $w=e$, the intersection is empty.)
Fomin and Shapiro conjectured that that the cells $U_-^{w, \geq 0} \cap H$ form a regular CW decomposition of $U_-^{\geq 0} \cap H$, whose closure relations are given by Bruhat order on $S_n$. 

For example, when $n=3$, the intersection $U_-^{\geq 0} \cap H$ is
\[ \left\{ \begin{sbm} 1&0&0 \\ x&1&0 \\ y & 1-x & 1 \\ \end{sbm} : 0 \leq y \leq x(1-x) \right\} .\]
This region is depicted in the figure below:\\
\centerline{\includegraphics[height=0.50 in]{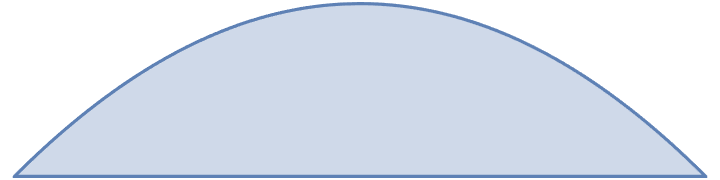}}
Note that each cell $U_-^{w, \geq 0} \cap H$ is an open ball of dimension $\ell(w)-1$, and its closure is a closed ball. For example, the line segment $\{y=0,\ 0 <x< 1 \}$ is $U_-^{231, \geq 0}$ and the parabolic segment $\{ y = x(1-x),\ 0<x<1 \}$ is  $U_-^{312, \geq 0} \cap H$.
The conjecture of Fomin and Shapiro was originally proved by Hersh~\cite{Hersh}, and can now be deduced as a consequence of the work of Galashin, Karp and Lam.
\end{remark}

%
%
%
%
%
%

\subsection{Positivity in partial flag manifolds} \label{sec:PartialPositive}
Lustzig developed a theory of total positivity in $G/P$ for a general reductive $G$ and parabolic $P \supseteq B_+$.
Thus, Lusztig's theory applies to partial flag manifolds. 
Let $0 < k_1 < k_2 < \cdots < k_r < n$. Let $F_{k_1} \subset F_{k_2} \subset \cdots \subset F_{k_r}$ be a partial flag with dimension vectors $(k_1, k_2, \ldots, k_r)$. 
Then $F_{k_1} \subset F_{k_2} \subset \cdots \subset F_{k_r}$  is totally nonnegative in the sense of Lusztig if and only if the partial flag $F_{k_1} \subset F_{k_2} \subset \cdots \subset F_{k_r}$ can be embedded in a totally nonnegative complete flag $F_1 \subset F_2 \subset \cdots \subset F_{n-1}$.

If our partial flag manifold is a Grassmannian, then our flags just have a single subspace $V$.
By Proposition~\ref{prop:SubspaceInFlag}, a subspace $V$ is in the totally nonnegative part of $G(k,n)$, in Lusztig's sense, if and only if the Pl\"ucker coordinates of $V$ are nonnegative.
We note that this is not true for other partial flag manifolds.

\begin{eg} \label{eg:Flag13Positive}
Let $V_1 \subset V_3$ be a partial flag in $\RR^4$. So the Pl\"ucker coordinates of $V_1$ and $V_3$ are $[ \Delta_1 : \Delta_2 : \Delta_3 : \Delta_4]$ and $[ \Delta_{123} : \Delta_{124} : \Delta_{134} : \Delta_{234}]$ with $\Delta_1 \Delta_{234} - \Delta_2 \Delta_{134} + \Delta_3 \Delta_{124} - \Delta_4 \Delta_{123}=0$. 
The condition that the individual subspaces $V_1$ and $V_3$ are totally nonnegative says that the Pl\"ucker coordinates of $V_1$ and $V_3$ are nonnegative; assume this from now on.
We will show that $(V_1, V_3)$ can be completed to a totally nonnegative flag $(V_1, V_2, V_3)$ if and only if, in addition to these conditions, we have $\Delta_1 \Delta_{234} - \Delta_2 \Delta_{134} \geq 0$. 

This can be understood geometrically in $\PP^3(\RR)$. Let $e_1$, $e_2$, $e_3$, $e_4$ be the standard basis of $\RR^4$ and let $x_1$, $x_2$, $x_3$, $x_4$ be the dual coordinates.
Let $T$ be the tetrahedron of points in  $\PP^3(\RR)$ with $x_j \geq 0$, so the $e_i$ are the vertices of $T$.
The condition that $V_1$ is totally nonnegative says that $\PP(V_1)$ is a point in $T$. 
The condition that $V_3$ is totally nonnegative says that the plane $\PP(V_3)$ crosses $T$ on the edges $\overline{e_1 e_2}$, $\overline{e_2 e_3}$, $\overline{e_3 e_4}$ and $\overline{e_4 e_1}$.

Thus, the intersection of $\PP(V_3)$ with $T$ is a quadrilateral whose sides, in cyclic order, are the intersections of $\PP(V_3)$ with $x_1=0$, $x_2=0$, $x_3=0$ and $x_4=0$; call this quadrilateral $Q$.
The vertices of $Q$, in cyclic order, are $y_{12} : = [ \Delta_{134} : \Delta_{234} : 0 : 0]$, $y_{23} : = [ 0 : \Delta_{124} : \Delta_{134} : 0]$, $y_{34} : = [ 0 : 0 : \Delta_{123} : \Delta_{124}]$, $y_{14} := [\Delta_{123} : 0:0: \Delta_{234} ]$.

We now want to ask when we can find a projective line $\PP(V_2)$, passing through $\PP(V_1)$ and lying in $\PP(V_3)$, with nonnegative Pl\"ucker coordinates.
The condition that $\PP(V_2)$ has nonnegative Pl\"ucker coordinates states that the line $\PP(V_2)$ passes through the quadrilateral $Q$ on the sides $\{ x_1=0 \}$ and $\{ x_4 = 0 \}$. 
The figure below shows what we expect to see inside the plane $\PP(V_3)$ if $(V_1, V_2, V_3)$ is a totally nonnegative flag; the gray shaded area is where the $x_i$ coordinates are positive. 

\centerline{
\begin{tikzpicture}
\node [label={[xshift=0.2cm, yshift=-0.6cm]$y_{12}$}] (y12) at (0,0) {};
\node [label={[xshift=-0.2cm, yshift=0cm]$y_{23}$}] (y23) at (1,2) {};
\node [label={[xshift=0.5cm, yshift=-0.5cm]$y_{34}$}] (y34) at (5,3) {};
\node [label={[xshift=0.4cm, yshift=-0.6cm]$y_{14}$}] (y14) at (5,0) {};
\filldraw [fill = gray] (y12.center) -- (y23.center)  -- (y34.center) -- (y14.center) -- (y12.center);
\draw [very thick] (-0.75, -0.3) -- (4.9, 3.5) node [pos=0.85, above] {$\PP(V_2)$};
\node [label={[xshift=0.75cm, yshift=-0.5cm]$\PP(V_1)$}] (p) at (2.0625, 1.6) {$\bullet$};
\draw [shorten >= -1.50cm, shorten <=-0.50cm] (y12.center) -- (y23.center) node [midway, left] {$x_4=0$};
\draw [shorten >= -0.50cm, shorten <=-2.50cm]  (y23.center) -- (y34.center)  node [pos=0.35, above] {$x_1=0$};
\draw [shorten >= -0.50cm, shorten <=-0.50cm] (y34.center) -- (y14.center) node [midway, right] {$x_2=0$};
\draw [shorten >= -1.50cm, shorten <=-1.50cm]  (y14.center) -- (y12.center) node [midway, below] {$x_3=0$};
\end{tikzpicture}}

Such a line $\PP(V_2)$ can be found if $\PP(V_1)$ is in the triangle with vertices $y_{12}$, $y_{23}$, $y_{34}$,  and cannot be found if $\PP(V_1)$ is in the triangle with vertices $y_{34}$, $y_{14}$, $y_{12}$. 
Some computation shows that $\PP(V_1)$ is in the former triangle if and only if $\Delta_1 \Delta_{234} - \Delta_2 \Delta_{134} \geq 0$.
\end{eg} 

\begin{Problem}
Suppose that we did make the naive definition in $\Fl(k_1, k_2, \ldots, k_r; n)$: Say that a partial flag is totally nonnegative if it is made from totally nonnegative subspaces, and decompose this semi-algebraic set according to which Pl\"ucker coordinates are zero and are nonzero. Would we still get a cell complex? A regular cell complex? Are there nice descriptions of the cells we get? See Bloch and Karp~\cite{BK1, BK2} for some partial answers to these questions.
\end{Problem}

Let $\Pio_u^w$ be a projected Richardson variety; by Proposition~\ref{OpenRichardsonModel}, we can choose $u$ and $w$ such that the projection $\pi : \Ro_u^w \to \Pio_u^w$ is an isomorphism.
We define $\Pio_u^{w, \geq 0}$ by $\Pio_u^{w, \geq 0} := \pi(\Ro_u^{w, \geq 0})$. 
So $\Pio_u^{w, \geq 0}$ is isomorphic to $\RR_{>0}^{\ell(w) - \ell(u)}$. 
The analogue of Theorem~\ref{thm:CellComplex} holds in the partial flag case as well.
This result was, likewise, proved by Galashin, Karp and Lam~\cite{GKL}.
\begin{theorem} \label{thm:ProjectedCellComplex}
$\Fl_n(k_1, k_2, \ldots, k_r)^{\geq 0}$ is a regular CW complex with cells $\Pio_u^{w, \geq 0}$. 
\end{theorem}

We close with some historical remarks.

\begin{remark}
We should emphasize the role of Lauren Williams in promoting the question of whether $\Fl_n^{\geq 0}$ and $G(k,n)^{\geq 0}$ are regular CW complexes.
She asked this question in many contexts, and proved many partial results towards it; we will describe several of them.

Together with Postnikov and Speyer, Williams showed that the cells of $G(k,n)^{\geq 0}$ form a CW complex~\cite{PSW}; Rietsch and Williams then generalized this result to all partial flag manifolds~\cite{RW1}.
Rietsch and Williams also showed~\cite{RW2} that the closure of any cell in a totally nonnegative partial flag manifold is contractible.

Williams~\cite{WilliamsShelling} also showed that the poset of cells in totally nonnegative partial flag manifolds is thin and shellable, which implies that it is so-called CW-poset~\cite{Bjorner}. 
A CW-poset is a poset which is the poset of cells in some regular CW-complex.
\end{remark}

\begin{remark} \label{PostnikovHistory}
Independent of Lusztig, Rietsch and the other researchers cited here, Alex Postnikov studied the structure of $G(k,n)^{\geq 0}$.
His results were described in a manuscript which was privately circulated in 2002, posted to the arXiv in 2006 and, as of January 2024, has not been submitted for publication~\cite{Post}.
We caution the reader that the version on Postnikov's website, which we have cited in the bibliography, is more up to date than the version on the arXiv.
Postnikov coined the term ``positroid" for the combinatorial objects indexing the cells of $G(k,n)^{\geq 0}$; we will discuss positroids further in Sections~\ref{sec:Positroid} and~\ref{sec:Plabic}.
\end{remark}

\section{Positroids} \label{sec:Positroid}

In this section, we will discuss combinatorial properties of projected Richardson varieties in $G(k,n)$, which are also known as \newword{positroid varieties}.
This subject was pioneered by Postnikov~\cite{Post}, but we will start by following the approach of Knutson, Lam and Speyer~\cite{KLS2}.
In the following section, we will discuss positivity for positroids, at which point we will follow Postnikov more closely.
See Remark~\ref{PostnikovHistory} for more on the history of Postnikov's work.
%
\subsection{What follows from the general $G/P$ theory}
Let $\preceq_k$ be the partial order on $S_n$ where $u \preceq_k w$ if there is a chain $u = v_0 \covered v_1 \covered v_2 \covered \cdots \covered v_{\ell} = w$ with $v_0[k] \prec v_1[k] \prec \cdots \prec v_{\ell}[k]$.
This is the $P$-Bruhat order  in the case $P = S_k \times S_{n-k}$; we call $\preceq_k$ the $k$-Bruhat order on $S_n$.
See Section~\ref{sec:ProjRichBackground} for several basic examples of $k$-Bruhat order.

Let $\tcQ(k,n)$ be the set of ordered pairs $(u,w)$ in $S_n \times S_n$ with $u \preceq_k w$. We partially order $\tcQ(k,n)$ by reverse containment; in other words, $(u,w) \preceq (u', w')$ if $u \preceq u' \preceq w' \preceq w$.
Let $\cQ(k,n)$ be the quotient of $\tcQ(k,n)$ by the equivalence relation that $(u,w) \sim (u', w')$ if there is an $x \in S_k \times S_{n-k}$ such that $u' = ux$ and $w' = wx$, and partially order $\cQ(k,n)$ by the relation that the equivalence class of $(u_1,w_1)$ is $\preceq$ the equivalence class of $(u_2, w_2)$ if there are $(u'_1, w'_1) \sim (u_1, w_1)$ and $(u'_2, w'_2) \sim (u_2, w_2)$ with $u'_1 \preceq u'_2 \preceq w'_2 \preceq w'_1$.
One can also embed $\cQ(k,n)$ as a subposet as $\tcQ(k,n)$ -- it is the subposet of pairs $(u,w)$ where $w$ is minimal in the coset $w (S_k \times S_{n-k})$. 

The projection map $\pi : R_u^w \to G(k,n)$ is birational onto its image if and only if $u \preceq_k w$, and its image depends only on the equivalence class of $(u,w)$ in $\cQ(k,n)$. 
The image of this projection is the projected Richardson variety $\Pi_u^w$, which we also call a \newword{positroid variety} in this context. 
In this case, the image $\pi(\Ro_u^w)$ also  depends only on the equivalence class of $(u,w)$; we denote it by $\Pio_u^w$ and call it the \newword{open positroid variety}.
We have $G(k,n) = \bigsqcup_{(u,w) \in \cQ(k,n)} \Pio_u^w$ and $\Pi_u^w = \bigsqcup_{(u', w') \succeq (u,w)} \Pio_{u'}^{w'}$.
Again, all of this is general properties of projected Richardsons, and we have not yet invoked anything special to the Grassmannian. 


We can also make statements about the positive real points: $\Pio_u^w \cap G(k,n)_{\geq 0}$ is an open ball of dimension $\ell(w) - \ell(u)$ and $\Pi_u^w \cap G(k,n)_{\geq 0}$ is a closed ball of this dimension~\cite{GKL}. The cells $\Pio_u^w \cap G(k,n)_{\geq 0}$ form a regular CW decomposition of $G(k,n)_{\geq 0}$. 
Here we are again saying things which are true for all projected Richardsons, but the Grassmannian case is nicer in a key way: $G(k,n)_{\geq 0}$ can be defined simply as the set of $k$-planes with nonnegative Pl\"ucker coordinates, whereas this definition does not work in other partial flag manifolds, see Example~\ref{eg:Flag13Positive}.

We now turn to things which are special to the Grassmannian.

\subsection{Affine permutations} \label{sec:AffinePerms}

We define the affine symmetric group $\tS_n$ to be the group of bijections $f : \ZZ \to \ZZ$ satisfying 
\[ f(i+n) = f(i)+n . \]
The group operation is composition.

There is a group homomorphism $\displace: \tS_n \to \ZZ$ given by
\[ \displace(f) = \frac{1}{n} \sum_{j=1}^n (f(j)-j). \]
We'll call $\displace(f)$ the \newword{displacement} of $f$.
We'll write $\tS^k_n$ for the set of permutations of displacement $k$.
So $\tS^0_n$ is a normal subgroup of $\tS_n$ and $\tS^k_n$ is a coset for $\tS^0_n$.

The group $\tS^0_n$ is a Coxeter group, with generators $\ts_1$, $\ts_2$, \dots, $\ts_n$ defined by
\[ \ts_i(j) = \begin{cases}
j+1 & j \equiv i \bmod n \\
j-1 & j \equiv i+1 \bmod n \\
j & \text{otherwise} \\
\end{cases}. \]
In this way, $\tS^0_n$ is the affine Coxeter group of type $\tA_{n-1}$.
We embed $S_n$ (which is the Coxeter group of type $A_{n-1}$) into $\tS^0_n$ by $s_i \mapsto \ts_i$.

\begin{remark}
The Coxeter group $\tA_{n-1}$ has rank $n$, not $n-1$. The notation should really be $\widetilde{A_{n-1}}$; it means ``take the root system $A_{n-1}$ and apply the construction which makes an affine Coxeter group", \textbf{not} ``the  rank $n-1$ Coxeter group of type $\tA$".
\end{remark}

For any integer $k$, define $\zeta_k(i) := i+k$. 
We use $\zeta_k$ to transfer the Bruhat order on $\tS^0_n$ to a partial order on $\tS^k_n$: For $f$ and $g \in \tS^k_n$, we define $f \preceq g$ if and only if $\zeta_k^{-1} f \preceq \zeta_k^{-1} g$ in the Bruhat order on $\tS^0_n$.
(Conjugation by $\zeta_k$ is an automorphism of Bruhat order on $\tS^0_n$, so we would get the same partial order if we multiplied on the right.)
We gave a model for Bruhat order on $S_n$ earlier in terms of rank matrices; we will give a similar model for Bruhat order on $\tS_n$ in Section~\ref{sec:CyclicRank}. 
First, though, we explain why we have introduced $\tS_n$.

We define a \newword{bounded affine permutation} to be a permutation $f$ in $\tS^k_n$ satisfying $i \leq f(i) \leq i+n$ for all $i$, and we order the bounded affine permutations by the induced order from $\tS^0_n$. 
We write $\Bound(k,n)$ for the subposet of bounded affine permutations in $\tS_n^k$.
Let $\omega_k$ be the affine permutation
\[ \omega_k(i) = \begin{cases} i+n & i \equiv 1,2,\ldots,k \bmod n \\ i & i \equiv k+1, k+2, \ldots, n \bmod n \end{cases}. \] 
\begin{prop}
The formula $(u,w) \mapsto u \omega_k w^{-1}$ is an isomorphism of posets from $\cQ(k,n)$ to the bounded affine permutations in $\tS^k_n$. To be more precise, if $(u_1,w_1)$ and $(u_2, w_2)$ are pairs in $\tcQ(k,n)$ which are equivalent in $\cQ(k,n)$, then $u_1 \omega_k w_1^{-1} = u_2 \omega_k w_2^{-1}$, and the resulting map $\cQ(k,n) \to \tS^k_n$ is an isomorphism of posets onto its image. 
\end{prop}

\begin{proof} See \cite[Theorem 3.16]{KLS2}. \end{proof}

\begin{eg}
In the first diagram, we have depicted the Hasse diagram of $\cQ(1,3)$; recall that the one line notation $a_1 a_2 \cdots a_n$ means that $i \mapsto a_i$ for $1 \leq i \leq n$.
We have listed all pairs in each equivalence class. (Compare  to Example~\ref{eg:kBruhat}.)
\[
\xymatrix{
(123, 123) \equiv (132,132) \ar@{-}[d] \ar@{-}[dr] & (213,213) \equiv (231, 231) \ar@{-}[dl] \ar@{-}[dr]   & (312,312) \equiv (321, 321) \ar@{-}[dl] \ar@{-}[d]  \\
(123, 213) \equiv (132, 231)  \ar@{-}[dr] & (132, 312)  \ar@{-}[d] & (213, 312) \equiv (231, 321)  \ar@{-}[dl] \\
& (123, 312) \equiv (132, 321) & \\
}
\]
In the second diagram, we have  depicted the corresponding elements of $\tS^1_3$. Here we use ``window notation": $[b_1 b_2 \cdots b_n]$ means that $i+rn \mapsto b_i+rn$ for all $1 \leq i \leq n$ and all $r \in \ZZ$. 
\[
\xymatrix@C=1.4in{
\ar@{-}[d] \ar@{-}[dr] [423] &  \ar@{-}[dl] \ar@{-}[dr]  [153] &  \ar@{-}[dl] \ar@{-}[d] [126] \\
 \ar@{-}[dr] [243] &        \ar@{-}[d] [324]         &      \ar@{-}[dl] [135] \\
& [234] & \\
}
\]
\end{eg}

\begin{remark}
Postnikov~\cite{Post} considers ``decorated" permutations: A \newword{decorated permutation} is a permutation $w$ in $S_n$ each of whose fixed points are colored either $-1$ or $+1$. 
There is an easy bijection between decorated permutations and $\bigsqcup_{k=0}^n \Bound(k,n)$: The bounded affine permutation $f: \ZZ \to \ZZ$ is sent to the permutation of $\ZZ/n \ZZ$ given by reducing modulo $n$, and we color each fixed point $i$ with $-1$ or $1$ according to whether $f(i) =i$ or $f(i) = i+n$. 
The statistic which Postnikov calls ``anti-exceedances" is our $k$.
The observation that decorated permutations should be encoded as affine permutations is due to Knutson, Lam and Speyer~\cite{KLS2}.
\end{remark}

\subsection{Cyclic rank matrices} \label{sec:CyclicRank}

As we described in Remark~\ref{rem:RankMatrices}, the Bruhat order on $S_n$ can be described in terms of ``rank matrices", which keep track of the cardinalities $\# \left( w[i] \cap [j] \right)$ for $w \in S_n$.
In this section, we will give a combinatorial model for Bruhat order on $\tS_n$ and relate it to the geometry of the Grassmannian.
The combinatorial description of Bruhat order on $\tS^0_n$ is due to Bj\"orner and Brenti~\cite{BB1, BB2}; 
the connection to geometry is due to Knutson, Lam and Speyer~\cite{KLS2}.
Note that the group which is called $\tS_n$ in~\cite{BB1} is the group which we call $\tS^0_n$.

Let $f$ be a permutation in $\tS^k_n$.
For any integers $i$, $j$, set
\[ r_{ij}(f) =  k - \# \{ a < i : f(a) > j \} . \]

\begin{prop}
For $f$ and $g \in \tS^k_n$, we have $f \preceq g$ if and only if $r_{ij}(f) \geq r_{ij}(g)$ for all $i$, $j$.
\end{prop}

\begin{proof}
See~\cite[Theorem 8.3.7]{BB2}. That theorem works with $f[i,j] = \# \{  a  \leq i : f(a) \geq j \}$, but the translation to $r_{ij}(f)$ is easy. 
\end{proof}

The following results are all straightforward, and can be found in~\cite[Section 3]{KLS2}.

\begin{prop}
We can recover $f$ from the array $r_{ij}(f)$ as follows: We have $f(i) = j$ if and only if $r_{ij}(f)+1 = f_{(i-1)j}(f) = r_{i(j+1)}(f)  = r_{(i-1)(j+1)}(f)$.
\end{prop}

\begin{prop} \label{PermutationConditions}
Given an array $r_{ij}$ of integers, this array is $r_{ij}(f)$ for a permutation $f \in \tS^k_n$ if and only if
\begin{enumerate}
\item For all $(i,j)$, we have $r_{ij} \leq r_{i(j+1)} \leq r_{ij}+1$ and $r_{ij} \leq r_{(i-1)j} \leq r_{ij}+1$.
\item We have $r_{(i+n)(j+n)} = r_{ij}$.
\item For all $(i,j)$, if $r_{ij} = r_{i(j+1)} = r_{(i-1)j}$ then $r_{(i-1)(j+1)} = r_{ij}$.
\item There is an integer $B$ such that $r_{ij} = j-i+1$ for $j < i-B$ and $r_{ij} = k$ for $j > i+B$.
\end{enumerate}
\end{prop}

\begin{prop} \label{BoundedCondition}
Let $f \in \tS^k_n$. Then $f$ is a bounded affine permutation, if and only if $r_{ij}(f) =0$ for $j \leq i-1$ and $r_{ij}(f) = k$ for $j \geq i+n-1$. 
\end{prop}

In other words, $\Bound(k,n)$ is in bijection with arrays of integers obeying the conditions~ of Propositions~\ref{PermutationConditions} and~\ref{BoundedCondition}.
We will call such an array a \newword{cyclic rank matrix}. 
Because of Proposition~\ref{BoundedCondition}, when we draw cyclic rank matrices, we will only depict the entries $r_{ij}(f)$ for $i-1 \leq j \leq i+n-1$, so the top entry of each column will be $0$ and the bottom entry will be $k$.

\begin{eg} \label{G24Example}
Let $k=2$, $n=4$ and consider the bounded affine permutation
\[ f(i) = \begin{cases}
 i & i \equiv 1 \bmod 4 \\ 
  i+2 & i \equiv 2 \bmod 4 \\ 
  i+4 & i \equiv 3 \bmod 4 \\
  i+2 & i \equiv 4 \bmod 4 \\
  \end{cases}. \]

In the array below, we depict $r_{ij}(f)$ for $1 \leq i,j \leq 8$. The index $i$ enumerates the columns, numbered from left to right; the index $j$ enumerates the rows, numbered from top to bottom. 
(This convention is transpose to that of~\cite{KLS2}, but it matches our convention of drawing permutation matrices with a $1$ in position $(f(i), i)$.)
The positions $(f(i), i)$ corresponding to the affine permutation are boxed. So the number in position $(i,j)$ is $k=2$ minus the number of boxed entries strictly to the lower left of $(i,j)$.

\[ \begin{bmatrix}
\fbox{0} & 0 &&& &&\phantom{\fbox{2}}&\phantom{\fbox{2}}\\
1 & 1 &0&& &&&\\
2&2&1&0 &&&\\
2&\fbox{2}&2&1&0 &&\\
2& 2&2&1&\fbox{0}&0 &&\\
2&2& 2& \fbox{1} & 1 & 1 & 0 &\\
2&2& \fbox{2} & 2&2&2&1&0 \\
2&2&2&2& 2&\fbox{2}&2&1 \\
\end{bmatrix} \]
\end{eg}

We now explain why we have chosen the particular formula for $r_{ij}(f)$ that we have.
\begin{theorem}
Let $(u,w) \in \cQ(k,n)$ and let $f \in \tS^k_n$ be the corresponding affine permutation.
Let $L$ be a $k$-plane in $G(k,n)$. 
The plane $L$ is in $\Pio_u^w$ if and only if, for all $i \leq j$, the coordinate projection of $L$ onto $\Span(e_i, e_{i+1}, \ldots, e_j)$ has dimension $r_{ij}(f)$.
The plane $L$ is in $\Pi_u^w$ if and only if, for all $i \leq j$, the coordinate projection of $L$ onto $\Span(e_i, e_{i+1}, \ldots, e_j)$ has dimension $\leq r_{ij}(f)$.
In both cases, the indices $i$, $i+1$, \dots, $j$ are taken modulo $n$.
\end{theorem}

\begin{proof} See \cite[Theorem 5.1]{KLS2}. \end{proof}

In short, \textbf{knowing which positroid cell $L$ lies in is equivalent to knowing the dimensions of the coordinate projection of $L$ onto $\Span(e_i, e_{i+1}, \ldots, e_j)$ for all $i \leq j$}, where indices are modulo $n$.

\begin{eg}
We continue studying the affine permutation $f$ in Example~\ref{G24Example}. 
Consider a $2$-plane $L$ which is the image of a $4 \times 2$ matrix with rows $\vec{v}_1$, $\vec{v}_2$, $\vec{v}_3$, $\vec{v}_4$.
We extend the notation $\vec{v}_i$ to be periodic in $i$ modulo $4$.
The conditions $r_{11}(f) = 0$ and $r_{02}(f) = 1$ say that, for $L$ in the corresponding positroid variety, we should have $\vec{v}_1 = \vec{0}$, and we should have $\dim \text{Span}(\vec{v}_0, \vec{v}_1, \vec{v}_2) = 1$. 
In other words, $\vec{v}_0$ must be parallel to $\vec{v}_2$.
The closed positroid variety $\Pi$ is the variety where these two conditions hold, and the open positroid variety $\Pio$ is the variety where these things occur and $\vec{v}_0$, $\vec{v}_2 \neq \vec{0}$. 
\end{eg}

\begin{remark}
If we asked more strongly to know the dimension of the coordinate projection of $L$ onto every coordinate subspace, this would be equivalent to studying the matroid associated to $L$.
In general, fixing a matroid $\cM$ and studying the set of all $L \in G(k,n)$ which realize $\cM$ gives horrible algebraic varieties, so it is surprising that looking only at the coordinate subspaces in consecutive positions gives nice varieties. The strata of $G(k,n)$ corresponding to matroids are sometimes called GGMS strata, as they were introduced in~\cite{GGMS}. For precise statements of ``GGMS strata can be horrible", see~\cite{Mnev1, Mnev2, SturmfelsMurphy1, SturmfelsMurphy2, Vakil, LeeVakil}.
\end{remark}

\subsection{Grassmann necklaces}
Grassmann necklaces were introduced by Postnikov~\cite[Section 16]{Post}. All the results in this section are easy and can be found in \cite[Sections 3 and 5]{KLS2}.

Let $r_{ij}$ be a cyclic rank matrix for a permutation in $\Bound(k,n)$.
As we go down column $i$, we have $0 = r_{i(i-1)} \leq r_{ii} \leq r_{i(i+1)} \leq \cdots \leq r_{i(i+n-1)} = k$, and $r_{ij} - r_{i(j-1)}  \in \{ 0,1 \}$ for each $j$. 
So there must be $k$ indices $j$ in $\{ i,i+1,\ldots, i+n-1 \}$ for which $r_{ij} > r_{i(j-1)}$.
Let $\tilde{I}_i$ be the set of such indices.


\begin{eg}
We continue pursuing Example~\ref{G24Example}. We have
\[ \tilde{I}_1 = \{ 2,3 \},\ \tilde{I}_2 = \{ 2,3 \},\ \tilde{I}_3 = \{ 3,4 \},\ \tilde{I}_4 = \{ 4,7 \} . \]
\end{eg}

Clearly, we can recover $r_{ij}(f)$ and hence $f$ from the $\tilde{I}_i$.
In fact, it is straightforward to go directly from the $\tilde{I}$ to the permutation $f$:
\begin{prop} \label{TildeNecklaceToPerm}
If $i \in \tilde{I}_i$, then $\tilde{I}_{i+1} = \tilde{I}_i \setminus \{ i \} \cup \{ f(i) \}$. If $i \not\in \tilde{I}_i$, then $f(i) = i$.
\end{prop}

\begin{prop} \label{TildeNecklaceCondition}
Let $\tilde{I}_i$ be a sequence of $k$-element sets of integers, with $\tilde{I}_i \subseteq \{ i,i+1, \ldots, i+n-1 \}$. 
Then $\tilde{I}_i$ corresponds to a permutation in $\Bound(k,n)$ if and only if, for all $i$, we have $\tilde{I}_i \setminus \{ i \}  \subseteq \tilde{I}_{i+1}$. 
\end{prop}

We have put a tilde over $\tilde{I}_i$ because the more standard thing to do is to reduce the elements of $\tilde{I}_i$ modulo $n$. 
Let $I_i$ be the subset of $[n]$ obtained by reducing $\tilde{I}_i$ modulo $n$. 
\begin{eg}
In our running example, we have
\[ I_1 = \{ 2, 3 \},\ I_2 = \{ 2,3 \},\ I_3 = \{ 3,4 \},\ I_4 = \{ 3,4 \} . \]
\end{eg}
Clearly, we can recover $\tilde{I}_i$ from $I_i$. Here are Propositions~\ref{TildeNecklaceToPerm} and~\ref{TildeNecklaceCondition} rewritten in terms of the $I_i$:
\begin{prop}
Let $(I_1, I_2, \ldots, I_n)$ correspond to a bounded permutation $f$. 
If $i \not\in I_i$, then $f(i) = i$ and $I_{i+1} = I_i$. 
If $i \in I_i$ and $I_i \neq I_{i+1}$, then $I_{i+1} = I_i \setminus \{ i \} \cup \{ f(i) \}$. 
If $i \in I_i$ and $I_i = I_{i+1}$, then $f(i) = i+n$.
\end{prop}

\begin{prop} \label{NecklaceCondition}
Let $(I_1, I_2, \ldots, I_n)$ be a sequence of $k$-element subsets of $[n]$. Then $(I_1, I_2, \ldots, I_n)$ corresponds to a permutation in $\Bound(k,n)$ if and only if, for all $i$, we have $I_i \setminus \{ i \} \subseteq I_{i+1}$ (including that $I_n \setminus \{ n \} \subseteq I_1$). 
\end{prop}

A sequence $(I_1, I_2, \ldots, I_n)$ of $k$-element subsets of $[n]$ obeying the conditions of Proposition~\ref{NecklaceCondition} is called a \newword{Grassmann necklace}.
So bounded affine permutations are also in bijection with Grassmann necklaces.

We now explain the geometric significance of Grassmann necklaces.
Let $L$ be a point of the Grassmannian $G(k,n)$. Then $L$ lies in a unique open positroid variety $\Pio$; let $(I_1, I_2, \ldots, I_n)$ be the corresponding Grassmann necklace.
The set $I_1$ encodes the ranks $r_{11}$, $r_{12}$, $r_{13}$, \dots.
Looking at the standard description of the Schubert decomposition of the Grassmannian, we see that $L$ is in the Schubert cell $\Xo_{I_1}$. 
More generally, let $\rho$ be the permutation $1 \mapsto 2 \mapsto 3 \mapsto \cdots \mapsto n \mapsto 1$ of $[n]$ and let $\sigma$ be the automorphism $e_1 \mapsto e_2 \mapsto e_3 \mapsto \cdots \mapsto e_n \mapsto e_1$ of $\AA^n$.
Then we similarly get the $L$ is in the permuted Schubert cell $\sigma^{j-1} \Xo_{\rho^{j-1}(I_j)}$.

In short, \textbf{knowing which positroid cell $L$ lies in is equivalent to which Schubert cell the spaces $L$, $\sigma(L)$, $\sigma^2(L)$, \dots, $\sigma^{n-1}(L)$ lie in.}

\begin{remark}
Instead of reading $r_{ij}(f)$ by columns, we could read by rows. Let $\tilde{J}_j = \{ i : r_{(i-1)j} = r_{ij} +1 \}$, so $\tilde{J}_j \subseteq \{ j-n+1, \ldots, j-1, j \}$; let $J_j$ be the reduction of $\tilde{J}_j$ modulo $n$.
We call $J_j$ the \newword{reverse Grassmann necklace} of $f$. We have $J_j \setminus \{ j \} \subseteq J_{j-1}$, and the reverse Grassmann necklace encodes which opposite Schubert cell the spaces $L$, $\sigma(L)$, $\sigma^2(L)$, \dots, $\sigma^{n-1}(L)$ lie in. 
\end{remark}

\subsection{The cohomology class of a positroid variety}

Before describing the cohomology class of a positroid variety, we should explain why we have said nothing about the cohomology class of a Richardson variety in $\Fl_n$. 
The cohomology (and also Chow) ring of $\Fl_n$ is the quotient of $\ZZ[x_1, x_2, \ldots, x_n]$ by the ideal generated by the positive degree homogeneous symmetric polynomials in the $x_i$, and it is a meaningful question to ask what class in $H^{\ast}(\Fl_n)$ corresponds to the fundamental class $[R_u^w]$.
This question is both straightforward and impossible.
The intersection of $X_u$ and $X^w$ is transverse, so $[R_u^w] = [X_u] [X^w]$.
The class of $X_u$ is represented by the Schubert polynomial $\mathfrak{S}_u(x_1, \ldots, x_n)$, and the class of $X^w$ is represented by the Schubert polynomial $\mathfrak{S}_{w_0 w}(x_1, \ldots, x_n)$. 
So the class of $[R_u^w]$ is represented by the product  $\mathfrak{S}_u(x_1, \ldots, x_n) \mathfrak{S}_{w_0 w}(x_1, \ldots, x_n)$. 

This is the sense in which the answer is straightforward. In a different sense, the answer is impossible: The usual goal of determining the cohomology class of a subvariety of $\Fl_n$ is to write it as a linear combination of Schubert classes: That is to say, to find the coefficients $c_v$ for which $[R_u^w] = \sum c_v [X_v]$. Equivalently, we want to write $\mathfrak{S}_u(x_1, \ldots, x_n) \mathfrak{S}_{w_0 w}(x_1, \ldots, x_n)$ as a sum $\sum c_v \mathfrak{S}_v(x_1, \ldots, x_n)$.
Multiplying Schubert polynomials, and expressing the result in the Schubert basis, is one of the most famous open problems in algebraic combinatorics.
This is why we have chosen not to discuss it here.

We can ask the same questions about the cohomology class of a positroid variety in $H^{\ast}(G(k,n))$. 
Here there are better answers! We will describe two of them, one by Knutson, Lam and Speyer~\cite{KLS2} and one by Bergeron and Sottile~\cite{BS2}.

Let $\Lambda$ be the ring of symmetric polynomials in infinitely many variables. 
The cohomology (and also Chow) ring of $G(k,n)$ is the quotient of $\Lambda$ by the ideal generated by $e_m(x)$ for $m > k$ and $h_m(x)$ for $m > n-k$. 
(Here $e$ and $h$, as usual, denote the elementary and the complete homogeneous symmetric polynomials.)
A basis of $H^{\ast}(G(k,n))$ is given by the Schur polynomials $s_{\lambda}(x_1, \ldots, x_k)$ where $\ell(\lambda) \leq k$ and $\ell(\lambda^T) \leq n-k$. 
Knutson, Lam and Speyer show that the cohomology class $[\Pi_f]$ is represented by the \newword{affine Stanley symmetric function} $F_f$, first introduced in~\cite{Lam06}.
We will now define the affine Stanley symmetric function.

We define an element $c$ of $\tS^0_n$ to be \newword{cyclically decreasing} if it has a reduced word $w = s_{i_1} s_{i_2} \cdots s_{i_{\ell}}$ such that (1) no letter appears more than once in $s_{i_1} s_{i_2} \cdots s_{i_{\ell}}$ and (2) if $s_i$ and $s_{i+1}$ both appear in $s_{i_1} s_{i_2} \cdots s_{i_{\ell}}$, then $s_{i+1}$ appears before $s_i$ (including that $s_1$ must precede $s_n$). There are $2^n-1$ cyclically decreasing elements of $\tS^0_n$, indexed by the proper subsets of $\{ s_1, s_2, \ldots, s_n \}$ because, for each proper subset of  $\{ s_1, s_2, \ldots, s_n \}$, there is a unique way to order it, up to commutation, to obey condition (2). We note that the identity element is considered cyclically decreasing.

Given $g \in \tS^0_n$, we define a \newword{cyclically decreasing factorization} of $g$ to be a sequence $c_1$, $c_2$, \dots, of cyclically decreasing elements, all but finitely many of which are the identity, such that $g = c_1 c_2 c_3 \cdots$ and $\ell(g) = \sum \ell(c_i)$. 
We define
\[ F_f(x) = \sum_{g = c_1 c_2 c_3 \cdots} \left( x_1^{\ell(c_1)} x_2^{\ell(c_2)} x_3^{\ell(c_3)} \cdots \right) \]
where the sum is over all cyclically decreasing factorizations of $f$. 

\begin{eg}
All cyclically decreasing factorizations of $s_1 s_2$ are of the form $(e,\ldots, e, s_1, e,\ldots, e, s_2, e,\ldots)$.
So $F_{s_1 s_2} = \sum_{i<j} x_i x_j = s_{11}(x)$. 
We will express this more briefly by omitting the $e$'s and using parentheses: All cyclically decreasing factorizations of $s_1 s_2$ have the form $(s_1)(s_2)$. 
By contrast, the element $s_2 s_1$ has cyclically decreasing factorizations of two forms: $(s_2)(s_1)$ and $(s_2 s_1)$. 
So $F_{s_2 s_1} = \sum_{i<j} x_i x_j + \sum_k x_k^2 = s_2(x)$. 
\end{eg}

\begin{eg} \label{nonSchurPositiveStanley}
This example is taken from~\cite[Section 7.4]{KLS2}.
Consider the element $g$ of $\tS^0_4$ with window notation $[-1, 4, 1, 6]$. We list the forms of the  cyclically decreasing factorizations of $g$ below:
\[\begin{array}{c}
 (s_1)(s_3)(s_2)(s_4),\  (s_1)(s_3)(s_4)(s_2),\  (s_3)(s_1)(s_2)(s_4),\  (s_3)(s_1)(s_4)(s_2) \\
  (s_1 s_3)(s_2)(s_4),\ (s_1 s_3)(s_4)(s_2),\ (s_1)(s_3 s_2)(s_4),\ (s_3)(s_1 s_4)(s_2),\ (s_1)(s_3)(s_2 s_4),\ (s_3)(s_1)(s_2 s_4) \\
 (s_1 s_3)(s_2 s_4) \\
 \end{array} \]
 The polynomial $F_g$ is given below:
 \[4 \sum_{i<j<k<\ell} x_i x_j x_k x_{\ell} +2 \sum_{i<j<k} x_i^2 x_j x_k + 2 \sum_{i<j<k} x_i x_j^2 x_k+2  \sum_{i<j<k} x_i x_j x_k^2 +  \sum_{i<j} x_i^2 x_j^2 . \]
 We compute that this
 \[ 4 m_{1111} + 2 m_{211} + m_{22} = s_{22} + s_{211} - s_{1111}. \]
 We note that $F_g$ is not necessarily Schur positive!
\end{eg}

We now state Theorem 7.1 of~\cite{KLS2}:
\begin{theorem}
Let $f$ be a bounded affine permutation in $\tS^k_n$ and let $g = \omega_k^{-1} f$. Then the symmetric polynomial $F_g$ represents the class $[\Pi_f]$ in $H^{\ast}(G(k,n))$.
\end{theorem}

\begin{eg}
Let $f$ be the permutation which is given in one line notation by $[1,6,3,8]$, for $(k,n) = (2,4)$. The corresponding positroid variety is the single point of $G(2,4)$ corresponding to $\Span(e_1, e_3)$.
Then $t_k^{-1} f$ is the permutation $w$ in Example~\ref{nonSchurPositiveStanley}. Although $F_g$ is not Schur positive, its image in $H^{\ast}(G(2,4))$ is $s_{22}$, which is the class of a single point, as claimed.
\end{eg}

\begin{eg}
In Example~\ref{eg:G24Complex1ForBS} we considered the positroid variety $\Delta_{23}=\Delta_{14} = 0$.
The bounded affine permutation is $[4,3,6,5]$, so $\omega_2^{-1} [4,3,6,5] = [2,1,4,3]$, which corresponds to $g = s_1 s_3$. 
The  cyclically decreasing factorizations are of the  forms $(s_1)(s_3)$, $(s_3)(s_1)$ and $(s_1 s_3)$, so
\[ F_w = 2 m_{11} + m_2 = s_{11} + s_2 . \]
Sure enough, the cohomology class of this positroid variety is $(s_1)^2 = s_{11} + s_2$.
\end{eg}

\begin{eg} \label{DivisorExample1}
In Example~\ref{eg:G24Complex2ForBS} we considered the positroid variety $\Delta_{14} = 0$.
The bounded affine permutation is $[3,4,6,5]$, so $\tau_2^{-1} [4,3,6,5] = [1,2,4,3]$, which corresponds to $g = s_3$. 
So $F_w = \sum x_i = s_1$. Sure enough, the cohomology class of this positroid variety is $s_1$.
\end{eg}

We now describe the work of Bergeron and Sotille.
Let $u \preceq_k w$ be permutations in $S_n$, and let $I$ be a $k$-element subset of $[n]$.
To compute the class of $\Pi_u^w$ in $H^{\ast}(G(k,n))$, we need to compute the intersection $[\Pi_u^w] \cap [X_{I}]$ in $H^{\ast}(G(k,n))$.
We will pull this back to an intersection in $\Fl_n$.
Since $\pi : R_u^w \to \Pi_u^w$ is birational, we can instead compute the intersection $[R_u^w] \cap [\pi^{-1}(X_{I})]$ in $H^{\ast}(\Fl_n)$. 
The preimage $\pi^{-1}(X_{I})$ is the Schubert variety $X_{v(I)}$ where  $v(I)$ is the unique permutation with $v(I)[k] = I$, such that $v(I)_1 < v(I)_2 < \cdots < v(I)_k$ and $v(I)_{k+1} < v(I)_{k+2} < \cdots < v(I)_n$.

So we want to compute $[X_u] \cap [X^w] \cap [X_{v(I)}]$ in $H^{\ast}(\Fl_n)$. 
Reassociating the product, we want to extract the coefficient of $\mathfrak{S}_w(x)$ in $\mathfrak{S}_u(x) \mathfrak{S}_{v(I)}(x)$.

The permutation $v(I)$ is a very particular type of permutation -- it has at most one descent. Such a permutation is called a \newword{Grassmannian permutation}.
This means that the corresponding Schubert polynomial $\mathfrak{S}_{v(I)}(x)$ is a Schur polynomial, namely, : $\mathfrak{S}_{v(I)}(x_1, x_2, \ldots, x_n) = s_{\lambda}(x_1, x_2, \ldots, x_k)$ where $\lambda_j = i_{k+1-j} - (k+1-j)$, for $I = \{ i_1 < i_2 < \cdots < i_k \}$.
So the corresponding combinatorial problem is to multiply the Schubert polynomial $\mathfrak{S}_u(x_1, x_2, \ldots, x_n)$ by the Schur polynomial $s_{\lambda}(x_1, x_2, \ldots, x_k)$. 

There is an extensive literature on multiplying Schuberts by Schurs, in which~\cite{BS2} is particularly relevant.
Other relevant work is~\cite{MPP}, giving positive rules for  multiplying Schuberts by Schurs in particular cases.

Given $u \preceq_k w$, Bergeron and Sottile compute a symmetric polynomial $S_{[u,w]_k}$  such that the coefficient of $s_{\lambda}$ in $S_{[u,w]_k}$ is  $[X_u] \cap [X^w] \cap [X_{v(I)}]$.
This polynomial is a sum over maximal chains $u=v_0 \prec_k v_1 \prec_k \cdots \prec_k v_m =w$ from $u$ to $w$ in the $\prec_k$ order. (Here $m  = \ell(w) - \ell(u)$.)
For each cover $v_{i-1} \prec_k v_{i}$, we have $v_{i} = (a_i b_i) v_{i-1}$ for some $1 \leq a_i < b_i \leq n$; we write this chain as $v_0 \xrightarrow{b_1} v_1 \xrightarrow{b_2} v_2 \xrightarrow{b_3} \cdots \xrightarrow{b_m} v_m$. 
Given a sequence $(b_1, b_2, \cdots b_m)$ of labels, we define the quasisymmetric function $Q_{b_1 b_2 \cdots b_m}$ by 
\[ Q_{b_1 b_2 \cdots b_m} = \sum_{\substack{i_1 \leq i_2 \leq \cdots \leq i_m  \\ b_j > b_{j+1} \implies i_j < i_{j+1}}} x_{i_1} x_{i_2} \cdots x_{i_m} . \]
Then 
\[ S_{[u,w]_k} = \sum_{u=v_0 \xrightarrow{b_1} v_1 \xrightarrow{b_2} v_2 \xrightarrow{b_3} \cdots \xrightarrow{b_m} v_m=w} Q_{b_1 b_2 \cdots b_m} \]
where the sum runs over all saturated chains from $u$ to $w$ in $\prec_k$.

\begin{eg}
We consider the interval $[2143, 3412]$ in $\prec_2$, depicted earlier in Example~\ref{eg:G24Complex1ForBS}.
The two maximal chains are $2143 \xrightarrow{3} 3142 \xrightarrow{4} 3412$ and  $2143 \xrightarrow{4} 2413 \xrightarrow{3} 3412$. 
So the Bergeron-Sottile polynomial is
\[ Q_{34} + Q_{43} = \sum_{i \leq j} x_i x_j + \sum_{i<j} x_i x_j = s_{2} + s_{11}. \]
The corresponding positroid variety $\Pi$ is $\Delta_{14} = \Delta_{23} = 0$, so $[\Pi] \cap s_2 = 1$ and $[\Pi] \cap s_{11}  = 1$ as claimed.
\end{eg} 

\begin{eg} \label{DivisorExample2}
We consider the interval $[2134, 3412]$ in $\prec_2$, from Example~\ref{eg:G24Complex2ForBS}.
The two maximal chains are $2134 \xrightarrow{3} 2314 \xrightarrow{4} 2413 \xrightarrow{3} 3412$ and  $2134 \xrightarrow{3} 3124 \xrightarrow{2} 3214 \xrightarrow{4} 3412$.
So the Bergeron-Sottile polynomial is
\[ Q_{343} + Q_{324} = \sum_{i \leq j<k} x_i x_j x_k + \sum_{i<j \leq k}  x_i x_j x_k = s_{21} . \]
The positroid variety is $\Delta_{14} = 0$, which has cohomology class $s_{1}$. 
So $[\Pi] \cap s_{21}  =1$, as claimed.
\end{eg} 

We note that we considered the same positroid variety in Examples~\ref{DivisorExample1} and~\ref{DivisorExample2}, but the corresponding symmetric polynomials are of different degrees, $3$ and $1$.
The degree of the affine Stanley polynomial corresponds to the codimension of the positroid variety, and the degree of the Bergeron-Sottile polynomial corresponds to the dimension of the positroid variety.

\section{Plabic graphs} \label{sec:Plabic}
We have earlier described Deodhar's parametrization of Richardson varieties.
In the Grassmannian case, Postnikov has discovered a different, much more flexible, way to parametrize positroid varieties.
The theory presented in this section was pioneered by Postnikov~\cite{Post}.

\subsection{Plabic graphs and the boundary measurement map with positive real weights}
Let $D$ be a closed disc in $\RR^2$ and let $\partial D$ be its boundary.
Let $G$ be a bipartite graph, with a  black and white coloring, embedded in $D$.
We write $\partial G$ for the vertices of $G$ on $\partial D$ and call these the \newword{boundary vertices}; we call the other vertices  \newword{interior vertices}.
We will assume that the boundary vertices are entirely black and we fix a numbering of them by $[n]$ in clockwise order.
Let there be $m+k$ white interior vertices and $m$ black interior vertices.
We call such a graph $G$ a \newword{plabic graph}.

A \newword{perfect matching} of $G$ is a collection of edges $M$ of $G$ such that every interior vertex lies on exactly one edge of $M$ and every boundary vertex lies on at most one edge of $M$.
So there are $k$ boundary vertices covered by $M$; we write $\partial M$ for the corresponding subset of $[n]$ and call $\partial M$ the \newword{boundary of $M$}. 
We will be interested in counting perfect matchings with fixed boundary. More generally, we will have a weighting function that assigns a nonzero \newword{weight} $w(e)$ to each edge $e$. 
We define the \newword{weight of a perfect matching $M$} to be $\prod_{e \in M} w(e)$, and, for $I$ a $k$-element subset of $[n]$, we put
\[ D_I := \sum_{\partial M=I} w(M). \]
Evaluating $D_I$ when all the weights are $1$ is, thus, counting the perfect matchings of $G$ with boundary $I$, and evaluating $D_I$ is general is performing a weighted count of such matchings. 
Combinatorialists have long been interested in this sort of problem, we mention \cite{Ciucu, Kenyon, Propp} as some relevant surveys.

\begin{eg} \label{eg:BasicG24}
Consider the graph below, with $k=2$ and $n=4$.
\bigskip

\centerline{
\begin{tikzpicture}
\draw[fill=none](0,0) circle (2.0) ; 
\draw (2,0) -- (0,1) node [midway,above] {$q$};
\draw (2,0) -- (0,-1) node [midway,below] {$r$};
\draw (-2,0) -- (0,1)  node [midway,above] {$p$};
\draw (-2,0) -- (0,-1) node [midway,below] {$s$};
\draw (0,2) -- (0,1)  node [midway,right] {$t$};
\draw (0,-2) -- (0,-1)  node [midway,left] {$u$};
\node (1) at (-2.5,0) {$1$};
\node (2) at (0,2.5) {$2$};
\node (3) at (2.5,0) {$3$};
\node (4) at (0,-2.5) {$4$};
\draw[fill=black](2,0) circle (0.1) ; 
\draw[fill=black](-2,0) circle (0.1) ; 
\draw[fill=black](0,2) circle (0.1) ; 
\draw[fill=black](0,-2) circle (0.1) ; 
\draw[fill=white](0,1) circle (0.1) ; 
\draw[fill=white](0,-1) circle (0.1) ; 
\end{tikzpicture}}

Then we have
\[ D_{12} = st,\ D_{13} = pr+qs,\ D_{14} = pu,\ D_{23} = rt,\ D_{24} = tu,\ D_{34} = qu . \]
\end{eg}

\begin{remark}  \label{multicoloredBoundary}
The assumption that the boundary vertices are all black is for convenience. 
More generally, let us continue to assume that $G$ is a bi-colored graph embedded in $D$, and continue to use the term ``boundary vertices" for vertices of $G$ on $\partial M$ and ``interior vertices" for the other vertices, but now allow the boundary vertices to be colored both black and white. 
Let $\partial_B(G)$ be the set of black boundary vertices and let $\partial_W(G)$ be the set of white boundary vertices. 
We suppose that there are $p$ white boundary vertices, $n-p$ black boundary vertices, $w$ white interior vertices and $b$ black interior boundary vertices and we define $k = w-b+p$.
For a perfect matching $M$, we define
\[ \partial(M) = \{ v \in \partial_W(G) : v \ \text{not covered by $M$} \} \cup \{ v \in \partial_B(G) : v \ \text{covered by $M$} \} .\]
Then it is an easy exercise to check that $|\partial(M)| = k$. 

All the results in this section remain true with this more general definition. This added generality is not usually useful, because one can reduce to the assumption that all boundary vertices are black in the following manner: Let $G$ be a graph as above. Define a new graph $G'$ as follows: For each white boundary vertex $v$ of $G$, nudge $v$ into the interior of $G$ and add a new black boundary vertex $v'$ with an edge from $v$ to $v'$. If we are working with weighted edges, give this new edge weight $1$. Then the perfect matching of $G$ are in easy bijection with the perfect matchings of $G'$, and this bijection preserves the weight.
\end{remark}

The connection to the Grassmannian comes from the following theorem:
\begin{theorem} \label{thm:BoundaryPlucker}
For any values of the weights $w(e)$ in the ground field $\kappa$, either all of the $D_I$ are zero, or else there is a $k$-plane $L$ in $\AA^n$ with $\Delta_I(L) = D_I$ for all $I$.
\end{theorem}

\begin{eg}
We observe that the values of $D_I$ in Example~\ref{eg:BasicG24} obey the Pl\"ucker relation
\[ (pr+qs)(tu) = (st)(qu) + (pu)(rt). \]
An explicit $2$-plane with these Pl\"ucker coordinates is the image of
\[ \begin{bmatrix}  p&-s \\ t&0 \\ q&r \\ 0&u \\ \end{bmatrix}. \]
\end{eg}

 Theorem~\ref{thm:BoundaryPlucker} is implicit in  Kasteleyn's  ``permanent-determinant" method for computing $D_I$~\cite{Kast}.
 Propp~\cite[Section 6]{ProppDominoShuffle} and Kuo~\cite{Kuo} note that the $D_I$ obey the $3$-term Pl\"ucker relations but don't label them as $3$-term Pl\"ucker relations.
 Postnikov, Speyer and Williams~\cite{PSW}, building on Talaska~\cite{Talaska}, point on that Postnikov's parametrization of $G(k,n)^{\geq 0}$~\cite{Post}, can be written in terms of dimers.
The earliest sources to directly state (and prove) Theorem~\ref{thm:BoundaryPlucker} may be Lam~\cite{LamNotes} and Speyer~\cite{SpeyerKast}.
While the author may be biased here, he suggests that~\cite{SpeyerKast} is a clear presentation and makes clear the sense in which this result is essentially implicit in Kasteleyn's work.

To avoid degeneracies, we assume from now on that $G$ has at least one perfect matching, so the $D_I$ are not all identically zero.

We write $\Edge(G)$ for the set of edges of $G$. Thus, the map $\mu$ sending the set of weights $(w(e))_{e \in \Edge(G)}$ to the $\binom{n}{k}$ values $D_I$ is essentially a map $\GG_m^{\Edge(G)} \to G(k,n)$. 
We write ``essentially'' because it is possible that all of the $D_I$ could be zero for some weights; we will return to this issue in Section~\ref{sec:twist}.
Since the polynomials defining $\mu$ are sums of monomials with nonnegative coefficients, the map $\mu$ takes $\RR_{> 0}^{\Edge(G)}$ to $G(k,n)_{>0}$ and, if there is a matching $M$ with $\partial(M) = I$, then the generating function $D_I$ is positive everywhere on $\RR_{> 0}^{\Edge(G)}$.
Following Postnikov, we call $\mu$ the \newword{boundary measurement map}.
We now begin to describe Postnikov's results.
\begin{Theorem} \cite{Post}
Let $G$ be any bipartite planar graph as above. Then there is a positroid cell $\Pio$ in $G(k,n)$ such that the image of $\mu$ is $\Pio_{>0}$. 
\end{Theorem}

The map $\mu: \RR_{>0}^{\Edge(G)} \longrightarrow \Pio_{>0}$ is rarely bijective, for the following reason. Let $v$ be any internal vertex of $G$. If we rescale all of the edges incident on $v$ by a common factor $t$, then every $D_I$ is multiplied by $t$.
Since the $D_I$ are homogeneous coordinates on $\Pio$, rescaling the edge weights in this way will not change the output of $\mu$.
So, whenever $G$ has internal vertices, the fibers of $\mu$ have positive dimension.

To address this issue, let $\Vertices(G)$ be the set of internal vertices of $G$. 
We will define two weights $w$ and $w'$ in $\GG_m^{\Edge(G)}$ to be \newword{gauge equivalent} if there is a function $t : \Vertices(G) \to \GG_m$ such that, for every edge $e$ with endpoints $(v,w)$, we have $w'(e) = t(v) t(w) w(e)$. 
Let $T$ be the torus of weight functions modulo gauge equivalence.
The dimension of $T$ is one less than the number of faces of $G$~\cite[Lemma 11.1]{Post}. 
There are several natural ways to coordinatize $T$; see  \cite[Lemma 11.2]{Post} and \cite[Proposition 5.5]{MullerSpeyerTwist} for two valuable choices.

We write $T_{>0}$ for positive valued edge weights modulo gauge equivalence.
Then $\mu$ induces a map $T_{>0} \to G(k,n)_{\geq 0}$, which we will also call $\mu$.

\begin{eg} \label{eg:BasicG24GaugeFixed}
In Example~\ref{eg:BasicG24}, we can use gauge transformations to uniquely set $t=u=1$. So we can consider $(p,q,r,s)$ as coordinates on $T$. (Note that this means that $T$ is $4$-dimensional, one less than the number of faces of $G$.) 
The map $T \to G(2,4)$ sends $(p,q,r,s)$ to the image of 
\[ \begin{bmatrix}  p&-s \\ 1&0 \\ q&r \\ 0&1 \\ \end{bmatrix}. \]
In other words, we have $(\Delta_{12}, \Delta_{13}, \Delta_{14}, \Delta_{23}, \Delta_{24}, \Delta_{34}) = (s, pr+qs, p, r, 1, q)$
If $(p,q,r,s)$ are positive reals, this a diffeomorphic parametrization of the big positroid cell in $G(2,4)_{>0}$ by $\RR_{>0}^4$.
If we consider complex points (or real points without positivity, or points over some other field), then this is an open inclusion into the largest positroid cell; its image is the open set $\{ \Delta_{24} \neq 0\}$ in this cell.
\end{eg}

\begin{eg} \label{eg:G36Hexagon}
Consider the plabic graph shown below with $k=3$ and $n=6$:

\centerline{
\begin{tikzpicture}
\draw[fill=none] (0,0) circle (2.0) ; 
\node [label={[label distance=0.2cm]180:$1$}] (1) at (-2,0) {};
\node  [label={[label distance=0.2cm]120:$2$}] (2) at (-1, 1.732) {};
\node  [label={[label distance=0.2cm]60:$3$}] (3) at (1, 1.732) {};
\node  [label={[label distance=0.2cm]0:$4$}] (4) at (2,0) {};
\node  [label={[label distance=0.2cm]300:$5$}] (5) at (1, -1.732) {};
\node   [label={[label distance=0.2cm]240:$6$}] (6) at (-1, -1.732) {};
\node (p) at (0.5,0.866) {};
\node (q) at (0.5,-0.866) {};
\node (r) at (-1,0) {};
\draw (2.center)--(p.center) node [midway,above] {$x_2$}; \draw (3.center)--(p.center); \draw (4.center)--(p.center) node [midway,above] {$x_3$};
\draw (4.center)--(q.center) node [midway,below] {$x_4$}; \draw (5.center)--(q.center); \draw (6.center)--(q.center) node [midway,below] {$x_5$};
\draw (6.center)--(r.center) node [midway,left] {$x_6$}; \draw (1.center)--(r.center); \draw (2.center)--(r.center) node [midway,left] {$x_1$};
\draw[fill=white] (p) circle (0.1) ;  \draw[fill=white] (q) circle (0.1) ; \draw[fill=white] (r) circle (0.1) ; 
\draw[fill=black] (1) circle (0.1) ;  \draw[fill=black] (2) circle (0.1) ; \draw[fill=black] (3) circle (0.1) ; \draw[fill=black] (4) circle (0.1) ;  \draw[fill=black] (5) circle (0.1) ; \draw[fill=black] (6) circle (0.1) ; 
\end{tikzpicture}}

We can use the gauge transformations to normalize the unlabeled  edges to $1$. The Pl\"ucker coordinates, computed using the boundary measurement formula, are as follows:
\[ \Delta_{123}=\Delta_{345}=\Delta_{156}=0 \] 
\[ \Delta_{234}= x_1 x_4  , \quad \Delta_{456}= x_3 x_6    , \quad  \Delta_{126}= x_2 x_5   \] 
\[ \Delta_{124}=x_2 x_4 , \quad \Delta_{245} = x_1 x_3, \quad \Delta_{346}=x_4 x_6 , \quad \Delta_{146} = x_3 x_5, \quad \Delta_{256}=x_2 x_6 , \quad \Delta_{236} = x_1 x_5\]
\[ \Delta_{134} = x_4 , \quad \Delta_{136} = x_5 , \quad \Delta_{356} = x_6 , \quad \Delta_{235} = x_1 , \quad \Delta_{125} = x_2 , \quad \Delta_{145}  = x_3 \]
\[ \Delta_{135} = 1 , \quad \Delta_{246} = x_1 x_3 x_5 + x_2 x_4 x_6 \]
This point of $G(3,6)$ is realized as the image of the matrix
\[ \begin{bmatrix} 
1&0&0\\
x_1&x_2&0 \\
0&1&0 \\
0&x_3&x_4 \\
0&0&1\\
x_6&0&x_5\\
\end{bmatrix} .\]
If $x_1$ through $x_6$ are positive reals, then this is a diffeomorphic parameterization of the positroid cell in $G(3,6)_{>0}$ corresponding to the Grassmann necklace $(124, 234, 346, 456, 562, 612)$, with corresponding bounded affine permutation $f(i) = i+2$ for $i$ odd and $f(i) = i+4$ for $i$ even.

The image is the open set $\{ \Delta_{135} \neq 0 \}$ in this cell; the complement of this open set in $\Pio$ consists of points where rows $1$, $3$ and $5$ are parallel to each other.
\end{eg}

Example~\ref{eg:BasicG24GaugeFixed} and~\ref{eg:G36Hexagon} are both examples of \newword{reduced plabic graphs}.
\begin{TD} \cite{Post} \label{ReducedMotivation}
Let $G$ be a plabic graph, and let $\Pio$ be the corresponding positroid variety. Then the following are equivalent:
\begin{enumerate}
\item We have $\dim T = \dim \Pio$. In other words, the number of faces of $G$ is $\dim \Pio +1$.
\item The boundary measurement map $\mu : T_{>0} \to \Pio_{>0}$ is bijective.
\item The boundary measurement map $\mu : T_{>0} \to \Pio_{>0}$ is a diffeomorphism.
\end{enumerate}
In this case, we call $G$ \newword{reduced}.
\end{TD}

The word ``reduced" is meant by analogy to ``reduced word" and is supposed to suggest an analogy between the use of reduced words to parametrize various objects in the flag manifold and the use of reduced plabic graphs to parametrize positroid cells.
In the author's opinion, the analogy is most close to the parametrization of totally positive unipotent cells, as discussed in Remark~\ref{rem:Unipotent}.
As discussed in that remark, let $U_-$ be the group of lower triangular matrices with $1$'s on the diagonal. Let $U_-^{\geq 0}$ be the sub-semigroup of matrices in $U_-$ with nonnegative minors. 
We have the following analogies:
\begin{itemize}
\item $U_-^{\geq 0}$ is stratified into $n!$ pieces indexed by the symmetric group (see Remark~\ref{rem:Unipotent} and \cite{FominShapiro}); $G(k,n)^{\geq 0}$ is indexed into pieces indexed by bounded affine permutations.
\item Given any permutation $w$ and any word  $s_{i_1} s_{i_2} \cdots s_{i_a}$ with Demazure product $w$, we  parametrize the unipotent cell by $(t_1, t_2, \ldots, t_a) \mapsto y_{i_1}(t_1) y_{i_2}(t_2) \cdots y_{i_a}(t_a)$.
 Given any bounded affine permutation and any plabic graph $G$, we parametrize the  positroid cell by the boundary measurement map.
\item The parametrizations are bijective if and only if the word/graph is reduced.
\end{itemize}

To make the analogy stronger, we now describe the analogues of the presentation of the symmetric group by generators and relations.
Postnikov constructs a collection of \newword{transformation} and \newword{reduction} moves, which transform one plabic graph into another and transform a collection of weights on the first graph into a collection of weights on the second graph, in such a way as to preserve the boundary measurements. These maps are well defined for positive real weights, but are only rational maps when considered for more general weights.
Transformation moves preserve the dimension of $T$; reduction moves reduce it.
See~\cite[Section 12]{Post} for details.

\begin{theorem}
Let $\Pio$ be a positroid cell, let $G$ be a reduced plabic graph for $\Pio$ and let $G'$ be any plabic graph for $\Pio$. Then there is a series of reduction and transformation moves turning $G'$ into $G$.
\end{theorem}

We note that we never need to apply reduction moves in the reverse direction and that, if $G$ and $G'$ are both reduced, then $\dim T = \dim T'$ so they are joined purely by transformation moves.
This is analogous to the theorem that, when transforming a word into a reduced word, one uses the substitution $s_i s_i \leadsto 1$, together with braid and commutation relations, but not $1 \leadsto s_i s_i$. 

\begin{remark} \label{Padding}
We said before that plabic graphs should be thought of as a generalization of reduced words in the symmetric group, and the plabic parametrization of positroid varieties as a generalization of the Chevalley parametrization of unipotent cells.
We can now be more precise: Let $w$ be a permutation in $S_m$
Define the bounded affine permutation $\widetilde{w}$ in $\tS^m_{2m}$ by
\[ \widetilde{w}(i) = \begin{cases} 2m+1-i & 1 \leq i \leq m \\ w^{-1}(i-m)+2m & m+1 \leq i \leq 2m \end{cases} . \]
The unipotent cell for $w$ is closely related to the positroid cell for $\widetilde{w}$.
If $s_{i_1} s_{i_2} \cdots s_{i_a}$ is any reduced word for $w$, then we can make a reduced plabic graph for $\widetilde{w}$ as follows: Take the disc $D$ to be a rectangle, put boundary vertices $1$, $2$, \dots, $m$ on the left hand side of the rectangle and vertices $m+1$, $m+2$, \dots, $2m$ on the right, and draw horizontal lines joining $i$ to $2m+1-i$ for $1 \leq i \leq m$. Then, for each letter $s_{i_k}$ in the word, draw a line segment between the horizontal line at height $i_k$ and the one at height $i_k +1$, making the top of the line segment white and the bottom black. Finally, if vertices of the same color wind up neighboring each other along one of the horizontal segments, then contract that same color segment down to a single point.
Then, choosing the correct coordinates, Postnikov's parametrization of the positroid cell for $\widetilde{w}$ using this plabic graph corresponds to the parametrization of the corresponding unipotent cell as a product of Chevalley generators.
See~\cite[Section 18]{Post} and~\cite[Section 12]{OPS} for more.
\end{remark}

\begin{remark}
More generally, let $u \preceq_k w$ and choose a reduced word $s_{i_1} s_{i_2} \cdots s_{i_a}$ for $w$. Deodhar parametrizes corresponding positroid cell.
Karpman~\cite{Karpman} tells us how to convert the data of $(u, w, (i_1, i_2, \ldots, i_a))$ into a reduced plabic graph (the ``bridge diagram") which recovers the Deodhar parametrization.
\end{remark}

\subsection{Zig-zag paths}

Suppose that one is given a plabic graph $G$ and one wants to know whether $G$ is reduced.
Let $\Pio$ be the corresponding positroid variety, with affine permutation $f$. 
By Theorem/Definition~\ref{ReducedMotivation}, $G$ is reduced if and only if the number of faces of $G$ is $\dim \Pio+1 = k(n-k)-\ell(f)+1$.
The number of faces of $G$, and the length $\ell(f)$, are both combinatorial.
However, at the moment, we do not have a combinatorial rule for determining $f$ from $G$.
In this section, we will partially address this gap and give a combinatorial criterion for when $G$ is reduced.
To see why we write that the issue is only partially addressed, see Problem~\ref{DemazureZigZag}.

Let $G$ be a plabic graph. 
For simplicity, we assume throughout this section that all of the boundary vertices of $G$ are degree $1$. 
We can always reduce to this case as follows: If $v$ is a boundary vertex, then push $v$ into the interior of the disc $D$ and add a path $(v'', v', v)$, making $v''$ into a new boundary vertex. The edges $(v', v)$ and $(v'', v')$ are given weight $1$.

A \newword{zig-zag path} is a directed path that travels along the edges of $G$ such that, when ever the path comes to a white vertex, it turns as far left as possible and, whenever the path comes to a black vertex, it turns as far right as possible. 
We continue a zig-zag path in both directions until either (1) it forms a closed loop or (2) it runs into a boundary vertex in such a way that the resulting right turn would carry it out of the disc.

\begin{eg} \label{eg:ZZHexagon}
We have drawn two of the six zig-zag paths for Example~\ref{eg:G36Hexagon}: 

\centerline{
\begin{tikzpicture}
\draw[fill=none] (0,0) circle (2.0) ; 
\node [label={[label distance=0.2cm]180:$1$}] (1) at (-2,0) {};
\node  [label={[label distance=0.2cm]120:$2$}] (2) at (-1, 1.732) {};
\node  [label={[label distance=0.2cm]60:$3$}] (3) at (1, 1.732) {};
\node  [label={[label distance=0.2cm]0:$4$}] (4) at (2,0) {};
\node  [label={[label distance=0.2cm]300:$5$}] (5) at (1, -1.732) {};
\node   [label={[label distance=0.2cm]240:$6$}] (6) at (-1, -1.732) {};
\node (p) at (0.5,0.866) {};
\node (q) at (0.5,-0.866) {};
\node (r) at (-1,0) {};
\draw (2.center)--(p.center) node [midway,above] {$x_2$};  \draw (4.center)--(p.center) node [midway,above] {$x_3$};
\draw (4.center)--(q.center) node [midway,below] {$x_4$}; \draw (5.center)--(q.center); \draw (6.center)--(q.center) node [midway,below] {$x_5$};
\draw (6.center)--(r.center) node [midway,left] {$x_6$}; \draw (2.center)--(r.center) node [midway,left] {$x_1$};
\draw[fill=white] (p) circle (0.1) ;  \draw[fill=white] (q) circle (0.1) ; \draw[fill=white] (r) circle (0.1) ; 
\draw[fill=black] (1) circle (0.1) ;  \draw[fill=black] (2) circle (0.1) ; \draw[fill=black] (3) circle (0.1) ; \draw[fill=black] (4) circle (0.1) ;  \draw[fill=black] (5) circle (0.1) ; \draw[fill=black] (6) circle (0.1) ; 
 \draw[very thick, -{Latex[length=4mm,width=3mm]}] (1.center)--(r.center);  \draw[very thick, -{Latex[length=4mm,width=3mm]}] (r.center)--(2.center); \draw[very thick, -{Latex[length=4mm,width=3mm]}] (2.center)--(p.center); \draw[very thick, -{Latex[length=4mm,width=3mm]}] (p.center)--(3.center);
 \draw[very thick, -{Latex[length=4mm,width=3mm]}] (6.center)--(q.center); \draw[very thick, -{Latex[length=4mm,width=3mm]}] (q.center)--(4.center);
\end{tikzpicture}}
One path goes from $1$ to $3$ and the other goes from $6$ to $2$. 
 The reader can convince themself that, in general, if $i$ if odd, then the zig-zag path starting at $i$ ends at $i+2 \bmod 6$ and, if $i$ is even, then the path starting at $i$ ends at $i+4 \bmod 6$, consistent with the bounded affine permutation for the corresponding positroid cell. 
\end{eg}

Postnikov~\cite[Theorem 13.2]{Post} used the zig zag paths to give a criterion for a plabic graph $G$ to be reduced. 
For simplicity, we assume that $G$ does not have any vertices of degree $1$, except those directly incident to boundary vertices; we can always reduce to this case because if $G$ is a graph with a degree $1$ vertex $v$, incident to another interior vertex $w$, then we can delete $v$ and $w$ from $G$ to make a modified plabic graph $G'$ whose matching are in obvious bijection with the matchings of $G$.

\begin{Theorem} \label{ReducedCriterion}
Let $G$ be a plabic graph which has no degree $1$ vertices incident to boundary vertices. Then $G$ is reduced if and only if
\begin{enumerate}
\item Every zig-zag path of $G$ joins two points in $\partial(G)$; there are no closed loops.
\item If $\alpha$ is a zig-zag path joining boundary vertex $i$ to boundary vertex $j$ for $j \neq i$, then $\alpha$ does not pass through any edge twice.
\item If $\alpha$ and $\beta$ are two distinct  zig-zag paths, and $e$ and $f$ are edges occurring in $\alpha$ and in $\beta$, then $e$ and $f$ occur in opposite orders in $\alpha$ and in $\beta$.
\item If $\alpha$ is a zig-zag path joining boundary vertex $i$ to itself, then either vertex $i$ is an isolated vertex, or else borders a single interior vertex.
\end{enumerate}
Moreover, if these conditions hold, then the bounded affine permutation $f$ can be obtained as follows: The zig-zag path starting at $i$ ends at $f(i) \bmod n$. If the zig-zag path starting at $i$ also ends at $i$, then $f(i)=i$ if $i$ is an isolated vertex and $f(i) = i+n$ if $i$ has an internal neighbor.
\end{Theorem}

\begin{remark} \label{BoundaryConvention}
We have stuck to our convention that vertices in $\partial(G)$ are black. If we allowed white vertices and black vertices, we could reduce to the case that $G$ has no vertices of degree $1$ at all, and the rule would be that,  if the zig-zag path starting at $i$ also ends at $i$, then $f(i)=i$ if $i$ is black and $f(i) = i+n$ if $i$ is white.
\end{remark}

We have written the theorem to match the presentation of Postnikov's Theorem 13.2 as closely as possible.  Here is an alternate presentation which combines conditions (1), (2) and (3):
Let $G$ be a plabic graph with no vertices of degree $1$; we can allow boundary vertices of both colors. Let $e$ be an edge of $G$, and let $\gamma_1$ and $\gamma_2$ be the zig-zag paths which start passing through $e$ in the two possible directions and travel forward indefinitely until we hit the boundary of $G$, or infinitely if we don't. 
\begin{theorem} \label{ZZReduced}
With the above notation and conventions, $G$ is reduced if and only if, $\gamma_1(e)$ and $\gamma_2(e)$ have no edges in common except for the initial edge $e$.
\end{theorem}
This presentation is similar to that of~\cite[Theorem 5.5]{Bocklandt}.

\begin{remark}
The terminology ``zig-zag path" originated in Kenyon's lectures on the dimer model~\cite{Kenyon}, where it is credited to joint work with Schlenker.
It was then taken up in the quiver representation literature, where it can be found in work of Hanany and Vegh~\cite{HananyVegh}, Broomhead~\cite{Broomhead}, Mozkovoy and Reineike~\cite{MozRein} and Bocklandt~\cite{Bocklandt}. 
In particular, Bocklandt showed that the condition analogous to Theorem~\ref{ZZReduced} (for graphs on a torus, rather than a disc), is equivalent to many algebraic conditions on the path algebra of the quiver with potential. 

Postnikov's own chosen term is ``trip", and he calls the permutation $f$ the ``trip permutation".
\end{remark}

\begin{remark}
A close relative of the zig-zag path is the ``alternating strand diagram". To draw a strand, take a zig-zag path and perturb it to make an oriented curve $\gamma$ such that the white vertices lie to the right of $\gamma$ and the interior black vertices lie to the left of $\gamma$.
As an example, we draw the alternating strand diagram for  Example~\ref{eg:G36Hexagon}: 

\centerline{
\begin{tikzpicture}
\draw[fill=none] (0,0) circle (2.0) ; 
\node [label={[label distance=0.2cm]180:$1$}] (1) at (-2,0) {};
\node  [label={[label distance=0.2cm]120:$2$}] (2) at (-1, 1.732) {};
\node  [label={[label distance=0.2cm]60:$3$}] (3) at (1, 1.732) {};
\node  [label={[label distance=0.2cm]0:$4$}] (4) at (2,0) {};
\node  [label={[label distance=0.2cm]300:$5$}] (5) at (1, -1.732) {};
\node   [label={[label distance=0.2cm]240:$6$}] (6) at (-1, -1.732) {};
\node (p) at (0.5,0.866) {};
\node (q) at (0.5,-0.866) {};
\node (r) at (-1,0) {};
\draw (2.center)--(p.center); \draw (3.center)--(p.center); \draw (4.center)--(p.center);
\draw (4.center)--(q.center); \draw (5.center)--(q.center); \draw (6.center)--(q.center);
\draw (6.center)--(r.center); \draw (1.center)--(r.center); \draw (2.center)--(r.center);
\draw[fill=white] (p) circle (0.1) ;  \draw[fill=white] (q) circle (0.1) ; \draw[fill=white] (r) circle (0.1) ; 
\draw[fill=black] (1) circle (0.1) ;  \draw[fill=black] (2) circle (0.1) ; \draw[fill=black] (3) circle (0.1) ; \draw[fill=black] (4) circle (0.1) ;  \draw[fill=black] (5) circle (0.1) ; \draw[fill=black] (6) circle (0.1) ; 
\draw[dashed,  -{Latex[length=3mm,width=2.25mm]}] plot [smooth, tension=1] coordinates { (1) (-1,0.866) (-0.25,1.299) (3)};
\draw[dashed,  -{Latex[length=3mm,width=2.25mm]}] plot [smooth, tension=1] coordinates { (3) (1.25,0.433) (1.25,-0.433) (5)};
\draw[dashed,  -{Latex[length=3mm,width=2.25mm]}] plot [smooth, tension=1] coordinates { (5) (-0.25,-1.299) (-1,-0.866) (1)};
\draw[dashed,  -{Latex[length=3mm,width=2.25mm]}] plot [smooth, tension=1] coordinates { (4) (1.25,0.433) (0.25, 0.433) (-0.25,1.299)  (2)};
\draw[dashed,  -{Latex[length=3mm,width=2.25mm]}] plot [smooth, tension=1] coordinates { (2) (-1,0.866) (-0.5, 0) (-1,-0.866)  (6)};
\draw[dashed,  -{Latex[length=3mm,width=2.25mm]}] plot [smooth, tension=1] coordinates { (6) (-0.25,-1.299) (0.25, -0.433)  (1.25, -0.433)  (4)};
\end{tikzpicture}}

 So each edge $e$ has two strands which cross at the midpoint of $e$; each white vertex has a collection of strands circling it clockwise and each interior black vertex has a collection of strands circling it counter-clockwise. 
Given a collection of oriented curves in the disc, the collection comes from a plabic graph if and only if, as you travel along each strand $\gamma$, the strands crossing $\gamma$ alternately come in from the left and from the right.
Postnikov defines such a collection of strands to be an ``alternating strand diagram". The plabic graph can be recovered from the alternating strand diagram: The white vertices correspond to the regions which are circled in a clockwise direction from the strands and the interior black vertices correspond to the regions. See~\cite[Section 14]{Post} for details.
\end{remark}

\begin{remark}
For the reduced plabic graph described in Remark~\ref{Padding}, the zig-zag paths starting at $i$ for $1 \leq i \leq m$ go horizontally directly across to $2m+1-i$, and the zig-zag paths starting at $i$ for $m+1 \leq i \leq 2m$ form a wiring diagram for the word $s_{i_1} s_{i_2} \cdots s_{i_a}$.
\end{remark}

We close with a problem that the author considers an embarrassing gap in the field.
If $G$ is a reduced plabic graph, then Theorem~\ref{ReducedCriterion} tells us how to (1) detect that $G$ is reduced and (2) find the bounded affine permutation corresponding to $G$. 
If $G$ is a non-reduced plabic graph, then Theorem~\ref{ReducedCriterion} will detect that $G$ is not reduced but will \text{not} give us a rule to extract the bounded affine permutation.
\begin{problem} \label{DemazureZigZag}
If $G$ is a non-reduced plabic graph, giving a parametrization of the positroid cell $\Pio_{>0}$, give a combinatorial rule to extract the bounded permutation of $\Pio$ from the zig-zag paths of $G$.
\end{problem}
The author posed this problem on Mathoverflow~\cite{MO375807} and, as of January 2024, he has received no answers.

\subsection{The twist and its consequences} \label{sec:twist}

For any plabic graph $G$, with corresponding positroid variety $\Pio$, we have described a boundary measurement map $\mu : T_{>0} \to \Pio_{>0}$. 
One would like to, more generally, define a map of complex varieties $\mu: T \to \Pio$, where $T$ is the complex torus of edge weights modulo gauge transformation, and $\Pio$ is the complex open positroid variety.
The awkward point is that, once edge weights can be complex (or even just negative), there can be cancellation in the sums defining the boundary measurement map, so the Pl\"ucker coordinate $\Delta_I$ may be given by a nonempty sum, and hence not identically zero on $\Pio$, but may still be zero somewhere on $T$.
Indeed, in example~\ref{eg:BasicG24GaugeFixed}, we saw the boundary measurement map where $(p,q,r,s) \in \GG_m^4$ is mapped to $(\Delta_{12}, \Delta_{13}, \Delta_{14}, \Delta_{23}, \Delta_{24}, \Delta_{34}) = (s, pr+qs, p, r, 1, q)$.
So we can have $\Delta_{13} = pr+qs$ equal to $0$ even though $p$, $q$, $r$ and $s$ are nonzero.
It is therefore unclear that $\mu$ will land in $\Pio$, and it is not even clear that $\mu$ is well-defined as a map to $G(k,n)$: Maybe there are some points of $T$ where all of the Pl\"ucker cordinates vanish.

\begin{eg}
Indeed, this can happen quite easily for non-reduced $G$. 
The graph below is a non-reduced graph for the big cell in $G(1,2) = \PP^1$.

\centerline{
\begin{tikzpicture}
\draw[fill=none] (0,0) circle (1.5) ; 
\node [label={[label distance=0.2cm]180:$1$}] (1) at (-1.5,0) {};
\node  [label={[label distance=0.2cm]0:$2$}] (2) at (1.5,0) {};
\node (a) at (0,0) {};
\draw (1.center) to [bend right]  node [midway, below] {$q$} (a.center); \draw (1.center) to [bend left] node [midway, above] {$p$} (a.center);
\draw (a.center) to [bend right ] node [midway, below] {$s$}  (2.center)  ; \draw (a.center) to [bend left] node [midway, above] {$r$} (2.center)  ;
\draw[fill=black] (1) circle (0.1);
\draw[fill=white] (a) circle (0.1);
\draw[fill=black] (2) circle (0.1);
\end{tikzpicture}}

In terms of the homogeneous coordinates on $\PP^1$, the boundary measurement map is $(p,q,r,s) \mapsto (p+q : r+s)$.
If $p+q=r+s=0$, then this does not give a well defined point of $\PP^1$.
\end{eg}

Surprisingly, when $G$ is reduced, these problems do not occur! 
The following is the main result of Muller and Speyer's paper~\cite{MullerSpeyerTwist}: 
\begin{theorem} \label{TwistMain}
Let $G$ be a reduced plabic graph with corresponding positroid variety. Then the boundary measurement map $\mu : T \to \Pio$ is an open inclusion. 
\end{theorem}

\begin{eg}
It is \textbf{not} always true that the image of the boundary measurement map can be described by the non-vanishing of certain Pl\"ucker coordinates. 
The following example is taken from~\cite[Appendix~A.3]{MullerSpeyerTwist}. Consider the reduced plabic graph shown below, corresponding to the big positroid cell in $G(3,6)$:

\centerline{
\begin{tikzpicture}
\draw[fill=none] (0,0) circle (2.0) ; 
\node [label={[label distance=0.2cm]0:$1$}] (1) at (2,0) {};
\node [label={[label distance=0.2cm]300:$2$}] (2) at (1,-1.732) {};
\node [label={[label distance=0.2cm]240:$3$}] (3) at (-1,-1.732) {};
\node [label={[label distance=0.2cm]180:$4$}] (4) at (-2,0) {};
\node [label={[label distance=0.2cm]120:$5$}] (5) at (-1,1.732) {};
\node [label={[label distance=0.2cm]60:$6$}] (6) at (1,1.732) {};
\node (a) at (1,0) {};
\node (b) at (0.5,-0.866) {};
\node (c) at (-0.5,-0.866) {};
\node (d) at (-1,0) {};
\node (e) at (-0.5,0.866) {};
\node (f) at (0.5,0.866) {};
\node (b2) at (0.75,-1.299) {};
\node (d2) at (-1.5,0) {};
\node (f2) at (0.75,1.299) {};
\draw (a.center) -- (b.center) -- (c.center) -- (d.center) -- (e.center) -- (f.center) -- (a.center) ;
\draw (a.center) -- (1.center); \draw (c.center) -- (3.center); \draw (e.center) -- (5.center);
\draw (b.center) -- (b2.center) -- (2.center); \draw (d.center) -- (d2.center) -- (4.center);  \draw (f.center) -- (f2.center) -- (6.center);
\draw (1.center) -- (f2.center); \draw (3.center) -- (b2.center); \draw (5.center) -- (d2.center); 
\draw[fill=black] (1) circle (0.1); \draw[fill=black] (2) circle (0.1); \draw[fill=black] (3) circle (0.1); \draw[fill=black] (4) circle (0.1); \draw[fill=black] (5) circle (0.1); \draw[fill=black] (6) circle (0.1);
\draw[fill=black] (b) circle (0.1); \draw[fill=black] (d) circle (0.1); \draw[fill=black] (f) circle (0.1); 
\draw[fill=white] (b2) circle (0.1); \draw[fill=white] (d2) circle (0.1); \draw[fill=white] (f2) circle (0.1); 
\draw[fill=white] (a) circle (0.1); \draw[fill=white] (c) circle (0.1); \draw[fill=white] (e) circle (0.1); 
\end{tikzpicture}}

The image of the boundary measurement is the locus 
\[ \Delta_{123} \Delta_{234} \Delta_{345} \Delta_{456} \Delta_{156} \Delta_{126} \Delta_{125} \Delta_{134} \Delta_{356} ( \Delta_{124} \Delta_{356} - \Delta_{123} \Delta_{456}) \neq 0. \]
The final binomial is not a product of Pl\"ucker coordinates
\end{eg}

As one might guess, one proves Theorem~\ref{TwistMain} by constructing a rational inverse to $\mu$. 
This inverse requires two constructions: The \newword{twist map} and the \newword{face labeling of $G$}.
We will describe both constructions but leave the details to~\cite{MullerSpeyerTwist}.

The \newword{twist map} is an automorphism $\tau$ of $\Pio$.
The \newword{face labeling} is a set of Pl\"ucker coordinates $\cI(G)$, one for each face of $G$.
Let $\Omega$ be the open subset of $\Pio$ where the coordinates in $\cI(G)$ are nonzero. Let $\sigma : \Omega \to \GG_m^{\#\cI(G)-1}$ be the map taking a point of $\Omega$ to those Pl\"ucker coordinates indexed by the face labels.
(The $-1$ is because the Pl\"ucker coordinates are homogeneous coordinates.)
Muller and Speyer show:
\begin{theorem}
The twist of the image of the boundary measurement map is $\Omega$; in other words, $\tau(\mu(T)) = \Omega$. The maps $T \overset{\mu}{\longrightarrow} \mu(T) \overset{\tau}{\longrightarrow} \Omega \overset{\sigma}{\longrightarrow} \GG_m^{\#\cI(G)-1}$ are all isomorphisms.
\end{theorem}
For the big positroid cell, these results were found earlier by Marsh and Scott~\cite{MarshScott}.
We now describe the twist and the face labeling.

Let $M$ be an $n \times k$ matrix of rank $k$, with rows $\vec{v}_1$, $\vec{v}_2$, \ldots, $\vec{v}_n$.
The column span of $M$ gives a point of the Grassmannian $G(k,n)$; let it lie in the positroid variety with decorated permutation $f$. For simplicity, we assume that $f$ has no fixed points, which means that none of the $\vec{v}_i$ are zero.
Associated to $f$ is the Grassmann necklace $(I_1, I_2, \ldots, I_n)$ and, using our hypothesis that $f(i) \neq i$, we have $i \in I_i$ for each $i$.
We define $\vec{w}_i$ to be the unique vector such that
\[ \vec{w}_i \cdot \vec{v}_j = \begin{cases} 1 &i=j \\ 0 & j \in I_i \setminus \{ i \}. \end{cases}. \]
We define $\tau(M)$ to be the matrix with rows $\vec{w}_1$, $\vec{w}_2$, \dots, $\vec{w}_n$.

\begin{eg} \label{eg:TwistG36Hexagon}
We work with the matrix from Example~\ref{eg:G36Hexagon}.
\[ \begin{bmatrix} 
1&0&0\\
x_1&x_2&0 \\
0&1&0 \\
0&x_3&x_4 \\
0&0&1\\
x_6&0&x_5\\
\end{bmatrix} .\]
 We have $I_1 = \{ 1,2,4 \}$, so $\vec{w}_1$ is defined by the equations $\vec{w}_1 \cdot \begin{sbm} x_1&x_2&0 \\ \end{sbm} = \vec{w}_1 \cdot \begin{sbm} 0&x_3&x_4 \\ \end{sbm} = 0$ and $\vec{w}_1 \cdot \begin{sbm} 1&0&0 \end{sbm}=1$, giving
 $\vec{w}_1 = \begin{sbm} 1&\ \ -\tfrac{x_1}{x_2}&\ \  \tfrac{x_1 x_3}{x_2 x_4} \\ \end{sbm}$. 
 Similarly, $I_2 = \{ 2,3,4 \}$, so $\vec{w}_2$ is defined by the equations $\vec{w}_2 \cdot \begin{sbm} 0&1&0 \\ \end{sbm} = \vec{w}_1 \cdot \begin{sbm} 0&x_3&x_4 \\ \end{sbm} = 0$ and $\vec{w}_2 \cdot \begin{sbm} x_1&x_2&x_3 \end{sbm}=1$, giving $\vec{w}_2 = \begin{sbm} \tfrac{1}{x_1} & 0 & 0 \end{sbm}$. Continuing in this manner, we compute
 \[ \tau \left(  \begin{bmatrix} 
1&0&0\\
x_1&x_2&0 \\
0&1&0 \\
0&x_3&x_4 \\
0&0&1\\
x_6&0&x_5\\
\end{bmatrix} \right) =
\begin{bmatrix} 
1 & -\tfrac{x_1}{x_2} & \tfrac{x_1 x_3}{x_2 x_4} \\[0.1 cm]
\tfrac{1}{x_1} &0&0 \\[0.1 cm]
\tfrac{x_3 x_5}{x_4 x_6} &1 & -\tfrac{x_3}{x_4} \\[0.1 cm]
0&\tfrac{1}{x_3} & 0 \\[0.1 cm]
-\tfrac{x_5}{x_6} &\tfrac{x_1 x_5}{x_2 x_6} & 1\\[0.1 cm]
0&0&\tfrac{1}{x_5} \\
\end{bmatrix}.\]
\end{eg}

We defined $M \to \tau(M)$ as a map from $k \times n$ matrices to $k \times n$ matrices, but it descends to an map $G(k,n) \to G(k,n)$. 
To see why, notice that $M$ and $M'$ have the same row span if and only if $M' = Mg$ for some $g \in \GL_k$. 
Then note that $\tau(Mg) = \tau(M) (g^T)^{-1} $. 
On each open positroid stratum $\Pio$ of $G(k,n)$, the twist $\tau : \Pio \to G(k,n)$ is an algebraic morphism. 
Though not obvious, $\tau$ maps each positroid cell to itself.

\begin{eg}
In Example~\ref{eg:TwistG36Hexagon}, we worked with the positroid cell where $\Delta_{123} = \Delta_{345} = \Delta_{156} = 0$. The reader may check that these minors also vanish in the twisted matrix that we computed.
\end{eg}

\begin{TD}
For each positroid cell $\Pio$, the twist map $\tau$ descends to an algebraic automorphism of $\Pio$, also called the \newword{twist} and denoted $\tau$. 
\end{TD}

\begin{remark}
The inverse of the twist is constructed similarly, using the reverse Grassmann necklace.
\end{remark}

We now describe the second ingredient of our construction, the face labeling. 	
Let $G$ be a reduced plabic graph for a positroid cell in $G(k,n)$, with corresponding bounded permutation $f$.
For $i \in [n]$, let $\gamma_i$ be the zig-zag path from $f^{-1}(i)$ to $i$. 
Since $G$ is reduced, each $\gamma_i$ is a single path separating the disc $D$ into two connected components.
For $F$ a face of $G$, let $I(F)$ be the set of $i$ such that $\gamma_i$ lies to the left hand side of $\gamma_i$. 
(If $f(i) = i$, then $i$ does not lie in any $I(F)$; if $f(i) = i+n$, then $i$ lies in every $I(F)$.)

\begin{eg} \label{eg:FaceLabelsHexagon}
In Example~\ref{eg:ZZHexagon}, we drew two of the zig-zag paths for the plabic graph shown below. Here we redraw those two zig-zag paths, $\gamma_3$ and $\gamma_4$, and label the faces of the graph.

\centerline{
\begin{tikzpicture}
\draw[fill=none] (0,0) circle (2.0) ; 
\node [label={[label distance=0.2cm]180:$1$}] (1) at (-2,0) {};
\node  [label={[label distance=0.2cm]120:$2$}] (2) at (-1, 1.732) {};
\node  [label={[label distance=0.2cm]60:$3$}] (3) at (1, 1.732) {};
\node  [label={[label distance=0.2cm]0:$4$}] (4) at (2,0) {};
\node  [label={[label distance=0.2cm]300:$5$}] (5) at (1, -1.732) {};
\node   [label={[label distance=0.2cm]240:$6$}] (6) at (-1, -1.732) {};
\node (p) at (0.5,0.866) {};
\node (q) at (0.5,-0.866) {};
\node (r) at (-1,0) {};
\draw (2.center)--(p.center);  \draw (4.center)--(p.center);
\draw (4.center)--(q.center); \draw (5.center)--(q.center); \draw (6.center)--(q.center);
\draw (6.center)--(r.center); \draw (2.center)--(r.center);
\draw[fill=white] (p) circle (0.1) ;  \draw[fill=white] (q) circle (0.1) ; \draw[fill=white] (r) circle (0.1) ; 
\draw[fill=black] (1) circle (0.1) ;  \draw[fill=black] (2) circle (0.1) ; \draw[fill=black] (3) circle (0.1) ; \draw[fill=black] (4) circle (0.1) ;  \draw[fill=black] (5) circle (0.1) ; \draw[fill=black] (6) circle (0.1) ; 
 \draw[very thick, -{Latex[length=4mm,width=3mm]}] (1.center)--(r.center);  \draw[very thick, -{Latex[length=4mm,width=3mm]}] (r.center)--(2.center); \draw[very thick, -{Latex[length=4mm,width=3mm]}] (2.center)--(p.center); \draw[very thick, -{Latex[length=4mm,width=3mm]}] (p.center)--(3.center);
 \draw[very thick, -{Latex[length=4mm,width=3mm]}] (6.center)--(q.center); \draw[very thick, -{Latex[length=4mm,width=3mm]}] (q.center)--(4.center);
 \node (124) at (-1.333,-0.577) {$124$};
  \node (234) at (-1.333, 0.577) {$234$};
   \node (346) at (0.167,1.433) {$346$};
      \node (456) at (1.167,0.866) {$456$};
            \node (256) at (1.167,-0.866) {$256$};
               \node (126) at (0.167,-1.433) {$126$};
               \node (246) at (0,0) {$246$};
\end{tikzpicture}}

\end{eg}

\begin{remark}
The face which borders the portion of $\partial(D)$ between boundary vertex $i-1$ and boundary vertex $i$ will receive the label $I_i$, the $i$-th element of the Grassmann necklace. 
\end{remark}

\begin{remark}
We have assigned the path from $f^{-1}(i)$ to $i$ the label $i$, so-called ``target labeling". We could also assign the label $i$ to the path from $i$ to $f(i)$, so-called ``source labeling". 
Both source and target labelings are important, see~\cite{MullerSpeyerTwist} for a presentation that includes both and see~\cite{FraserSB} for relations between them.
\end{remark}

\begin{eg}
We work out Muller and Speyer's result for our running example. We want to analyze the composite $T \overset{\mu}{\longrightarrow} \mu(T) \overset{\tau}{\longrightarrow} \Omega \overset{\sigma}{\longrightarrow} \GG_m^{\#\cI(G)-1}$. The map $\mu$ is computed in Example~\ref{eg:G36Hexagon} and the map $\tau$ is computed in Example~\ref{eg:TwistG36Hexagon}. 
The map $\sigma$ is given by evaluation of the face minors from Example~\ref{eg:FaceLabelsHexagon}, namely, $(\Delta_{124} : \Delta_{234} : \Delta_{346} : \Delta_{456} : \Delta_{256} : \Delta_{126} : \Delta_{246})$. Computing the corresponding minors of the matrix in  Example~\ref{eg:TwistG36Hexagon}, we get
\[ (\tfrac{1}{x_2 x_4} : \tfrac{1}{x_1 x_4} : \tfrac{1}{x_4 x_6} : \tfrac{1}{x_3 x_6} : \tfrac{1}{x_2 x_4} : \tfrac{1}{x_2 x_5} : \tfrac{1}{x_1 x_3 x_5} ) . \]
As the reader can easily see, this is a monomial map from one $6$-dimensional torus to another. (Recall that the $\Delta$'s are homogeneous coordinates, so only defined up to common ratio.) 
A bit more work checks that this monomial map is invertible, which is the result. (For example, $x_1 = \tfrac{1/(x_2 x_5) \cdot 1/(x_3 x_6)}{1/(x_2 x_6) \cdot 1/(x_1 x_3 x_5)}$.) 
For an explicit formula for the inverse, and much more, see~\cite{MullerSpeyerTwist}.
\end{eg}

\begin{remark}
It would be interesting to investigate the twist map $\tau : \Pio \to \Pio$ from the perspective of complex dynamical systems.
\end{remark}


\begin{remark}
Lam and Galashin~\cite{GLCluster}, building on previous work of~\cite{Scott, MullerSpeyerTwist, SSbW} construct a cluster structure on $\Pio$, for which the Pl\"ucker coordinates $\Delta_I$, as $I$ runs over the set of face labels, are the cluster variables.
\end{remark}

\section{Acknowledgments}
I have learned these topics through many conversations with many mathematicians. 
I recall learning about Gr\"obner bases, toric varieties and degeneration from Bernd Sturmfels and Ezra Miller, learning about Frobenius splitting and Lie theory from Allen Knutson, learning about total positivity and positroids from Lauren Williams and learning about perfect matchings from Jim Propp. 
I have also learned from my students: Giwan Kim taught me how to think about standard monomial theory, and Gracie Ingermanson taught me how to think about Deodhar pieces.
I appreciate comments on the manuscript from Dave Anderson, Terrence George, Daoji Huang, Steven Karp,  Allen Knutson, Thomas Lam, George Lusztig, Oliver Pechenik, Minh-Tam Trinh and Lauren Williams.

I am grateful to my editors for their patience with the lengthy process of writing this manuscript, and to my referees for their extremely careful reading.

I was supported by NSF grants DMS-1854225 and DMS-1855135 while writing this manuscript.


\end{document}